\documentstyle[11pt,amsfonts]{article}

\newtheorem{thm}{Theorem}[section]
\newtheorem{propose}[thm]{Proposition}
\newtheorem{lemma}[thm]{Lemma}
\newtheorem{cor}[thm]{Corollary}
\newtheorem{rmk}[thm]{Remark}
\newtheorem{defn}[thm]{Definition}

\renewcommand{\P}{\mbox{$\Bbb P$}}   
\newcommand{\V}{\mbox{$\Bbb V$}}     
\newcommand{\A}{\mbox{$\Bbb A$}}     
\newcommand{\Q}{\mbox{$\Bbb Q$}}     
\newcommand{\C}{\mbox{$\Bbb C$}}     
\newcommand{\Z}{\mbox{$\Bbb Z$}}     
\newcommand{\HH}{\mbox{$\Bbb H$}}
\newcommand{\ssp}[1]{\mbox{$\scriptscriptstyle {#1}$}}

\newcommand{\Gal}{{\rm Gal}\,}

\newcommand{\Hom}{{\rm Hom}\,}      
\newcommand{\RHom}{{\rm RHom}\,}    
\newcommand{\RG}{{\rm R}\Gamma\,}   

\newcommand{\im}{{\rm im}\,}        
\newcommand{\coker}{{\rm coker}\,}  
\newcommand{\gr}{{\rm gr}\,}        
\newcommand{\dlog}{{\rm dlog}\,}    

\newcommand{\Pic}{{\rm Pic}\,}      
\newcommand{\bPic}{\mbox{$\bf Pic$}} 
\newcommand{\Lie}{{\rm Lie}\,}      
\newcommand{\Spec}{{\rm Spec}\,}    
\newcommand{\Ext}{{\rm Ext}\,}      
\newcommand{\rank}{{\rm rank}\,}    

\newcommand{\Cl}{{\rm Cl}\,}
\newcommand{\G}{\mbox{$\Bbb G$}}    
\newcommand{\Alb}{{\rm Alb}\,}
\renewcommand{\exp}{{\rm exp}\,}

\newcommand{\nin}{\not\in}
\newcommand{\by}[1]{\stackrel{#1}{\rightarrow}}
\newcommand{\longby}[1]{\stackrel{#1}{\longrightarrow}}

\newcommand{\mod}[1]{\mbox{$\displaystyle\mid #1 \mid$}}
\newcommand{\implies}{\mbox{$\Rightarrow$}}
\newcommand{\tensor}{\otimes}
\newcommand{\into}{\hookrightarrow}

\renewcommand{\tilde}{\widetilde}
\renewcommand{\hat}{\widehat}

\renewcommand{\leadsto}{\mbox{ $\longmapsto$ }}
\newcommand{\df}{\mbox{\,$\stackrel{\ssp{\rm def}}{=}$}\,}
\newcommand{\ie}{{\it i.e.\/},\ }
\newcommand{\cf}{{\it cf.\/}\ }
\newcommand{\eg}{{\it e.g.\/},\  }

\newcommand{\et}{\mbox{\scriptsize{\'{e}t}}}

\newcommand{\s}[1]{\mbox{\scriptsize{$#1$}}}

\renewcommand{\bar}{\overline}
\newcommand{\srhom}{\mbox{${\cal RH}om\,$}}

\newcommand{\veq}{\mbox{\large $\parallel$}}  
\newcommand{\sZ}{\mbox{\scriptsize{$\Z$}}}   
\newcommand{\sC}{\mbox{\scriptsize{$\C$}}}   
\newcommand{\sQ}{\mbox{\scriptsize{$\Q$}}}   
\newcommand{\longto}{\longrightarrow}

\newcommand{\limdir}[1]{{\displaystyle{\mathop{\rm
lim}_{\buildrel\longrightarrow\over{#1}}}}\,}
\newcommand{\liminv}[1]{{\displaystyle{\mathop{\rm
lim}_{\buildrel\longleftarrow\over{#1}}}}\,}
\newcommand{\onto}{\mbox{$\to\!\!\!\to$}}

\newcommand{\boxtensor}{{\Box\kern-9.03pt\raise1.42pt\hbox{$\times$}}}
\newcommand{\supp}{{\rm supp}\,}
\renewcommand{\d}{\mbox{\LARGE $\cdot $}}

\newcommand{\Div}{{\rm Div}\,}
\newcommand{\arc}[1]{{#1}\kern-16pt\raise7pt\hbox{\Large{$\frown$}}}

\newcommand{\longinto}{\lhook\joinrel\kern-3pt\hbox to
100pt{\rightarrowfill}}           

\newcommand{\propsubset}{\mbox{$\textstyle{
\subseteq_{\kern-5pt\raise-1pt\hbox{\mbox{\tiny{$/$}}}}}$}}

\newcounter{elno}                

\newenvironment{proof}{{\bf Proof}:\quad
                     }{\hfill$\odot$\par\vspace{1cm}}

\newcounter{example}[section]
\def\theexample{\thesection.\arabic{example}}

\newcounter{exercise}[section]
\def\theexercise{\thesection.\arabic{exercise}}

\newcommand{\cA}{{\cal A}}

\newcommand{\cD}{{\cal D}}
\newcommand{\cE}{{\cal E}}
\newcommand{\cF}{{\cal F}}

\newcommand{\cK}{{\cal K}}
\newcommand{\ccL}{{\cal L}}

\newcommand{\cO}{{\cal O}}
\newcommand{\cP}{{\cal P}}
\newcommand{\cQ}{{\cal Q}}

\newcommand{\cS}{{\cal S}}

\newcommand{\cU}{{\cal U}}

\newcommand{\cZ}{{\cal Z}}

\catcode`\@=11
\def\opn#1#2{\def#1{\mathop{\kern0pt\fam0#2}\nolimits}} 
\def\underrightarrow{\mathpalette\underrightarrow@}
\def\underrightarrow@#1#2{\vtop{\ialign{$##$\cr
 \hfil#1#2\hfil\cr\noalign{\nointerlineskip}%
 #1{-}\mkern-6mu\cleaders\hbox{$#1\mkern-2mu{-}\mkern-2mu$}\hfill
 \mkern-6mu{\to}\cr}}}

\def\underleftarrow{\mathpalette\underleftarrow@}
\def\underleftarrow@#1#2{\vtop{\ialign{$##$\cr
 \hfil#1#2\hfil\cr\noalign{\nointerlineskip}#1{\leftarrow}\mkern-6mu
 \cleaders\hbox{$#1\mkern-2mu{-}\mkern-2mu$}\hfill
 \mkern-6mu{-}\cr}}}
\def\:{\colon}
\let\oldtilde=\tilde
\def\tilde#1{\mathchoice{\widetilde{#1}}{\widetilde{#1}}%
{\indextil{#1}}{\oldtilde{#1}}}
\def\indextil#1{\lower2pt\hbox{$\textstyle{\oldtilde{\raise2pt%
\hbox{$\scriptstyle{#1}$}}}$}}
\def\pnt{{\raise1.1pt\hbox{$\textstyle.$}}}
\let\amp@rs@nd@\relax
\newdimen\ex@
\newdimen\bigaw@
\newdimen\minaw@
\newdimen\minCDaw@
\newif\ifCD@
\def\minCDarrowwidth#1{\minCDaw@#1}

\def\@CD{\def\A##1A##2A{\llap{$\vcenter{\hbox
 {$\scriptstyle##1$}}$}\Big\uparrow\rlap{$\vcenter{\hbox{%
$\scriptstyle##2$}}$}&&}%
\def\V##1V##2V{\llap{$\vcenter{\hbox
 {$\scriptstyle##1$}}$}\Big\downarrow\rlap{$\vcenter{\hbox{%
$\scriptstyle##2$}}$}&&}%
\def\={&\hskip.5em\mathrel
 {\vbox{\hrule width\minCDaw@\vskip3\ex@\hrule width
 \minCDaw@}}\hskip.5em&}%
\def\verteq{\Big\Vert&&}%
\def\noarr{&&}%
\def\vspace##1{\noalign{\vskip##1\relax}}\relax\let\amp@rs@nd@&\iffalse}\fi
 \CD@true\vcenter\bgroup\relax\let\\=\cr\iffalse}\fi\tabskip\z@skip\baselineskip20\ex@
 &\hfill$\m@th##$\hfill\cr}
\def\@endCD{\cr\egroup\egroup}
\def\>#1>#2>{\amp@rs@nd@\setbox\z@\hbox{$\scriptstyle
 \;{#1}\;\;$}\setbox\@ne\hbox{$\scriptstyle\;{#2}\;\;$}\setbox\tw@
 \hbox{$#2$}\ifCD@
 \global\bigaw@\minCDaw@\else\global\bigaw@\minaw@\fi
 \ifdim\wd\z@>\bigaw@\global\bigaw@\wd\z@\fi
 \ifdim\wd\@ne>\bigaw@\global\bigaw@\wd\@ne\fi
 \ifCD@\hskip.5em\fi
 \ifdim\wd\tw@>\z@
 \mathrel{\mathop{\hbox to\bigaw@{\rightarrowfill}}\limits^{#1}_{#2}}\else
 \mathrel{\mathop{\hbox to\bigaw@{\rightarrowfill}}\limits^{#1}}\fi
 \ifCD@\hskip.5em\fi\amp@rs@nd@}
\def\<#1<#2<{\amp@rs@nd@\setbox\z@\hbox{$\scriptstyle
 \;\;{#1}\;$}\setbox\@ne\hbox{$\scriptstyle\;\;{#2}\;$}\setbox\tw@
 \hbox{$#2$}\ifCD@
 \global\bigaw@\minCDaw@\else\global\bigaw@\minaw@\fi
 \ifdim\wd\z@>\bigaw@\global\bigaw@\wd\z@\fi
 \ifdim\wd\@ne>\bigaw@\global\bigaw@\wd\@ne\fi
 \ifCD@\hskip.5em\fi
 \ifdim\wd\tw@>\z@
 \mathrel{\mathop{\hbox to\bigaw@{\leftarrowfill}}\limits^{#1}_{#2}}\else
 \mathrel{\mathop{\hbox to\bigaw@{\leftarrowfill}}\limits^{#1}}\fi
 \ifCD@\hskip.5em\fi\amp@rs@nd@}
\def\@CDS{\def\A##1A##2A{\llap{$\vcenter{\hbox
 {$\scriptstyle##1$}}$}\Big\uparrow\rlap{$\vcenter{\hbox{%
$\scriptstyle##2$}}$}&}%
\def\V##1V##2V{\llap{$\vcenter{\hbox
 {$\scriptstyle##1$}}$}\Big\downarrow\rlap{$\vcenter{\hbox{%
$\scriptstyle##2$}}$}&}%
\def\={&\hskip.5em\mathrel
 {\vbox{\hrule width\minCDaw@\vskip3\ex@\hrule width
 \minCDaw@}}\hskip.5em&}
\def\verteq{\Big\Vert&}
\def\novarr{&}
\def\noharr{&&}
\def\SE##1E##2E{\slantedarrow(0,18)(4,-3){##1}{##2}&}
\def\SW##1W##2W{\slantedarrow(24,18)(-4,-3){##1}{##2}&}
\def\NE##1E##2E{\slantedarrow(0,0)(4,3){##1}{##2}&}
\def\NW##1W##2W{\slantedarrow(24,0)(-4,3){##1}{##2}&}
\def\slantedarrow(##1)(##2)##3##4{%
\thinlines\unitlength1pt\lower 6.5pt\hbox{\begin{picture}(24,18)%
\put(##1){\vector(##2){24}}%
\put(0,8){$\scriptstyle##3$}%
\put(20,8){$\scriptstyle##4$}%
\end{picture}}}
\def\vspace##1{\noalign{\vskip##1\relax}}\relax\let\amp@rs@nd@&\iffalse}\fi
 \CD@true\vcenter\bgroup\relax\let\\=\cr\iffalse}\fi\tabskip\z@skip\baselineskip20\ex@
 &\hfill$\m@th##$\hfill\cr}
\def\@endCDS{\cr\egroup\egroup}
\newdimen\TriCDarrw@
\newif\ifTriV@
\newenvironment{TriCDV}{\@TriCDV}{\@endTriCD}

\def\@TriCDV{\TriV@true\def\TriCDpos@{6}\@TriCD}
\def\@TriCDA{\TriV@false\def\TriCDpos@{10}\@TriCD}
\def\@TriCD#1#2#3#4#5#6{%
\setbox0\hbox{$\ifTriV@#6\else#1\fi$}
\TriCDarrw@=\wd0 \advance\TriCDarrw@ 24pt
\advance\TriCDarrw@ -1em
\def\SE##1E##2E{\slantedarrow(0,18)(2,-3){##1}{##2}&}
\def\SW##1W##2W{\slantedarrow(12,18)(-2,-3){##1}{##2}&}
\def\NE##1E##2E{\slantedarrow(0,0)(2,3){##1}{##2}&}
\def\NW##1W##2W{\slantedarrow(12,0)(-2,3){##1}{##2}&}
\def\slantedarrow(##1)(##2)##3##4{\thinlines\unitlength1pt
\lower 6.5pt\hbox{\begin{picture}(12,18)%
\put(##1){\vector(##2){12}}%
\put(-4,\TriCDpos@){$\scriptstyle##3$}%
\put(12,\TriCDpos@){$\scriptstyle##4$}%
\end{picture}}}
\def\={\mathrel {\vbox{\hrule
   width\TriCDarrw@\vskip3\ex@\hrule width
   \TriCDarrw@}}}
\def\>##1>>{\setbox\z@\hbox{$\scriptstyle
 \;{##1}\;\;$}\global\bigaw@\TriCDarrw@
 \ifdim\wd\z@>\bigaw@\global\bigaw@\wd\z@\fi
 \hskip.5em
 \mathrel{\mathop{\hbox to \TriCDarrw@
{\rightarrowfill}}\limits^{##1}}
 \hskip.5em}
\def\<##1<<{\setbox\z@\hbox{$\scriptstyle
 \;{##1}\;\;$}\global\bigaw@\TriCDarrw@
 \ifdim\wd\z@>\bigaw@\global\bigaw@\wd\z@\fi
 \mathrel{\mathop{\hbox to\bigaw@{\leftarrowfill}}\limits^{##1}}
 }
 \CD@true\vcenter\bgroup\relax\let\\=\cr\iffalse}\fi
 \tabskip\z@skip\baselineskip20\ex@
 \ifTriV@
 \halign\bgroup
 &\hfill$\m@th##$\hfill\cr
#1&\multispan3\hfill$#2$\hfill&#3\\
&#4&#5\\
&&#6\cr\egroup%
\else
 \halign\bgroup
 &\hfill$\m@th##$\hfill\cr
&&#1\\%
&#2&#3\\
#4&\multispan3\hfill$#5$\hfill&#6\cr\egroup
\fi}
\def\@endTriCD{\egroup}

\title{\bf\sc Albanese and Picard 1-motives}
\author{by Luca {\sc Barbieri-Viale}
and Vasudevan {\sc Srinivas}}
\date{}
\begin{document}

\maketitle

\begin{abstract}
Let $X$ be an $n$-dimensional algebraic variety over a field of
characteristic zero. We describe algebraically defined Deligne
1-motives $\Alb^{+}(X)$,  $\Alb^{-}(X)$, $\Pic^{+}(X)$ and
$\Pic^{-}(X)$ which generalize the classical Albanese
and Picard varieties of a smooth projective variety. We compute
Hodge, $\ell$-adic and De Rham realizations proving Deligne's
conjecture for $H^{2n-1}$, $H_{2n-1}$, $H^1$ and $H_1$.

We investigate functoriality, universality, homotopical invariance
and invariance under formation of projective bundles. We compare
our cohomological and homological 1-motives for normal schemes.
For proper schemes, we obtain an Abel-Jacobi map from the
(Levine-Weibel) Chow group of zero cycles to our cohomological Albanese
1-motive which is the universal regular homomorphism to semi-abelian
varieties. By using this universal property we get ``motivic'' Gysin maps
for projective local complete intersection morphisms.

\end{abstract}

{\small \tableofcontents}

\section{Introduction}
This paper is motivated by Deligne's conjecture
that 1-motives obtained from the mixed Hodge structure on the
cohomology of an algebraic variety would be ``algebraically defined'' (see
\cite[\S 10.4.1]{D} and \cite{DM}). Deligne (\cite[\S 10.1.3]{D}) observed
that a torsion free mixed Hodge structure $H$ (\ie such that $H_{\sZ}$ is
torsion-free), which is of Hodge type $\{(0,0), (0,-1), (-1,0),
(-1,-1)\}$, and such that $\gr_{-1}^W(H)$ is polarizable, yields
{\it i)}\, a semi-abelian variety $G$, whose abelian quotient is the
abelian variety given by $\gr_{-1}^W(H)$, together with  {\it ii)}\, a
homomorphism $u$ from the lattice $L = \gr_{0}^W(H_{\sZ})$ to the group
$G(\C)$, induced by the canonical map $H_{\sZ}\to H_{\sC}$.

Deligne called such a complex of group schemes $[L\by{u} G]$ a {\em
1-motive over $\C$}, and showed that the category of 1-motives over $\C$
is equivalent to the category of torsion free mixed Hodge structures of
the above type. Thus any such 1-motive $M=[L\by{u} G]$ has a {\em Hodge
realization}\, $T_{Hodge}(M)$, \ie there is a unique (up to isomorphism)
torsion-free mixed Hodge structure $T_{Hodge}(M)$ such that $M$ can be
obtained from $T_{Hodge}(M)$ as above.  Deligne (\cite[\S 10.1.11]{D})
also defined the {\em $\ell$-adic}\, and {\em De Rham realizations}\, of a
1-motive $M$, denoted by $T_{\ell}(M)$ and $T_{DR}(M)$, respectively (see
Section~\ref{pre} for more details).

\subsection*{The conjecture}
Deligne's conjecture, in particular, is that if $X$ is an $n$-dimensional
algebraic variety over a field $k$ of characteristic 0, then there are
``algebraically defined'' 1-motives, also defined over $k$, compatible
with base change to extension fields, such that {\it i)}\, when $k=\C$,
their Hodge realizations are respectively isomorphic to the mixed Hodge
structures on
\[H^{2n-1}(X,\Z(n))/{\rm (torsion)},\;\;
H_{1}(X,\Z)/({\rm torsion}),\;\; H^{1}(X,\Z (1)),\;\;
H_{2n-1}(X,\Z (1-n))/({\rm torsion})\]
{\it ii)}\, if $k$ is algebraically closed, their $\ell$-adic and De Rham
realizations are naturally isomorphic to the corresponding $\ell$-adic and
De Rham (co)homology {\it iii)}\, the above identifications
are compatible with other structures, like comparison isomorphisms,
filtrations, Galois action, {\it etc\/}. Our goal in this paper is to
prove these statements. Moreover, we obtain some geometric properties of
our constructions.

We recall that the case $n=1$, \ie when $X$ is a curve, is already treated
by Deligne (\cite[\S 10.3]{D}), and the case when $X$ smooth and proper
corresponds to the well known trascendental descriptions
of the Albanese and Picard varieties. Our construction of Albanese and
Picard 1-motives generalizes these cases. In the general case,
for $n\geq 1$, we propose the following dictionary:
$$\begin{array}{lcl}
\mbox{\it Mixed Hodge Structure} & \hspace*{1cm} &\mbox{\it 1-Motive}\\
& &\\
H^{2n-1}(X,\Z (n)) & &\Alb^+(X)  \\
H_{2n-1}(X,\Z (1-n)) & & \Pic^-(X)\\
H^{1}(X,\Z (1)) & & \Pic^+(X)  \\
H_{1}(X,\Z) & & \Alb^-(X)  \\
\end{array}$$
Here, $\Alb^+(X)$ is the ``cohomological Albanese'' 1-motive obtained from
the mixed Hodge structure $H^{2n-1}(X,\Z (n))/({\rm torsion})$ and,
dually, $\Pic^-(X)$ is the ``homological Picard'' 1-motive obtained from
$H_{2n-1}(X,\Z (1-n))/{\rm (torsion)}$, {\it etc.}\, The 1-motive
$\Alb^+(X)$ is the Cartier dual of $\Pic^-(X)$ and $\Pic^+(X)$ is the
Cartier dual of $\Alb^-(X)$. In case $X$ is singular, we have that
$\Alb^+(X)\neq\Alb^-(X)$ in general, because of the possible failure of
Poincar\'e duality. When $n=1$, $\Alb^+(X)$ and $\Pic^+(X)$ (and dually,
$\Alb^-(X)$ and $\Pic^-(X)$) coincide.

We recall that the geometric definition of the ``cohomological Picard and
homological Albanese'' 1-motives of a smooth, but possibly non-proper
scheme $X$, goes back to Serre's explicit construction of its Albanese
variety, see \cite{SER}; in fact, Serre's Albanese variety was
defined as the Cartier dual of the 1-motive
$$\Pic^+(X)\df[\Div_Y^0(\bar X)\to \Pic^0(\bar X)]\mbox{
\hspace*{1cm}($X$ smooth)}$$
where $\bar X$ is a smooth compactification of $X$ with boundary $Y$,
$\Div_Y^0(\bar X)$ is the free abelian group of divisors which are
algebraically equivalent to zero and supported on $Y$, being mapped
canonically to $\Pic^0(\bar X)$. On the other hand, a geometric
construction of $\Alb^+$ or $\Pic^-$ for a smooth open variety is more
difficult and it appears to be new as well.

Following the construction in \cite{LI}, in the paper of Ramachandran \cite{RA}
a geometric construction of $\Pic^+(X)$ and $\Alb^-(X)$ was proposed for
varieties with a singular closed point obtained by collapsing a finite set
of closed points in a smooth open variety; in a subsequent paper, see
\cite{RA2}, he proposed, independently, definitions of Albanese and
Picard motives corresponding to our $\Pic^+(X)$ and $\Alb^-(X)$.
Ramachandran has recently announced \cite{RA3} a proof of the algebraicity
(up to isogeny) of certain 1-motives built out of $H^i(X,\Q(1))$ for
$i\leq \dim X +1$.

Apart from Ramachandran's work, a related paper by Carlson \cite{C} on
analogues of Abel's theorem for $H^2$ of some singular surfaces (see also
\cite{G}), and the ``Hodge-Lefschetz 1-motives'' considered in \cite{BR}
(see also \cite{BST} and  \cite{BSL}) we do not know any results on
Deligne's conjecture (1972).

\subsection*{The results}
Our definition of $\Pic^-(X)$ is roughly the following (see Sections~2.1
and 2.2 below for a more precise statement). Let $X$ be any
equidimensional algebraic variety over an algebraically closed field $k$
of characteristic zero. Let $f:\tilde X\to X$ be a resolution of
singularities and let $\bar X$ be a smooth compactification of $\tilde X$
with normal crossing boundary divisor $Y$.

Let $S$ be the singular locus of $X$ and let $\bar S$ be the closure of
$f^{-1}(S)$ in $\bar X$.  Then we let $\Div_{\bar S}^0(\bar X,Y)$ be the
group of divisors supported on $\bar S$ which are {\it i)}\, disjoint from
$Y$ (\ie are linear combinations of compact components of $f^{-1}(S)$),
and {\it ii)}\, are algebraically equivalent to zero relative to $Y$.  We
let $\Div_{\bar S/S}^0(\bar X,Y)$ be the subgroup of those divisors which
have vanishing push-forward (as Weil divisors) along $f$.

We can show the existence of a group scheme $\Pic (\bar X, Y)$
associated to isomorphism classes of pairs $(\ccL, \varphi)$ such that
$\ccL$ is a line bundle on $\bar X$ and $\varphi :\ccL\mid_Y\cong \cO_Y$
is a trivialization on $Y$. The connected component of the identity
$\Pic^0(\bar X, Y)$ is a semi-abelian variety, which can be represented
as an extension
$$ 0\to\frac{H^0(Y,\cO_Y^*)}{\im H^0(\bar X,\cO_{\bar X}^*)}\to \Pic^0(\bar
X, Y) \to \ker^0(\Pic^0(\bar X)\to\oplus_i\Pic^0(Y_i)) \to 0$$
where $Y=\cup_i Y_i$ is expressed as a union of (smooth) irreducible
components. The mapping which takes a divisor $D$ disjoint from $Y$ to the
class of the pair $(\cO_{\bar X}(D),1)$ in $\Pic (\bar X, Y)$ yields the
``homological Picard'' 1-motive
$$\Pic^-(X)\df [\Div_{\bar S/S}^0(\bar X,Y)\to\Pic^0(\bar X,Y)].$$
The ``cohomological Albanese'' 1-motive $\Alb^+(X)$ is defined to be the
Cartier dual of $\Pic^-(X)$ (see Section~3.1); a ``concrete'' description
of it is also given when $X$ is either smooth or proper.

The definition of $\Pic^{+}(X)$ is obtained by generalizing Serre's
construction of the generalized Albanese variety to smooth simplicial
schemes (see Sections~4.1 and 4.2 for the details). Let $X$ be a
variety over an algebraically closed field $k$ of characteristic
0. Let $X_{\d}$ be a smooth proper hypercovering of $X$, and consider a
proper smooth compactification $\bar X_{\d}$ with normal crossing boundary
$Y_{\d}$ (we refer to \cite[\S 6.2]{D} for the existence of such a
hypercovering and compactification). Let $\Div_{Y_{\d}}(\bar X_{\d})$ be
the subgroup of divisors on $\bar X_0$ supported on $Y_0$ with zero
pullback on $\bar X_1$, \ie by definition
$$\Div_{Y_{\d}}(\bar X_{\d})\df
\ker (\Div_{Y_0}(\bar X_0)\longby{d_0^*-d_1^*}\Div_{Y_1}(\bar X_1)).$$
We consider the simplicial Picard functor
$$T\leadsto \bPic (T\times \bar X_{\d})\df \HH^1(T\times \bar X_{\d},
\cO^*_{T\times \bar X_{\d}})$$
and we show that the associated sheaf $\bPic_{{\bar X_{\d}}/k}$ (with
respect to the $fpqc$-topology) is representable by a group scheme locally
of finite type over $k$, whose connected component of the identity
$\bPic^0_{{\bar X_{\d}}/k}$ over $k = \bar k$ is an extension of the
abelian scheme $\ker^0 (\Pic_{{\bar X_0}/k}^0\to\Pic_{{\bar X_1}/k}^0)$ by
the torus given by
$$\frac{\ker (H^0(\bar X_1,\cO_{\bar X_1}^*)\to
H^0(\bar X_2,\cO_{\bar X_2}^*))}{\im (H^0(\bar X_0,\cO_{\bar X_0}^*)\to
H^0(\bar X_1,\cO_{\bar X_1}^*))}.$$
Let $\Div_{Y_{\d}}^0(\bar X_{\d})$ denote the subgroup of those divisors
which are mapped to $\bPic^0_{{\bar X_{\d}}/k}(k)$ under the canonical
mapping. We then define the ``cohomological Picard'' 1-motive
of the variety $X$ as
$$\Pic^+(X)\df [\Div_{Y_{\d}}^0(\bar X_{\d})\to \bPic^0(\bar X_{\d})].$$
The ``homological Albanese'' 1-motive $\Alb^-(X)$ is defined to be the
Cartier dual of $\Pic^+(X)$ (see Section~5.1).

We show that $\Pic^{-}(X)$, $\Alb^{+}(X)$, $\Pic^{+}(X)$ and $\Alb^{-}(X)$
do have the appropriate Hodge, De Rham and $\ell$-adic realizations
(in Sections~2.4--2.6, 3.3, 4.3--4.5 and 5.3 respectively).
We mostly deal with the geometric case, \ie we consider varieties $X$
over an algebraically closed field $k$; the case when $k$ is
not algebraically closed is considered in Section~\ref{rational}.

We show as well that our definitions are functorial and independent of
choices of resolutions or compactifications (\eg see Section~2.3); we remark
(in Section~6.1) that $\Alb^+$ can be contravariant functorial only for
morphisms between varieties of the same dimension, and similarly $\Pic^-$ is
covariant functorial for such maps. We then show the homotopical
invariance of $\Pic^+$ (and hence dually of $\Alb^-$), and that $\Pic^-$
and $\Pic^+$ (and dually, the corresponding Albanese 1-motives)  are
invariant under formation of projective bundles (see Sections~6.2).

For proper $X$, we remark that our ``cohomological'' Albanese 1-motive
$\Alb^+(X)$ is a quotient of Serre's Albanese of the regular locus $X_{\rm
reg}$, \ie we have an extension
$$ 0\to T(S)\to\Alb^-(X_{\rm reg})\to\Alb^+(X)\to 0$$
where $T(S)$ is a torus whose character group is a sub-lattice of the
lattice of Weil divisors which are supported on the singular locus $S$.
Thus, if $X$ is also irreducible and normal, then $T(S)=0$, and further,
any non-zero Cartier divisor supported on the exceptional locus of a
resolution is not numerically equivalent to zero;  therefore,
$\Alb^-(X_{\rm reg})= \Alb^+(X)$ is an abelian variety which is isomorphic
to the Albanese variety $\Alb (\tilde X)$ of any resolution of
singularities $\tilde X$ of $X$. In general, $\Alb^-(X_{\rm reg})$ is a
torus bundle over $\Alb(\tilde X)$ whose pull-back to $X_{\rm reg}$ (under
a suitable Albanese mapping) is canonically trivialized. Thus, after
choosing appropriate base points, there is a (canonical) section $$a^-:
X_{\rm reg} \to \Alb^-(X_{\rm reg})$$ which is a universal morphism to
semi-abelian varieties in the sense of Serre \cite{SEM} (see Section~5.2).

We then show (in Section~6.3) that $a^-$ factors through rational
equivalence yielding a ``motivic'' Abel-Jacobi mapping
$$a^+: CH^n(X)_{\deg 0}\to \Alb^+(X)$$
from the Levine-Weibel ``cohomological'' Chow group \cite{LW} of
zero-cycles on a projective variety $X$. We also prove that $a^+$ is the
universal regular homomorphism to semi-abelian varieties (compare with
\cite{ESV} and \cite{BS}). By using this universal property we get
``motivic'' Gysin maps for projective local complete intersection
morphisms, \ie for such a morphism $g:X' \to X$ we get a push-forward
$$g_*^+:\Alb^+(X')\to\Alb^+(X)$$
and, dually, a pull-back $g^*_-:\Pic^-(X)\to\Pic^-(X')$.

We note that the isogeny classes of our 1-motives define objects in the
triangulated category of mixed motives of Voevodsky, since it contains
Deligne's 1-motives (tensor $\Q$). Therefore we can view our constructions
as determining ``Picard and Albanese mixed motives'' as well.

We finally remark that these results were previously announced in \cite{BSAP}.

\subsection*{Some further questions}

We expect purely algebraic proofs for the Lefschetz theorem on inclusions
of general complete intersections $g:Y\into X$ (\ie $g_*^+$ and $g^*_-$
would be isomorphisms in this case, if $\dim Y\geq 2$) as well as Roitman
theorems on torsion zero-cycles (\ie  $a^+$ would be an isomorphism on
torsion, see \cite{BPW} and \cite{BS} for the case $k=\C$, and \cite{SA}
for the homological case): these matters are of independent interest,
and we hope to treat them elsewhere. In the context of algebraic proofs,
it seems desirable as well to have such a proof that $\Pic^+(X)$ and
$\Alb^-(X)$ are independent of the choices of hypercovering and
compactification (see Remark~\ref{rempicplus}). In fact, the underlying
philosophy of the theory of 1-motives suggests that it should be possible
(or at least desirable) to obtain all constructions and properties
``intrinsically'', without recourse to the use of any specific realization
functor. From this point of view, another problem is to prove ``directly''
that the Gysin maps for projective local complete intersection morphisms
are independent of the factorization (see Remark~\ref{remgysin}).

It is natural to ask whether there is an analogue of our results in
positive characteristic. After the work of de Jong, there are smooth
proper hypercoverings in this context as well, which suggests that one
could possibly extend the definitions of $\Pic^+$ and $\Alb^-$ to this
case. However, since our definitions of $\Pic^-$ and $\Alb^+$, and the
proofs that they have the correct realizations, make use of resolution
of singularities and duality theory, it is not clear to us how these might
extend to positive characteristics. In positive characteristics, one
also needs to better understand what would play the role of the De Rham
realization.

More generally, since Deligne has defined a notion of a 1-motive over a
base scheme $S$, we could ask for the appropriate families $X\to S$ for
which it is possible to define $\Pic^+(X/S)$, $\Pic^-(X/S)$, and the
corresponding Albanese 1-motives, as 1-motives over $S$. Going still
further, one could speculate about possible analogous 1-motives in the
context of Arakelov geometry.

\subsection*{Acknowledgements}
The second author wishes to acknowledge support from the IFCPAR Project
``Geometry''(No. 1601-2), which made possible a visit to the Univ.
Paris~VII, concurrently with a visit of the first author to the Institut
Henri Poincar\'e, and enabled useful discussions on this work with O.
Gabber, L.~Illusie, F.~Oort and M.~Raynaud. In particular, we thank
Raynaud for providing us with the reference~\cite{BZ}. The first author is
deeply grateful to the Institut Henri Poincar\'e, the Italian CNR and
MURST, the European Union Science Plan ``K-Theory and Algebraic Groups''
and the Tata Institute of Fundamental Research for support and
hospitality.

\section*{Notations}
We are mainly concerned with schemes locally of finite
type over a base field $k$ of characteristic zero, which is assumed to be
algebraically closed in most of this paper; we will consider non
algebraically closed fields in Section~\ref{rational}. We tacitly assume
that our schemes are reduced and separated, unless explicitly mentioned
otherwise. A {\em variety}\, will be a reduced, separated $k$-scheme of
finite type. We will often tacitly identify a variety over $k=\bar k$ with
its set of closed points. The hypothesis of zero characteristic is
repeatedly used, often without explicit mention, for example via the
existence of resolutions of singularities.

We denote by $X_{\d}$ a simplicial $k$-scheme, whose components
$X_{i}$ are $k$-schemes, and we denote by $d^i_j: X_{i}\to X_{i-1}$,
$0\leq j\leq i$, the face maps; we omit upper indices if there is no risk
of misunderstanding, \eg we may write $d_0$ and $d_1$ for the two faces
map from $X_{1}$ to $X_{0}$. We will also sometimes identify a $k$-scheme
$X$ with the ``constant'' simplicial scheme $X_{\d}$ it defines, where
$X_n=X$ for all $n\geq 0$, and all face and degeneracy morphisms are the
identity; if $\pi:X_{\d}\to X$ is the augmentation, then we note that for
any sheaf of abelian groups $\cF$ on $X$, the canonical map $\cF\to {\bf
R}\,\pi_*(\pi^*\cF)$ is an isomorphism, and we have canonical isomorphisms
$H^*(X,\cF)\cong \HH^*(X_{\d},\pi^*\cF)$.

For a $\C$-variety $X$ we will denote by $H^*(X,\Z (\cdot))$
(resp. $H_*(X,\Z(\cdot))$) the singular cohomology (resp. homology)
group of the associated analytic space as well as the (Tate twisted) mixed
Hodge structure on it. Concerning mixed Hodge structures we will use
Deligne's notation \cite{D}:  in particular, we will denote by $W_iH$
the weight filtration on $H_{\sQ}$ (and if $H$ is torsion free, on
$H_{\sZ}$ as well), and by $F^iH$ the Hodge filtration on $H_{\sC}$.

For a simplicial scheme $X_{\d}$ and a simplicial abelian sheaf
$\cF_{X_{\d}}$ we will denote by $\HH^*(X_{\d},\cF_{X_{\d}})$ the
cohomology groups obtained from the right derived functor of the following left
exact functor
$$\cF_{X_{\d}}\leadsto \ker (\Gamma
(X_0,\cF_{X_{0}})\longby{d_0^*-d_1^*} \Gamma (X_1,\cF_{X_{1}})).$$
The same conventions as above apply to simplicial $\C$-schemes $X_{\d}$
and the mixed Hodge structure on $\HH^*(X_{\d},\Z (\cdot))$.

We denote duals by $(-)^{\vee}$ with the following conventions: if $G$
is a group scheme of additive type then $G^{\vee}$ is $\Hom (G,\G_a)$; if
$G$ is
a torus or it is locally constant and torsion free then $G^{\vee}$ is
$\Hom (G,\G_m)$; if $A$ is an abelian
variety then $A^{\vee}$ is $\Pic^0(A)$; if $H$ is a mixed Hodge structure
then $H^{\vee}$ is the internal $\Hom (H,\Z (1))$.

We denote by $\Div (X)$ the group of Weil divisors on an
equidimensional variety $X$. If $Y$, $Z$ are closed subschemes of $X$, we
denote by $\Div_Z(X)\subset \Div (X)$ the subgroup of divisors which are
supported on $Z$, and by $\Div(X,Y)\subset \Div (X)$ the subgroup of
divisors which have the support disjoint from $Y$; finally set
$\Div_Z(X,Y)=\Div_Z(X)\cap \Div(X,Y)$.

For any (possibly singular) variety $X$ we let denote by $CH_d(X)$ the
``homological'' (Fulton) Chow groups \cite{FU} of $d$-dimensional cycles
on $X$. We denote by $CH^n(X)$ the ``cohomological'' (Levine-Weibel)
Chow group \cite{LW} of zero-cycles supported on the regular locus of an
$n$-dimensional quasi-projective variety $X$ over an algebraically closed
field.

If $f:G_1\to G_2$ is a homomorphism of $k$-group schemes, $\ker^0 f$
will denote the identity component of the kernel of $f$.

\section{Preliminaries on 1-motives}\label{pre}
For the sake of exposition, and to fix notation and terminology, we
collect some general facts concerning 1-motives.

\subsection{Deligne's definition}
Let $S$ be any scheme. We will denote by $M=(L,A,T,G,u)$ a {\it 1-motive}\,
over
$S$, \ie an extension $G$ of an abelian scheme $A$ by a torus
$T$ over $S$, a group scheme $L$ which is, locally for the
\'etale topology on $S$, isomorphic to a finitely-generated
free abelian constant group, and an $S$-homomorphism
$L\by{u} G$ (see \cite[\S 10]{D}).

Diagrammatically a 1-motive $M=(L,A,T,G,u)$ can be represented as
$$ \begin{array}{ccc}
 & L & \\
& \ \ \downarrow u & \\
 1\to \ T& \to\ G\ \to & A\ \to 0
\end{array}$$
and can be regarded also as defining a complex of group schemes $M=[L\by{u}
G]$,
where $L$ is in degree $-1$ and $G$ is in degree $0$.

A group scheme $G$ which is an extension of an abelian scheme $A$ by a
torus $T$ is also usually called a {\it semi-abelian scheme}, and we are
not going to distinguish it from the 1-motive which it defines in a
canonical way (\ie by taking $L$ to be zero). The same convention applies
to the case of an abelian variety $A$ (identified with the 1-motive
$(0,A,0,A,0)$) or a torus. A lattice $L$ determines a 1-motive $[L\to 0]$,
which we denote by $L[1]$ (consistent with the notation when considered as
a complex of group schemes).

A {\em morphism of 1-motives}\, is a morphism of the corresponding complexes of
group schemes. Moreover, there is a natural full embedding of the category of
1-motives over $S$  into the derived category of bounded
complexes of sheaves for the $fppf$-topology on $S$ (\cf
\cite[Prop.2.3.1]{RAY}).

A 1-motive $M$ is canonically equipped with an increasing filtration by
sub-1-motives as follows:
$$W_i(M) =\left \{\begin {array}{rr} M & i\geq 0\\ G
& i= -1\\ T & i= -2\\ 0 & i \leq -3 \end{array} \right. $$
In particular we have $\gr_{-1}^W(M) = A$.

A complex of 1-motives is {\em exact} if it determines an exact sequence
of complexes of group schemes. For example, associated to any 1-motive
$M=(L,A,T,G,u)$ there is a functorial exact sequence of 1-motives
\begin{equation}\label{motexseq}
0\to G \to M\to L[1]\to 0
\end{equation}
where $L[1]=\gr_0^W(M)=[L\to 0]$.

\subsection{Hodge realization} We recall that the
Hodge realization $T_{Hodge}(M)$ ($T(M)$ for short) of a
1-motive $M$ over $k=\C$ (see the construction by Deligne in
\cite[10.1.3]{D}) is the mixed Hodge structure given by  the lattice
$T_{\sZ}(M)$
obtained by the pull-back of  $u : L\to G$ along $\exp : \Lie
(G) \to G$, with the weight filtration
$$W_iT(M) \df\left \{\begin
{array}{rr} T_{\sZ}(M) & i\geq 0\\
H_1(G) & i= -1\\ H_1(T) & i= -2\\ 0 & i \leq -3
\end{array} \right. $$
The Hodge filtration is defined by
$$F^0(T_{\sZ}(M)\tensor\C)\df \ker (T_{\sZ}(M)\tensor\C\to \Lie(G)),$$
whence $\gr_{-1}^WT(M)\cong H_1(A,\Z)$ as pure Hodge structures of weight
$-1$. The functor
$$M \leadsto T_{Hodge}(M)$$
is an equivalence between the category of 1-motives and the
full subcategory of torsion free $\Z$-mixed Hodge structures
of type
$$\{(0,0), (0,-1), (-1,0), (-1,-1)\}$$
such that $\gr_{-1}^W(H)$ is polarizable. In fact, Deligne
(\cite[\S 10.1.3]{D}) observed that such a torsion free
mixed Hodge structure $H$  yields
{\it (i)}\, an abelian variety $A$ with
$$A(\C) = \frac{\gr_{-1}^W(H_{\sC})}{H_{\sZ} + F^0\gr_{-1}^W(H_{\sC})}$$
{\it (ii)}\, an algebraic torus $T$ with character
group $\gr_{-2}^W(H_{\sZ})$, so that
$$T(\C) =\Hom (\gr_{-2}^W(H_{\sZ}), \C^*)$$
and
{\it (iii)}\, a complex algebraic group $G$ with
$$G(\C) = \frac{W_{-1}(H_{\sC})}{W_{-1}(H_{\sZ}) + F^0\cap
W_{-1}(H_{\sC})}$$
which is an algebraic extension of $A$ by $T$;
moreover, the canonical map $H_{\sZ}\to H_{\sC}$ yields
{\it (iv)}\, a homomorphism $u$ from the lattice $L =
\gr_{0}^W(H_{\sZ})$ to the group $G(\C)$. Deligne considered such a set of
data {\it (i)---(iv)}\, as defining a 1-motive over $\C$, and showed that
it is equivalent to the given mixed Hodge structure.

Thus any 1-motive $M=(L, A, T, G, u)$ over $\C$ has a Hodge realization
$T_{Hodge}(M)$ and, conversely, any such mixed Hodge structure yields,
canonically, a 1-motive. The exact sequence (\ref{motexseq}) gives rise to
an exact sequence of Hodge realizations
\begin{equation}\label{hodgeexseq}
0\to T_{Hodge}(G)\to T_{Hodge}(M)\to T_{Hodge}(L[1])\to 0
\end{equation}

For example, any abelian variety $A$ over $\C$ considered as a 1-motive (\ie we
regard $A$ as $(0,A,0,A,0)$), has Hodge realization $H_1(A,\Z)$; in
particular, for a non-singular complete variety $X$ over $\C$, the classical
Albanese variety $M=\Alb (X)$ has Hodge realization ($n = \dim X$)
$$T(M) = H_1(\Alb (X),\Z)\cong H_1(X,\Z)/{\rm (torsion)}\cong
H^{2n-1}(X,\Z (n))/{\rm (torsion)}$$
because of the canonical isomorphism $\Alb (X)\cong J^n(X)$ where
$J^n(X)$ is the cohomological (Griffiths) intermediate
jacobian; for a smooth projective variety the Hodge structures
on $H_1(X,\Z)$ and $H^{2n-1}(X,\Z (n))$ are canonically
isomorphic (by Poincar\'e duality) and they both correspond to the
Albanese
variety.

\subsection{$\ell$-adic and \'etale realization} Let $M=[L\by{u} G]$ be a
1-motive over $S$ which we consider as a complex of $fppf$-sheaves over
$S$ with $L$ in degree $-1$ and $G$ in degree $0$. For any fixed integer
$m$ we let $T_{\sZ/m}(M)$ be $H^{-1}(M/m)$ where $M/m$ is the
cone of the multiplication by $m$ on $M$. Then
$T_{\sZ/m}(M)$ is a finite group scheme which is flat over $S$, and is
\'etale if $S$ is defined over $\Z [\frac{1}{m}]$.

For $S=\Spec (k)$ and $k=\bar k$ we then have
$$T_{\sZ/m}(M)(k) = \frac{\{(x,g)\in L\times G(k)\mid
u(x)=-mg\}}{\{(mx,-u(x)) \mid x \in L\}}$$
If $\ell$ is a prime number then the
$\ell$-adic realization $T_{\ell}(M)$ is simply defined to be the inverse limit
over $\nu$ of $T_{\sZ/\ell^{\nu}}(M)$. $T_{\ell}(M)$ is the $\ell$-adic
Tate module of an $\ell$-divisible group. The $\ell$-adic realization of
an abelian variety $A$ is the $\ell$-adic Tate module of $A$; the
$\ell$-adic realization of a lattice $L$ is $L\otimes_{\Z} \Z_{\ell}$.

If $S=\Spec (k)$ and $k=\bar k$ is of characteristic zero then
$$\hat{T}(M)\df\liminv{m}T_{\sZ/m}(M)=\prod_{\ell} T_{\ell}(M).$$
We call $\hat{T}(M)$ the {\em \'etale realization} of $M$.
In particular, if $k=\C$ then $\hat{T}(M) = T_{\sZ}(M)\tensor\hat{\Z}$
because the complex $T_{\sZ}(M)\to \Lie G$ is quasi-isomorphic to
$M(\C)$, and therefore $M/m(\C)$ is quasi-isomorphic to
$(T_{\sZ}(M)\tensor\Z/m) [+1]$.

The exact sequence (\ref{motexseq}) of 1-motives yields a long exact
sequence of cohomology groups
$$H^{-2}(L[1]/m)\to H^{-1}(G/m)\to H^{-1}(M/m)\to
H^{-1}(L[1]/m)\to H^{0}(G/m)$$
where $H^{-2}(L[1]/m)=\ker (L\by{m}L)$ is clearly zero, and
$H^{0}(G/m)=\coker (G\by{m}G)$ vanishes since multiplication by $m$ is an
epimorphism of $fppf$-sheaves. In the sequence above we are left with
finite group schemes, and thus, by taking the inverse limit on $m$, the
sequence yields the following short exact sequence
\begin{equation}\label{etexseq}
0\to\hat{T}(G)\to\hat{T}(M)\to \hat{T}(L[1])\to 0
\end{equation}
The exact sequence above is clearly functorial with respect to maps
of 1-motives, since it is obtained from (\ref{motexseq}) by applying the
functor $\hat{T}$; it is the \'etale analogue of (\ref{hodgeexseq}).

We will later make use of the following fact.
\begin{propose}\label{fhat} The \'etale realization functor $\hat{T}$ from
the category of 1-motives over $k=\bar k$ to abelian groups is faithful,
and further, it reflects isomorphisms (\ie if $M\to M'$ is a map of
1-motives such that $\hat{T}(M)\cong\hat{T}(M')$ then $M\to M'$ is an
isomorphism in the category of 1-motives).
\end{propose}
\begin{proof}
Consider $M=[L\by{u} G]$, $M'=[L'\by{u'}G']$ and $f:M\to M'$. Now
$\hat{T}$ is clearly an additive functor; hence, in order to show that
$\hat{T}$ is
faithful, we just need to show that $\hat{T}(f)=0$ implies $f=0$. By making
use of the exact sequence (\ref{etexseq}) we can see that it is enough to
check
it seperately for maps of semi-abelian schemes or lattices. Since torsion
points are Zariski dense in a semi-abelian scheme over $k=\bar k$,
$\hat{T}(f)=0$ implies
$f=0$ for morphisms $f$ between semi-abelian schemes. Finally
$\hat{T}(L[1])=L\otimes\hat{\Z}$ which is clearly faithful.

If $M\to M'$ induces an isomorphism
$\hat{T}(M)\cong \hat{T}(M')$ then by (\ref{etexseq}) we have that
$\hat{T}(G)$ injects into $\hat{T}(G')$ and $\hat{T}(L[1])$ surjects onto
$\hat{T}(L'[1])$, therefore we have an extension of lattices
$$0\to L''\to L\to L'\to 0$$
Moreover by the snake lemma applied to the resulting diagram given by
(\ref{etexseq}) we get that
$$\hat{T}(L''[1])\cong \frac{\hat{T}(G')}{\hat{T}(G)}.$$
Now we have that $F = \ker(G\to G')$ is a finite group, since
$\hat{T}(G)\into
\hat{T}(G')$; we can see that
$$F\cong\frac{\hat{T}(G/F)}{\hat{T}(G)}\into \frac{\hat{T}(G')}{\hat{T}(G)}.$$
Thus $F=0$, since it injects into $\hat{T}(L''[1])$ which is torsion free.
If we let $G''$ denote the quotient of $G$ by $G'$, we then get the following
exact sequence of complexes
$$0\to [L''\to 0]\to [L\to G]\to [L'\to G']\to [0\to G'']\to 0$$
Applying $\hat{T}$ we have that the composition of the following maps
$$\hat{T}(L''[1])\to\hat{T}(M)\by{\cong}\hat{T}(M')\to\hat{T}(G'')$$
is the zero map as well as an isomorphism, therefore $\hat{T}(L''[1])=
\hat{T}(G'')=0$ whence $L''= G''=0$, \ie $M\by{\cong} M'$.
\end{proof}

\subsection{De Rham realization} The De Rham realization
of a 1-motive $M=[L\by{u} G]$ over an algebraically closed field
$k$ is obtained {\it via}\, Grothendieck's interpretation of
$H^1_{DR}$ (\cf \cite[\S 4]{MM}, \cite[10.1.7]{D} and \cite{DI}).
Consider $\G_a$ as a complex of $k$-group schemes concentrated in degree
$0$. Then, for any 1-motive $M$ over $k$, we have $\Hom (M,\G_a) =0$, and
there is an extension
$$0\to \Ext (L[1],\G_a)\to \Ext (M,\G_a)\to\Ext (G,\G_a)\to 0$$
where $\Ext (G,\G_a)$ is canonically identified with the Lie algebra
of the dual of the abelian variety $A$ (the abelian quotient of the
semi-abelian variety $G$), and $\Ext(L[1],\G_a)=\Hom(L,\G_a)$. Hence
the $k$-vector space $\Ext (M,\G_a)$ is finite dimensional.

By general arguments (\cf \cite{MM}, \cite{D}) $M$ has a universal
$\G_a$-extension $M^{\natural}$, in Deligne's notation \cite[10.1.7]{D},
where $M^{\natural}=[L\by{u^{\natural}} G^{\natural}]$ is a complex of
$k$-group schemes which is an extension of $M$ by the vector space
$\Ext (M,\G_a)^{\vee}$, considered as a complex in degree zero. In fact, we
have a diagram
$$\begin{array}{ccccccc}
&0&&0&&&\\
&\downarrow&&\downarrow&&&\\
0\to &\Ext (G,\G_a)^{\vee}&\to &\Ext (M,\G_a)^{\vee}&\to &
\Ext (L[1],\G_a)^{\vee}&\to 0\\
&\downarrow&&\downarrow&&\ \ \veq\mbox{\small def}&\\
0\to& G^{\sharp}&\to & G^{\natural}&\to & L^{\natural}&\to 0\\
&\downarrow&&\downarrow&&&\\
&G&\mbox{\large $=$}&G&&&\\
&\downarrow&&\downarrow&&&\\
&0&&0&&&
\end{array}$$
where $G^{\natural}$ is the push-out of the universal $\G_a$-extension
$G^{\sharp}$ of the semi-abelian variety $G$. The canonical map
$u^{\natural}:L\to G^{\natural}$ is such that the composition
$$L\by{u^{\natural}} G^{\natural}\to
L^{\natural}=\Ext(L[1],\G_a)^{\vee}=\Hom (L,\G_a)^{\vee}$$
is the natural evaluation map.

In particular we get the following extension
$$0\to \Ext (M,\G_a)^{\vee}\to G^{\natural}\to G\to 0$$
of group schemes. The De Rham realization of
$M$ is then defined as
$$T_{DR}(M)\df\Lie G^{\natural},$$
with the Hodge-De Rham filtration given by
$$F^0T_{DR}(M)\df \ker (\Lie G^{\natural}\to \Lie G)\cong
\Ext(M,\G_a)^{\vee}.$$
If $k=\C$ then the De Rham realization is compatible
with the Hodge realization, see \cite[\S 10.1.8, \S 10.3.15]{D}.
We also have an exact sequence
\begin{equation}\label{derhamexseq}
0\to T_{DR}(G)\to T_{DR}(M)\to T_{DR}(L[1])\to 0
\end{equation}
which is the sequence of Lie algebras associated to
$$0\to G^{\sharp}\to  G^{\natural}\to  L^{\natural}\to 0.$$
We may also view (\ref{derhamexseq}) as obtained by applying the
functor $T_{DR}$ to (\ref{motexseq}); thus (\ref{derhamexseq}) is the De
Rham version of (\ref{hodgeexseq}) and (\ref{etexseq}).

Let $X$ be a smooth projective variety over $k=\bar k$ of
characteristic 0, and let
$\Pic^{\natural}(X)$ be the group of isomorphism classes of pairs
$(\ccL,\nabla)$ where $\ccL$ is a line bundle on $X$ and $\nabla$ is an
integrable connection on $\ccL$. Then there is the following extension
$$0\to H^0(X,\Omega^1_X) \to \Pic^{\natural}(X)^0\to \Pic^0(X)\to 0$$
where $\Pic^{\natural}(X)^0$ is the the subgroup of those pairs
$(\ccL ,\nabla)$ such that $\ccL\in\Pic^0(X)$. The above extension is the
group of points of the universal $\G_a$-extension of the abelian variety
$\Pic^0_{X/k}$ and $\Lie \Pic^{\natural}(X)^0 = H^1_{DR}(X)(1)$
(\cf \cite[\S 4]{MM}), where the twist (1) indicates that the
indexing of the Hodge filtration is shifted by 1.
In general, for any abelian variety $A$, $A^{\natural}=
\Pic^{\natural}(A^{\vee})^0$, so that $A$ has De Rham realization
$T_{DR}(A) = H^1_{DR}(A)^{\vee}\df H_1^{DR}(A)$.

\subsection{Cartier duals} We now recall briefly the
construction by Deligne \cite[\S 10.2.11--13]{D} of the
dual 1-motive. The definition is motivated by the case of
1-motives over $\C$ where the Hodge realization has a
dual mixed Hodge structure which yields the dual 1-motive.
In fact, if $H$ is a torsion free mixed Hodge structure
of type $\{(0,0), (0,-1), (-1,0), (-1,-1)\}$ such that
$\gr_{-1}^W(H)$ is polarizable then $H^{\vee}=\Hom (H,\Z
(1))$ is again of the same kind; since any 1-motive $M$
over $\C$ corresponds (uniquely up to isomorphism) to such
an $H=T_{Hodge}(M)$ we can just set $T(M^{\vee})=T(M)^{\vee}$
as an implicit ``analytic'' definition for $M^{\vee}$.

In order to give an algebraic description of $M^{\vee}$ the yoga of
biextensions is needed: see \cite{MUB} for the definition of biextension
(\cf \cite[VII, (2.1)]{SGA7}). Let $M=[L\by{u}G]$ be a 1-motive
over a field, \ie
$$ \begin{array}{ccc}
 & L & \\
& \ \ \downarrow u & \\
 1\to \ T& \to\ G\ \to & A\ \to 0
\end{array}$$
where $G$ is an extension of an abelian variety $A$
by a torus $T$. Then $T^{\vee}$ is a lattice and the dual
abelian variety $A^{\vee}$ can be regarded as $\Ext
(A,\G_m)$; there is a canonical homomorphism
$v : T^{\vee}\to A^{\vee}$ by pushing out characters
$\chi : T\to \G_m$ along the given extension $G$ (\cf \cite{OG},
\cite{SEG}).
By construction, the Poincar\'e biextension
$\cP$ of $A\times A^{\vee}$ by $\G_m$ is trivial on $L\times
T^{\vee}$, \ie there is a bihomomorphism $$\psi : L\times
T^{\vee}\to (u\times v)^*\cP$$
Since $\Hom (L,A)\cong\Ext (A^{\vee},L^{\vee})$ the composite
homomorphism $L\by{u}G\onto A$ yields an extension $G^{u}\in \Ext
(A^{\vee},L^{\vee})$  and the Cartier dual $M^{\vee}$ is given
by  $$ \begin{array}{ccc}   &T^{\vee} & \\
& \ \ \downarrow u^{\vee} & \\
1\to \ L^{\vee}& \to\ G^{u}\ \to & A^{\vee}\ \to 0
\end{array}$$
where the lifting $u^{\vee}$ of $v$ from $A^{\vee}$ to $G^{u}$
is determined by the trivialization $\psi$.

The object $(L,T^{\vee},A,A^{\vee},u,v,\psi)$ is then sometimes
called the ``symmetric avatar'' of the 1-motive $M$; the symmetric avatar of
the Cartier dual is $(T^{\vee},L,A^{\vee},A,v,u,\psi^{t})$.

Finally, as is shown by Deligne \cite[\S 10.2]{D}, the Poincar\'e
biextension  yields pairings on realizations
$$T_{\sZ/m}(M)\otimes T_{\sZ/m}(M^{\vee})\to \mu_m$$
and
$$T_{DR}(M)\otimes T_{DR}(M^{\vee})\to k(1)$$
which are compatible, over the complex numbers, with the canonical
pairing induced by the duality between mixed Hodge structures (here
$k(1)$ is a 1-dimensional filtered $k$-vector space with filtration
$F^{-1}k(1)=k(1)$, $F^0k(1)=0$). We therefore can
see any given realization of the Cartier dual as being the appropriate dual
of that realization of the original 1-motive.

\section{Homological Picard 1-motive: $\Pic^-$}\label{hompic}

We first begin by introducing some notation and terminology needed below.
Let $X$ be an equidimensional variety over a field $k$ of characteristic
zero (not necessarily algebraically closed). Let $S\subset X$ be the
singular locus and let $f:\tilde X \to X$ be a resolution of
singularities. We let $\tilde S =f^{-1}(S)$ be the reduced
inverse image. Consider a smooth compactification of $\tilde X$, which we
denote by $\bar X$; let $Y=\bar X - \tilde X$ be the boundary, which we
assume to be a divisor in $\bar X$. Let $\bar{S}$ denote the Zariski
closure of $\tilde{S}$ in $\bar{X}$. We can arrange that the resolution
$\tilde{X}$ and compactification $\bar{X}$ are chosen so that $\bar{X}$ is
projective, and $\bar{S}+Y$ is a reduced normal crossing divisor in
$\bar{X}$; we call such a compactification $\bar{X}$ a {\em good normal
crossing compactification} (or good n.c. compactification) of the
resolution of $X$. For such a compactification to exist, the resolution
$f:\tilde{X}\to X$ must be chosen such that $\tilde{S}$ is a normal
crossing divisor.

\subsection{Relative Picard functor}
Associated to any pair $(V, Z)$ consisting of any $k$-scheme $V$ and a closed
sub-scheme $Z$, we have a natural long exact sequence
\begin{equation}\label{longpic}
\cdots\to H^0(V,\cO_{V}^*)\to H^0(Z,\cO_{Z}^*)\to
\Pic (V, Z)\to \Pic (V)\to \Pic (Z)\to\cdots
\end{equation}
induced by the surjection of Zariski (or \'etale) sheaves
$\G_{m,V}\to i_*\G_{m, Z}$ where $i: Z\into V$;
here
\[\Pic (V, Z) =\HH^1(V,\G_{m,V}\to i_*\G_{m,Z})\]
is the group of isomorphism classes of pairs $(\ccL, \varphi)$ such that
$\ccL$ is a line bundle on $V$ and  $\varphi :\ccL\mid_Z\cong \cO_Z$ is a
trivialization on $Z$ (\cf \cite{SUV}, \cite[\S 2]{RAS}, \cite[\S
8]{BLR}).

Now let $X$ be an equidimensional $k$-variety, and $\bar{X}$ a good normal
crossing compactification of a resolution of $X$, with boundary $Y$.

\begin{lemma}\label{grpic} Let $(\bar X, Y)$ be as above. The
$fpqc$-sheaf associated to the relative Picard functor
$$T\leadsto  \Pic (\bar X\times_k T, Y\times_k T)$$
is representable by a $k$-group scheme which is locally of finite
type over $k$. If $k$ is algebraically closed, its group of $k$-points is
$\Pic(\bar X,Y)$.
\end{lemma}
\begin{proof} See the Appendix~\ref{rep}.
\end{proof}

Now assume $k=\bar{k}$. Let $Y=\cup Y_i$, where $Y_i$ are the (smooth)
irreducible components of $Y$.
\begin{propose}\label{relpic} The sequence (\ref{longpic}) yields a
semi-abelian group scheme $\Pic^0(\bar X, Y)$ over $k= \bar k$, which
can be represented as an extension
\begin{equation}
1\to T(\bar X,Y)\to \Pic^0(\bar X, Y)\to A(\bar X,Y)\to 0
\end{equation}
where:
\begin{description}
\item[{\it (i)}] $\Pic^0(\bar X, Y)$ is the connected
component of the identity of $\Pic (\bar X,Y)$;
 \item[{\it (ii)}] $T(\bar X,Y)$ is the $k$-torus
$$T(\bar X,Y)\df \coker\left((\pi_{\bar{X}})_*\G_{m,\bar{X}} \to
(\pi_Y)_*\G_{m,Y}\right)$$ where
$\pi_{\bar{X}}:\bar{X}\to\Spec k$, $\pi_Y:Y\to\Spec k$ are the structure
morphisms;
\item[{\it (iii)}] $A(\bar X,Y)$ is the abelian variety
$$A(\bar X, Y) \df \ker^0(\Pic^0(\bar X)\to\oplus\Pic^0(Y_i))$$
which is the connected component of the identity of the kernel.
\end{description}
\end{propose}
\begin{proof} Everything follows from Lemma~\ref{grpic} combined with
(\ref{longpic}), by taking the connected components of the identity, once
we know the following.
\begin{equation}\label{kerpic}
\ker^0(\Pic^0(\bar X)\to\Pic^0(Y)) =
\ker^0(\Pic^0(\bar X)\to\oplus\Pic^0(Y_i)).
\end{equation}
Recall that $\ker^0$ denotes the  connected component of the
identity of the kernel. In order to prove (\ref{kerpic}) we
consider the normalization $\pi :\coprod Y_i\to Y$ and the
following commutative diagram
\begin{equation}\label{norpic}
\begin{TriCDV}
{\Pic^0(\bar X)}{\> \alpha >>}{\oplus \Pic^0(Y_i)}
{\SE \beta E E }{\NE E \pi^* E}{\Pic^0(Y)}
\end{TriCDV}
\end{equation}
Now, because of \cite[Expos\'e XII, Prop.2.3]{SGA6} (\cf \cite{BOU}) the
morphism $$\pi^*:\Pic (Y)\to \oplus\Pic (Y_i)$$ is
representable by an affine morphism. Then $\beta (\ker^0\alpha)=0$, since
$\ker^0\alpha$ is an abelian variety. Since we obviously have
$\ker\beta\subset\ker\alpha$, we must have $\ker^0\alpha =\ker^0\beta$ which
is the claimed equality (\ref{kerpic}).
\end{proof}

\subsection{Definition of $\Pic^-$}

Let $X$ be an equidimensional variety over $k=\bar{k}$ of characteristic
0. As before, let $\bar{X}$ be a good, normal crossing compactification of
a resolution $f:\tilde{X}\to X$ of $X$, with boundary divisor
$Y$. Let $D$ be any Weil (or equivalently Cartier) divisor on $\bar X$
such that $\supp (D)\cap Y =\emptyset$, \ie $D\in\Div (\bar X, Y)$; then
$(\cO_{\bar X}(D), 1)$ defines an element  $[D]\in\Pic (\bar X,
Y)$, where $1$ denotes the tautological section of $\cO_{\bar X}(D)$,
trivializing it on $\bar X - D$, and hence also on $Y$. We say that a
divisor $D\in\Div (\bar X, Y)$ is {\em algebraically equivalent  to zero
relative to $Y$}\, if $[D]\in \Pic^0(\bar X, Y)$ and we denote by
$\Div^0(\bar X, Y)\subset \Div (\bar X, Y)$ the subgroup of divisors
algebraically equivalent to zero relative to $Y$.

Let $\bar S$ be the closure of $\tilde S$ in $\bar X$; then $\bar S\cup Y$
has normal crossings as well, since $\bar{X}$ is ``good''.
Recall that $\Div_{\bar S}(\bar X, Y)\subset \Div (\bar X, Y)$ denotes the
group of divisors $D$ on $\bar X$ supported on $\bar S$ such that $\supp
(D)\cap Y=\emptyset$, \ie it is the free abelian group on the compact
irreducible components of $\tilde S$.  We have a
push-forward on Weil divisors $f_*: \Div_{\tilde S}(\tilde X)\to \Div_S
(X)$ and we let  $\Div_{\tilde S/S}(\tilde X)$ be the kernel of
$f_*$. We finally denote by $\Div^0_{\bar S/S}(\bar X, Y)$ the
intersection of  $\Div_{\tilde S/S}(\tilde X)$ with
$\Div_{\bar S}^0(\bar X, Y)$.  Thus  $\Div^0_{\bar S/S}(\bar X, Y)$ is the
group of divisors on $\bar X$ which are linear combinations of compact
divisorial components in $\tilde S$, which have trivial push-forward under
$f$ and which are algebraically equivalent to zero relative to $Y$.

\begin{defn}\label{defminus}
{\rm Let $X$ be an equidimensional variety over $k=\bar k$.
With the hypothesis and notation as above we define the following 1-motive
$$\Pic^-(X) \df [\Div^0_{\bar S/S}(\bar X, Y)\by{u}
\Pic^0(\bar X, Y)]$$
where $u(D) = [D]$. We call $\Pic^-(X)$ the {\it homological
Picard 1-motive\,} of $X$.

For any closed sub-scheme $Z\subset\bar X$ we define the following 1-motive
$$\Pic^+(\bar X - Z, Y)\df [\Div^0_{Z}(\bar X, Y)\by{u} \Pic^0(\bar X, Y)].$$
If $Z$ is the union of all compact components of divisors in $\tilde S$, we
then  remark that $\Pic^-(X)$ is a sub-1-motive of $\Pic^+(\bar X -Z,
Y)$.

If $X$ is an arbitrary $n$-dimensional variety over $k=\bar{k}$, let
$X^{(n)}$ denote the union of its $n$-dimensional irreducible components.
Define
$$\Pic^-(X)\df \Pic^-(X^{(n)}).$$
}\end{defn}

We next show that our definition of $\Pic^-(X)$ is independent of the
choices made, \ie of the resolution $\tilde X$ and compactification $\bar
X$ as above, when $X$ is equidimensional (\cf also
Remark~\ref{remellminus}).

\subsection{Independence of resolutions and compactifications}
For an equidimensional $k$-variety $X$ as above, consider two
resolutions of singularities $f':X'\to X$ and $f'':X''\to X$  of $X$, with
corresponding good compactifications $\bar X'$ and $\bar X''$. We then can
find a third resolution $f:\tilde X\to X$
dominating both $X'$ and $X''$, and choose a compactification $\bar X$ which is
a resolution of the closure of (the isomorphic image of) $X_{\rm reg}=X-S$
in $\bar X'\times \bar X''$, which is also a good normal crossing
compactification of $\tilde X$.

Hence, to prove independence of $\Pic^-(X)$ from the choices made, it
suffices to consider the following situation. Let $f_1:\tilde X_1\to
X$ be a resolution with good normal crossing compactification $\bar X_1$,
and let $f_2:\tilde X_2\to X$ be another one, with good
normal crossing compactification $\bar X_2$, such that we
have a morphism $f:\bar X_2\to\bar X_1$ whose restriction $f: \tilde
X_2\to\tilde X_1$ is a proper morphism of $X$-schemes, necessarily a
birational morphism. Under these conditions, we wish to show that
$\Pic^-(X)$ defined using either $\bar X_1$ or $\bar X_2$ coincide.

Let $Y_i = \bar X_i - \tilde X_i$ for $i = 1, 2$. We then clearly have a
morphism of 1-motives
$$[\Div^0_{\bar S_1/S}(\bar X_1, Y_1)\to
\Pic^0(\bar X_1, Y_1)]\to [\Div^0_{\bar
S_2/S}(\bar X_2, Y_2)\to\Pic^0(\bar X_2, Y_2)]$$
given by pulling back cycles and line bundles. It suffices to prove this
is an isomorphism of 1-motives.

We first claim that there is an isomorphism of semi-abelian varieties
$$\Pic^0(\bar X_1, Y_1)\cong  \Pic^0(\bar X_2, Y_2).$$
In fact we have the following diagram
$$\begin{array}{ccc}
\Pic^0(\bar X_1)&\cong &\Pic^0(\bar X_2)\\
\downarrow & &\downarrow \\
\Pic (Y_1)&\into &\Pic (Y_2)
\end{array}$$
where the bottom arrow is injective since $f_*(\cO_{Y_2})=
\cO_{Y_1}$ (because $Y_1$ is semi-normal, and $Y_2 \onto Y_1$ has connected
fibers). Thus the kernels of the restrictions are the same, and so,
regarding the relative $\Pic^0$ as an extension (by (\ref{longpic})) and
using Proposition~\ref{relpic}, the claim is  clear.

Now we have a splitting of the pullback map
$$\Div^0_{\bar S_1/S}(\bar X_1, Y_1)\by{f^*}\Div^0_{\bar S_2/S}(\bar X_2,
Y_2)$$
using proper push-forward $f_*$ of divisors; we thus have
$$\Div^0_{\bar S_1/S}(\bar X_1, Y_1)\oplus G\cong
\Div^0_{\bar S_2/S}(\bar X_2, Y_2)$$
where $G \df \{D\in \Div^0_{\bar S_2/S}(\bar X_2, Y_2)\mid f_*(D)=0\}$. We
will show that $G=0$. Since $\bar X_2$ and $\bar X_1$ are birational,
$$f_*:\Pic^0(\bar X_2)\by{\cong}\Pic^0(\bar X_1)$$
we note that if $D\in G$, then $D$ is linearly equivalent to zero
on $\bar X_2$. Thus $D = {\rm div} (r)_{\bar X_2}$ where $r$ is a rational
function on $\bar X_2$, and therefore also on $\bar X_1$. But $${\rm div}
(r)_{\bar X_1}=f_*(D)=0,$$ whence $r$ is constant.

\begin{rmk}{\rm We remark that in our definition of $\Pic^-$, we can allow
$f:\tilde  X\to X$ to be a birational proper morphism from a smooth variety,
which is not necessarily a resolution of singularities of $X$. In fact, for any
birational proper morphism $g:\tilde X'\to \tilde X$ between two such smooth
$X$-varieties we can choose compactifications such that $g$ induces a morphism
$\bar g:\bar X'\to \bar X$. By arguing as above we then
see that $\Pic^0 (\bar X, Y)\cong \Pic^0 (\bar X', Y')$ and
$\Div^0_{\bar S/S}(\bar X, Y)\cong \Div^0_{\bar S'/S}(\bar X', Y')$.}
\end{rmk}

\subsection{Hodge realization of $\Pic^-$}
In order to deal with the Hodge realization of $\Pic^-$ the following
results are needed.

\begin{lemma}\label{piccus}
Let $\bar X$ and $Y$ be as above, with $k=\C$. We then have the
following properties of (the group of $\C$-points of) $\Pic(\bar X,Y)$.
\begin{description}
\item[{\it a)}] There is an exact sequence
$$H^1(\bar X, Y;\Z (1))\to H^1(\bar X,\cO_{\bar X}(-Y))\to
\Pic (\bar X, Y)\by{c\ell}H^2(\bar X, Y;\Z (1)).$$
\item[{\it b)}] There is an isomorphism
$$\Pic^0(\bar X, Y)\cong \ker (\Pic (\bar X,
Y)\by{c\ell} H^2(\bar X, Y;\Z (1))).$$
\item[{\it c)}] There is an isomorphism
$$J^1(\bar X, Y)\df \frac{H^1(\bar X, Y;\C (1))}{F^0 + H^1(\bar X, Y;\Z (1))}
\cong \frac{H^1(\bar X,\cO_{\bar X}(-Y))}{H^1(\bar X, Y;\Z (1))}.$$
\item[{\it d)}] Under the isomorphism (induced by {\it a)--c)})
$$\Pic^0(\bar X,Y)\cong J^1(\bar X, Y)$$
the mapping $D\leadsto [D] =(\cO_{\bar X}(D), 1)$ from
$\Div^0(\bar X, Y)$ to $\Pic^0(\bar X,Y)$ is identified with the
extension class map, for the mixed Hodge structure, determined by the support
of $D$.
\item[{\it e)}] Let $Z$ be a closed sub-scheme $Z\subset\bar X$ such that
$Z\cap Y=\emptyset$. Then
$$T_{Hodge}(\Pic^+(\bar X - Z, Y))\cong H^1(\bar X - Z, Y;\Z (1)).$$
\end{description} \end{lemma}
\begin{proof} We first claim that
$$\frac{H ^1(\bar X,Y;\C(1))}{F^0}\cong H^1(\bar X, \cO_{\bar X}(-Y)).$$
To see this, we consider the twisted log De Rham complex
$\Omega^{\d}_{\bar X}(\log (Y))(-Y)$.
It is well-known (see \cite{ST}, page~4 for a quick proof) that its
hypercohomology groups are the relative cohomology groups
$H^*(\bar X,Y;\C)$, the Hodge-De Rham filtration on the relative
cohomology is given by the subcomplexes $\Omega^{\d\geq i}_{\bar
X}(\log (Y))(-Y)$, and the corresponding hypercohomology spectral
sequence degenerates at $E_1$. We then have
\begin{equation}\label{filt}
\frac{H^k(\bar X,Y,\C)}{F^i}\cong\HH^k(\bar X,
\Omega^{\d<i}_{\bar X}(\log (Y))(-Y)).
\end{equation}
The claimed isomorphism is obtained from (\ref{filt}) for $k=i=1$.

Let $\cO^*_{(\bar X_{\rm an} ,Y_{\rm an})}$ be the sheaf on $\bar X_{\rm an}$
given by the kernel of $\cO^*_{\bar X_{\rm an}}\to i_*\cO^*_{Y_{\rm an}}$ where
$i: Y\into \bar X$ is the inclusion. We have
$$\Pic (\bar X, Y)\cong H^1(\bar X_{\rm an}, \cO^*_{(\bar X_{\rm an},
Y_{\rm an})})$$ because of (\ref{longpic}) and GAGA.
We have that $\cO_{\bar{X}_{\rm an}}(-Y_{\rm an})$ is the kernel of
$\cO_{\bar{X}_{\rm an}}\to
i_*\cO_{Y_{\rm an}}$; therefore, by the exponential sequences on $\bar
X_{\rm an}$ and $Y_{\rm an}$, since $i_*$ is an exact functor, we get the
following induced relative exponential exact sequence of sheaves on
$\bar X_{\rm an}$
\begin{equation}\label{exprel}
0\to j_!(\Z(1))\to \cO_{\bar{X}_{\rm an}}(-Y_{\rm an})\to \cO^*_{(\bar
X_{\rm an} ,Y_{\rm an})} \to 0
\end{equation}
where $j_!$ is the extension by zero functor along
$j: \bar X_{\rm an} - Y_{\rm an}\into \bar X_{\rm an}$.

We then get the following exact sequence of cohomology groups
\begin{equation}\label{relsupcl}
\cdots \to H^1(\bar X, Y;\Z(1))\to H^1(\bar X, \cO_{\bar X}(-Y))\to \Pic
(\bar X, Y)\by{c\ell} H^2(\bar X, Y;\Z(1))\to \cdots
\end{equation}
The exact sequence in {\it a)}\, is then obtained.
Since $H^2(\bar X, Y;\Z(1))$ is finitely generated and $\ker c\ell$ is
divisible, we get {\it b)}. From (\ref{filt}) we then get {\it c)}.

Part {\it d)}\, is well known if $Y=\emptyset$ (\eg see \cite{C}). In
order to show part {\it d)}\, in general, we can proceed  as follows.
By considering relative Deligne-Beilinson cohomology
$H^*_{\cD}(\bar X,Y;\Z (*))$ we get a canonical cycle class map
$$c_1:\Pic (\bar X, Y)\to H^2_{\cD}(\bar X, Y;\Z (1))$$
Moreover, $c_1$ is an isomorphism, fitting into the following commutative
diagram with exact rows
$$ \begin{array}{ccccccc}
\Pic^0(\bar X, Y)&\into &\Pic (\bar X, Y)&\by{c\ell}& H^2(\bar X, Y;\Z(1))
&\to &
H^2(\bar X, \cO_{\bar X}(-Y)) \\
\cong\downarrow\quad & & \cong\downarrow\;c_1 & &\veq &
&\cong\downarrow\quad\\
J^1(\bar X, Y)&\into &H^2_{\cD}(\bar X, Y;\Z (1))&\to & H^2(\bar X, Y;\Z(1))
&\to & H^2(\bar X,Y;\C (1))/F^0
\end{array}$$
obtained from (\ref{filt}) and (\ref{relsupcl}).
For any closed sub-scheme $Z \subset \bar X$ with $Z\cap Y=\emptyset$, we
then have
the following commutative diagram of cohomology groups having exact rows
and columns
$$\begin{array}{ccccccc}
&&&&&&0\\
&&&&&&\downarrow\\
&&&&&&\Pic^0(\bar X, Y)\\
&&&&&&\downarrow\\
&&&&H^2_{\cD, Z}(\bar X, Y;\Z (1))&\to&H^2_{\cD}(\bar X, Y;\Z (1))\\
&&&&\downarrow&&\downarrow\\
H^1(\bar X,Y;\Z (1))&\to & H^1(\bar X -Z,Y;\Z (1)) &\to &H^2_Z(\bar X,Y;\Z
(1))&\to & H^2(\bar X,Y;\Z (1))\\
\downarrow&&\downarrow&&\downarrow&&\\
{\displaystyle\frac{H^1(\bar X,Y; ;\sC (1))}{F^0}}&\to &
{\displaystyle\frac{H^1(\bar X -Z,Y;\sC (1))}{F^0}}
&\to &
{\displaystyle\frac{H^2_Z(\bar X,Y;\sC (1))}{F^0}}& &\\
\downarrow&&&&&&\\
J^1(\bar X,Y)&&&&&&\\
\downarrow&&&&&&\\
0&&&&&&
\end{array}$$
Here $H^*_{\cD, Z}(\bar X, Y;\Z (\cdot))\cong H^*_{\cD, Z}(\bar X;\Z
(\cdot))$ is
the (relative) Deligne-Beilinson cohomology of $(\bar X, Y)$ with support
in $Z$.

Let $Z$ be the support of a divisor $D\in\Div(\bar X, Y)$, \ie  $Z\cap Y
=\emptyset$. We then have that $$H^*_Z(\bar X, Y;\Z(1))\cong H^*_Z(\bar
X,\Z(1)).$$ In particular: $H^1_Z(\bar X, Y;\Z(1))= 0$ and $H^2_Z(\bar X,Y;
\Z (1))$ is purely of type $(0,0)$; in fact, we have an isomorphism
$$H^2_{\cD,Z}(\bar X, Y;\Z (1))\cong H^2_Z(\bar X,Y;\Z
(1))\cong \Div_Z(\bar X, Y).$$
The claim {\it d)}\, then follows from a diagram chase in the diagram
above, using a general homological lemma \cite[Lemma~2.8]{BSL}.

Part {\it e)}\, then follows from the diagram as well, yielding the following
isomorphism, in the category of 1-motives over $\C$,
$$ \begin{array}{ccc}
\Div_{Z}^0(\bar X, Y) & \longby{e}& J^1(\bar X,Y)\\
\parallel & &\quad\downarrow\cong\\
\Div_{Z}^0(\bar X, Y)&\longby{u}&\Pic^0(\bar X,Y)
\end{array}$$
where $e$ denotes the extension class map determined by $H^1(\bar X - Z, Y;
\Z (1))$, regarded as an extension of mixed Hodge structures.
\end{proof}

For the following duality result we refer to the book of Spanier \cite{SP},
giving a  proof in the topological setting. In order to deduce such a duality
statement for different cohomology theories, as well as compatibilities
between them,  we are going to give a proof in the  Grothendieck-Verdier
duality style.
\begin{lemma}\label{duality} Let $M$ be a compact smooth $n$-dimensional
$\C$-variety. Let $A+B$ be a reduced normal crossing divisor
in $M$ such that $A\cap B =\emptyset$. Then there is a duality
isomorphism  $$H^r(M-A,B;\Z)\cong H_{2n-r}(M-B,A;\Z(-n))$$ in the category
of mixed Hodge structures. Moreover  {\it i)}\, this isomorphism is
functorial, \ie if $A'\subset A$, $B\subset B'$ and
$A'+B'$ is also a normal crossing divisor such that $A'\cap
B'=\emptyset$ then the following diagram
$$\begin{array}{ccc}
H^r(M-A',B';\Z)&\to & H^r(M-A,B;\Z)\\
\cong\downarrow\quad & &\quad\downarrow\cong\\
H_{2n-r}(M-B',A';\Z(-n))&\to & H_{2n-r}(M-B,A;\Z(-n))
\end{array}$$
commutes, in the category of mixed Hodge structures; finally {\it ii)}\, this
duality isomorphism is compatible with
the Poincar\'e-Lefschetz duality, \ie if $B=B'$ and $A'\subset A$ as above
then the following diagram, whose rows are long exact sequences, commutes
$$\hspace{-1cm}\begin{array}{rcl}
\cdots\to\quad H^r(M-A',B;\Z)\quad&\to \quad H^r(M-A,B;\Z) \quad\to &
H^{r+1}(M-A',M-A;\Z)\to\cdots\\
\cong\downarrow\quad\quad \quad &\cong\downarrow\quad  &
\quad\quad\quad\quad\cong\downarrow\\
\cdots\to H_{2n-r}(M-B,A';\Z(-n))&\to
H_{2n-r}(M-B,A;\Z(-n))\to &H_{2n-r-1}(A,A';\Z(-n))\to\cdots
\end{array}$$
in the category of mixed Hodge structures.
 \end{lemma}\begin{proof} Let $V = M - (A\cup B)$,
$V_A = M-A$ and $V_B = M - B$ be the corresponding open subsets; we
have a diagram \begin{equation}
\begin{array}{ccc}\label{square}
V &\by{\beta}& V_B \\
\s{\alpha}\downarrow& &\downarrow\s{\gamma}\\
V_A&\by{\delta}&M
\end{array}
\end{equation}
We let ${}^Ai:A\into V_B$, ${}^Bi:B\into V_A$ denote the closed imbeddings.
Let $\pi: M\to k$ be the structure
morphism. Because of the canonical exact sequence
\begin{equation}\label{exrel}
0\to\alpha_!\Z_V\to\Z_{V_A}\to {}^Bi_*\Z_B\to 0
\end{equation}
of sheaves on $V_A$ we have
$$\begin{array}{rl}
H^r(V_A,B)&\cong \Hom (\Z_{V_A}, \alpha_!\Z_V [r])\\
&\cong \Hom (\Z_{M}, \delta_*\alpha_!\Z_V [r])
\end{array}$$
where the $\Hom$ is taken in the derived category. Thus
$$\RHom (\Z_M,\delta_*\alpha_!\Z_V)\cong \RG
(M,\delta_*\alpha_!\Z_V)$$ computes the singular cohomology
of the pair $(M-A,B)$. Now we have
$$H_{2n-r}(M-B,A)/({\rm torsion})\cong\Hom (H^{2n-r}
(M-B,A),\Z)$$ as mixed Hodge structures. Similarly the
complex $\RG (M,\gamma_*\beta_!\Z_V)$ computes the cohomology
of the pair $(M-B,A)$ and we have
$$\begin{array}{rl}
\RHom (\RG (M,\gamma_*\beta_!\Z_V)[2n],\Z (-n))&\cong \RHom (
{\rm R}\pi_!(\gamma_*\beta_!\Z_V)[2n],\Z(-n))\\
&\cong \RHom (\gamma_*\beta_!\Z_V,\pi^!\Z[-2n](-n))\\
&\cong \RHom (\delta_!\alpha_*\Z_V,\Z_M)\\
&\cong\RG (M,\srhom (\delta_!\alpha_*\Z_V,\Z_M))
\end{array}$$
by using Grothendieck--Verdier duality, \ie ${\rm R}\pi_!$ is
left adjoint to $\pi^!$, where the dualizing complex $\omega_{M}\df
\pi^!\Z$ is given by $\pi^!\Z\cong\Z_M[2n](n)$, and
the obvious equality $\gamma_*\beta_! =\delta_!\alpha_*$. Now we
can argue that
$$\begin{array}{rl}
\srhom (\delta_!\alpha_*\Z_V,\Z_M) &\cong \delta_*\srhom
(\alpha_*\Z_V,\delta^!\Z_M)\\&\cong \delta_*\srhom
(\alpha_*\Z_V,\Z_{V_A})\\&\cong\delta_*\alpha_!\Z_V
\end{array}$$
where the last equality is given by the following isomorphism
\begin{equation}\label{biduality}
\alpha_!\Z_V\by{\cong}\srhom (\alpha_*\Z_V,\Z_{V_A})
\end{equation}
The isomorphism (\ref{biduality}) can be obtained from biduality for
constructible sheaves. In fact, let $\omega_{V_A}$ be the dualizing sheaf;
since $V_A$ is smooth $\omega_{V_A}[-2n](-n)\cong \Z_{V_A}$, therefore, by
biduality,
the formula (\ref{biduality}) is equivalent to
$$\srhom (\alpha_!\Z_V,\Z_{V_A})\cong \alpha_*\Z_V$$
which is clear since
$$\begin{array}{rl}
\srhom (\alpha_!\Z_V,\omega_{V_A}[-2n](-n)) &\cong \alpha_*\srhom
(\Z_V,\alpha^!\omega_{V_A}[-2n](-n))\\
&\cong \alpha_*\srhom (\Z_V,\omega_{V}[-2n](-n))\\
&\cong \alpha_*\omega_{V}[-2n](-n)\\
&\cong \alpha_*\Z_{V}
\end{array}$$
where we have used that $\alpha^!\omega_{V_A}$ is the dualizing sheaf on
$V$.

Summarizing, we have obtained the following isomorphism
$$\RHom (\RG (M,\gamma_*\beta_!\Z_V)[2n],\Z(-n))\cong \RG (M,
\delta_*\alpha_!\Z_V)$$ yielding the claimed duality
isomorphism of groups.

In order to show the compatibility of the above with
the mixed Hodge structures we consider the following
induced pairing in the derived category
\begin{equation}\label{pairing}
\gamma_*\beta_!\Z_V[2n](n)\stackrel{L}{\otimes}\delta_*\alpha_!
\Z_V \to\pi^!\Z
\end{equation}
This pairing is unique (up to a unique integer multiple); indeed,
we have
$$\begin{array}{l} \Hom
(\gamma_*\beta_!\Z_V[2n](n)\stackrel{L}{\otimes}\delta_*\alpha_!
\Z_V, \pi^!\Z)\\
\cong \Hom
(\gamma_*\beta_!\Z_V[2n](n)\stackrel{L}{\otimes}\srhom
(\delta_!\alpha_*\Z_V[2n](n), \pi^!\Z),\pi^!\Z)\\ \cong \Hom
(\gamma_*\beta_!\Z_V[2n](n),\srhom (\srhom
(\delta_!\alpha_*\Z_V[2n](n),\pi^!\Z),\pi^!\Z))\\ \cong \Hom
(\gamma_*\beta_!\Z_V[2n](n), \delta_!\alpha_*\Z_V[2n](n))\\ \cong
\Hom (\delta_!\alpha_*\Z_V[2n](n), \delta_!\alpha_*\Z_V[2n](n))\\
\cong \Hom (\Z_V(n), \Z_V(n))\\
\cong \Z
 \end{array}$$
where we have used the formula (\ref{biduality}), biduality
for the constructible sheaf $\delta_!\alpha_*\Z_V$ and the
standard formalism of derived categories. The same arguments
apply to the constant sheaves $\Q$ or $\C$.

By Saito's theory of mixed Hodge modules \cite{SAI}, \cite{SAII}, all of the
above constructions and isomorphisms can (after $\otimes \Q$) be ``lifted'' in
a natural way to the derived category of mixed Hodge modules. In particular,
we see that our duality isomorphism is compatible with the mixed Hodge
structures
as claimed.

We leave to the reader the analogous proofs of the assertions about
functoriality, and compatibility with Poincar\'e-Lefschetz duality.
\end{proof}

\begin{rmk}{\rm We remark that, for the truth of the Lemma~\ref{duality}
the assumption that $A\cap B = \emptyset$ is not really needed: it suffices
to assume that $A + B$ is a reduced normal crossing divisor on $M$, but the
proof in this case is a bit more involved.}
\end{rmk}

\begin{rmk}\label{remduality}
{\rm Let $\Omega_M^{\d}(\log (N))(-D)$ be the $\log$ De Rham
complex with terms
$$\Omega_M^{i}(\log (N))\tensor_{\cO_M} \cO_M(-D)$$
for $D$ any Weil divisor on $M$ and $N$ a reduced  normal crossing divisor
in $M$ which contains $\supp (D)$. Let $j:M-N\into M$ be the inclusion; we
then have
a quasi-isomorphism
$$j_!\C\by{\cong}\Omega_M^{\d}(\log (N))(-N).$$
In the notation of Lemma~\ref{duality} we remark that the following pairing
$$\begin{array}{c}
\Omega_M^{\cdot}(\log
(A+B))(-A)\otimes_{\sC}\Omega_M^{\cdot} (\log (A+B))(-B)\\
\downarrow \\
\Omega_M^{\cdot}(\log (A+B))(-A-B)\\
\downarrow\\
\Omega_M^{\cdot}
\end{array}$$
can be identified with the sheaf theoretic pairing given by (\ref{pairing})
(up to a unique scalar)
$$\gamma_*\beta_!\C_V\stackrel{L}{\otimes}
\delta_*\alpha_!\C_V\longto j_!\C_V\longto \C_M\; .$$
This is valid in the filtered derived category since (\cf \cite{ST})
we obtain the Hodge-De Rham filtrations by truncation of the above log De
Rham complexes; this gives a ``direct'' proof that the duality isomorphism
in the Lemma~\ref{duality} is compatible, after $\otimes \C$, with the
respective Hodge filtrations.

As remarked by Deligne \cite{DE}, the duality isomorphism is, after
$\otimes \Q$, also compatible with weight filtrations, because this is true
$\otimes \Q_{\ell}$, by comparison with \'etale cohomology, and the Weil
conjectures.

For another related compatibility argument using mixed Hodge complexes
we refer to \cite{ST}.}
\end{rmk}

We then have the following key result.
\begin{thm}\label{minus}
Let $X$ be a $\C$-variety of dimension $n$. Then
$$T_{Hodge}(\Pic^-(X))\cong H_{2n-1}(X,\Z(1-n)).$$
\end{thm}
\begin{proof} We first make a reduction to the case when $X$ is
equidimensional. Let $X^{(n)}\subset X$ be the union of the
$n$-dimensional irreducible components of $X$. Then by definition,
$\Pic^-(X)=\Pic^-(X^{(n)})$. On the other hand, the natural map
$H_{2n-1}(X^{(n)},\Z(1-n))\to H_{2n-1}(X,\Z(1-n))$ is an isomorphism of
mixed Hodge structures.

Now for equidimensional $X$, let $f:\tilde{X}\to X$ be a resolution, with
a good normal crossing compactification $\bar X$ with boundary $Y$. As
before, let $S$ be the singular locus of $X$, $\tilde{S}=f^{-1}(S)$, and
$\bar{S}\subset\bar{X}$ the Zariski closure of $\tilde{S}$.

Associated to the cartesian square
$$\begin{array}{ccc}
\tilde S &\into & \tilde X\\
\s{\alpha}\downarrow& &\downarrow\s{f}\\
S&\into &X \end{array}$$
there is a Mayer--Vietoris long exact sequence of mixed Hodge
structures on singular homology yielding the following extension
\begin{equation}\label{basicext}
0\to H_{2n-1}(\tilde X, \Z(1-n))\to H_{2n-1}(X, \Z(1-n))\to
L_X\to 0
\end{equation}
where
$$L_X = \ker H_{2n-2}(\tilde S, \Z(1-n))\to H_{2n-2}(\tilde X,
\Z(1-n))\oplus H_{2n-2}(S, \Z(1-n)).$$

Now we claim:
\begin{description}
\item[{\it (i)}]  $H_{2n-2}(\tilde S, \Z(1-n))\cong \Div_{\bar
S}(\bar X, Y)$,
\item[{\it (ii)}] $f_*:H_{2n-2}(\tilde S, \Z(1-n))\to
H_{2n-2}(S, \Z(1-n))$ is the proper push-forward of algebraic
cycles and
\item[{\it (iii)}] $H_{2n-i}(\tilde X, \Z(1-n))\cong H^i(\bar
X,Y)(1)$ as mixed Hodge structures.
\end{description}
In fact $H_{2n-2}(\tilde S, \Z(1-n))$ is the free
abelian group generated by the compact irreducible $(n-1)$-dimensional
components of $\tilde S$, and $H_{2n-2}(S, \Z(1-n))$ has a similar description.
Thus {\it (i) -- (ii)}\, are clear and {\it
(iii)}\, follows from Lemma~\ref{duality} because $\tilde X =
\bar X -Y$. Moreover we have that the mapping
$$H_{2n-2}(\tilde S, \Z(1-n))\longto H_{2n-2}(\tilde X,
\Z(1-n))$$
induced by the inclusion $\tilde S\into\tilde X$, is just the
cycle map relative to $Y$, \ie the following diagram
$$\begin{array}{ccccc}
\Div_{\bar S}(\bar X, Y) &\by{\cong} &
H_{2n-2}(\tilde S, \Z(1-n))& \longto & H_{2n-2}(\tilde X, \Z(1-n))\\
\downarrow & &\downarrow & &\quad\downarrow\cong\\
\Pic (\bar X, Y) &\mbox{\large =} & \Pic (\bar X, Y)&\longby{c\ell}&
H^2(\bar X, Y)
\end{array}$$
commutes.

Since, by definition, the kernel of $c\ell$ is $\Pic^0(\bar X, Y)$
(\cf Lemma~\ref{relpic}), the lattice $L_X$ is canonically
isomorphic to $\Div_{\bar S/S}^0(\bar X, Y)$. Moreover the exact
sequence (\ref{basicext}) modulo torsion is canonically
isomorphic to the following exact sequence
\begin{equation}
0\to H^1(\bar X, Y,\Z(1))/{\rm (torsion)}\to
H_{2n-1}(X,\Z(1-n))/{\rm (torsion)}\to L_X\to 0
\end{equation}
in the category of torsion free mixed Hodge structures. But
$H^1(\bar X,Y)(1)$ is torsion-free, by the universal coefficient
theorem in topology; hence so is $H_{2n-1}(X,\Z(1-n))$.
The Hodge structure on $L_X$ is pure of weight zero and type $(0,0)$; we
then have
$$W_{-i}(H_{2n-1}(X,\Z(1-n)))=
W_{-i}(H^1(\bar X, Y,\Z(1))),
\mbox{\hspace*{1cm}$i\geq 1$.}$$
We also have the following extension of mixed Hodge structures
\begin{equation}
0\to \frac{H^0(Y,\Z)}{\im H^0(\bar X,\Z)}\otimes \Z(1)\to
H^1(\bar X, Y,\Z(1))\to \ker (H^1(\bar X,\Z(1))\to
H^1(Y,\Z(1)))\to 0
\end{equation}
Thus the weight filtration of $H_{2n-1}(X,\Z(1-n))$ admits the following
description. Let
$$r=\rank \frac{H^0(Y,\Z)}{\im H^0(\bar X,\Z)}$$
then
$$W_{-2}(H_{2n-1}(X,\Z(1-n))) \cong
\Z(1)^{\oplus r}$$
and
$$W_{-1}(H_{2n-1}(X,\Z(1-n))) \cong
H^1(\bar X, Y,\Z(1)).$$
Since $H^1(\bar X,\Z(1))$ is pure of weight $-1$, we have
$$\ker (H^1(\bar X,\Z(1))\to H^1(Y,\Z(1))) =
\ker (H^1(\bar X,\Z(1))\to \oplus H^1(Y_i,\Z(1)))$$
whence
$$\gr_{-1}^{W} H_{2n-1}(X,\Z(1-n))\cong
\ker (H^1(\bar X,\Z(1))\to \oplus H^1(Y_i,\Z(1)))$$
and
$$\gr_{0}^{W} H_{2n-1}(X,\Z(1-n))\cong
L_X$$

Thus the 1-motive associated (by Deligne) to
$H_{2n-1}(X,\Z(1-n))$ is given by the following
$$ \begin{array}{ccc}
 & L_X & \\
& \ \ \downarrow \s{e} & \\ 1\to \ (\C^*)^{\oplus r}& \to\
J^1(\bar X,Y)\ \to & \ker^0 (J^1(\bar X)\to \oplus J^1(Y_i)) \to 0
\end{array}$$
($J^1(\bar X,Y)$ was defined in Lemma~\ref{piccus}(c); $J^1(\bar X)$,
$J^1(Y_i)$ are similarly defined).
Since
$$T_{Hodge}(L_X\to J^1(\bar X,Y))\cong H_{2n-1}(X,\Z(1-n))$$
by Deligne's construction, we are reduced to showing that
\begin{equation}
[L_X\to J^1(\bar X,Y)]\cong \Pic^-(X)
\end{equation}
in the category of 1-motives over $\C$.

By Lemma~\ref{piccus} (\cf Proposition~\ref{relpic}, where $T(\bar X,
Y)(\C)=(\C^*)^{\oplus r}$ and $A (\bar X, Y)$ is the above abelian
variety) we have that
$$\Pic^0(\bar X, Y)\cong J^1(\bar X, Y).$$
According to our definition of $\Pic^-(X)$ we are left to check that the
following
\begin{equation}\label{abel}
\begin{array}{ccc}
 L_X & \longby{e}& J^1(\bar X,Y)\\ \parallel & &\quad\downarrow\cong\\
\Div_{\bar S/S}^0(\bar X, Y)&\longby{u}&\Pic^0(\bar X,Y)
\end{array}
\end{equation}
commutes. We will deduce this from Lemma~\ref{piccus}(d)

The Mayer-Vietoris exact sequence yielding (\ref{basicext}) is
given by the following commutative diagram of mixed Hodge structures
$$\begin{array}{ccccccc}
H_{2n-1}(\tilde X)&\into& H_{2n-1}(\tilde
X,\tilde S)&\to& H_{2n-2}(\tilde S) &\to & H_{2n-2}(\tilde X)\\
\downarrow & &\cong\downarrow\quad & &\downarrow & &\downarrow\\
H_{2n-1}(X)&\into& H_{2n-1}(X,S)&\to& H_{2n-2}(S) &\to & H_{2n-2}(X)
\end{array}$$
which yields the following diagram of mixed Hodge structures
$$\begin{array}{ccccc} 0\to H_{2n-1}(\tilde X,\Z(1-n))&\to&
H_{2n-1}(X,\Z(1-n))&\to &\Div_{\bar S/S}^0(\bar X, Y)\to 0\\ \parallel &
&\downarrow & &\downarrow\\ 0\to H_{2n-1}(\tilde X,\Z(1-n))&\to&
H_{2n-1}(\tilde X,\tilde S)  &\to& \Div_{\bar S}^0(\bar X, Y)\to 0
\end{array}$$ where $$\ker (H_{2n-2}(\tilde S) \to H_{2n-2}(\tilde X)) \cong
\Div_{\bar S}^0(\bar X, Y).$$

Let $D\in L_X= \Div_{\bar S/S}^0(\bar X, Y)$; then $Z=\supp (D)$ is a
closed subset of $\bar X$ such that $Z\cap Y=\emptyset$, and $D$ is
homologically equivalent to zero relative to $Y$; we let
$$\Z[Z]^0\df \ker (\Div_{Z}^0(\bar X, Y)\to H^2(\bar X,Y;\Z(1))).$$
We have the following diagram of torsion-free mixed Hodge structures
\begin{equation}\label{extclass}
\begin{array}{ccccc} 0\to H^{1}(\bar X,Y;\Z(1)) &\to&
H^{1}(\bar X - Z,Y;\Z(1))&\to &\Z[Z]^0\to 0\\
\hspace*{0.5cm}\cong\downarrow\quad & &\downarrow &
&\downarrow\hspace*{0.5cm}\\
0\to H_{2n-1}(\tilde X,\Z(1-n))
&\to& H_{2n-1}(\tilde X,\tilde S; \Z(1-n))
&\to& \Div_{\bar S}^0(\bar X, Y)\to 0
\end{array}
\end{equation}
where the middle vertical mapping is obtained as follows. By
Lemma~\ref{duality} we have
$$H^{1}(\bar X - Z,Y;\Z(1))\cong H_{2n-1}(\bar X - Y,Z;\Z(1-n)).$$
Since $\tilde X = \bar X - Y$ and $Z\into \tilde S$ we have the following
canonical map of mixed Hodge structures
$$H_{2n-1}(\bar X - Y,Z;\Z(1-n))\into
H_{2n-1}(\tilde X,\tilde S;\Z(1-n)).$$
The claimed map is obtained by composition of the duality isomorphism and
the latter inclusion. Thus the diagram (\ref{extclass}) commutes by the
functoriality assertion in Lemma~\ref{duality}.

By diagram chase on (\ref{extclass}) one can then see that the image of $D$
under the mapping $$e :L_X\to J^1(\bar X, Y)$$ is the image of $[D]$ under
the extension class map $$\Z[Z]^0\to J^1(\bar X, Y)$$
determined by the top row of (\ref{extclass}). Thus (\ref{abel})  commutes by
Lemma~\ref{piccus} part {\it d)}. The Theorem~\ref{minus} is proved.
\end{proof}

\begin{rmk}{\rm In order to show that (\ref{abel}) commutes, which is the key
point in proving Theorem~\ref{minus}, one can instead choose $Z$ to be the
union of all compact components of $\tilde S$. By excision and duality we
then have
$$\displaylines{
H_{2n-1}(X)\into H_{2n-1}(X, S,\Z(1-n))\cong H_{2n-1}(\tilde X,
\tilde S,\Z(1-n))\cr
\cong H_{2n-1}(\tilde X, Z,\Z(1-n)) \cong  H^{1}(\bar X - Z, Y,\Z(1))\cr
}$$
and, comparing with the Mayer-Vietoris sequence (\ref{basicext}) we have the
following pull-back diagram (obtained as above from Lemma~\ref{duality})
\begin{equation}\label{pbhr}
\begin{array}{ccccccc}
0\to &H^{1}(\bar X, Y,\Z(1))&\to&H^{1}(\bar X - Z, Y,\Z(1))&\to&
\Div_{Z}^0(\bar X,Y)&\to 0 \\
&\veq & &\uparrow & &\uparrow & \\
0\to&H^1(\bar X, Y,\Z(1))&\to & H_{2n-1}(X,\Z(1-n)) &\to &
\Div_{\bar S/S}^0(\bar X,Y)&\to 0
\end{array}
\end{equation}
in the category of mixed Hodge structures. Therefore, the claimed
commutativity of (\ref{abel}), now follows directly from Lemma~\ref{piccus}
part {\it e)}, as $H^{1}(\bar X - Z, Y,\Z(1))$ is the Hodge realization of
$\Pic^+(\bar X- Z, Y)$ and $\Pic^-(X)$ is a sub-1-motive.} \end{rmk}

\subsection{\'Etale realization of $\Pic^-$}
Let $V$ be any $k$-scheme over a field $k$ (of characteristic zero).
For any pair $(V, Z)$ where $Z$ is a closed subscheme of $V$ we
denote by $i:Z\into V$ and $j:V-Z\into V$ the corresponding inclusions. We then
have that $\G_{m,V}\to i_*\G_{m,Z}$ is an epimorphism of \'etale sheaves,
and we let $\G_{m,(V,Z)}$ denote its kernel. Associated to any such  pair
there is an exact sequence
$$0\to j_!(\mu_{n})\to \G_{m,(V,Z)}\by{m} \G_{m,(V,Z)}\to 0$$
induced by multiplication by $m$ on $\G_m$'s and
the snake lemma, where as usual $\mu_m$ denotes the \'etale sheaf of
$m^{\rm th}$ roots of unity (\cf \cite{MI}, \cite{SUV}).  A ``relative
Hilbert's theorem 90'' is clearly available (\cf
\cite[Section~1]{SUV}).

\begin{propose}\label{rel90}
There is an isomorphism
$$H^1_{\et}(V, \G_{m,(V,Z)})\cong \Pic (V, Z)$$
\end{propose}
\begin{proof} From the Leray spectral sequence
along $\varepsilon :V_{\et}\to V_{\rm Zar}$ for the sheaf $\G_{m,(V,Z)}$ we
get a functorial map $\Pic (V, Z)\to H^1_{\et}(V, \G_{m,(V,Z)}).$ We can then
consider the long exact sequence (\ref{longpic}) and compare with the
corresponding sequence of \'etale cohomology groups. Since $H^i_{\et}(V,
\G_{m,V}) \cong H^i(V, \cO_{V}^*) $ and $H^i_{\et}(Z, \G_{m,Z})  \cong
H^i(Z, \cO_{Z}^*)$ for $i=0, 1$ we then get the result.
\end{proof}

\begin{propose} \label{relkum}
We have the following ``relative Kummer sequence''
$$0\to H^0(V, \G_{m,(V,Z)})/m \by{u}
H^1_{\et}(V,j_!(\mu_m))\by{p} \Pic (V, Z)_{m-{\rm tors}}\to 0
$$
where:
\begin{description} \item[-]
$H^1_{\et}(V,j_!(\mu_m))$ can be interpreted as the group of isomorphism
classes of triples $(\ccL, \varphi, \eta)$ given by a line bundle $\ccL$
on $V$, a trivialization $\varphi :\ccL\mid_Z\cong \cO_Z$ and an
isomorphism $\eta: \cO_{V}\by{\cong}\ccL^{\otimes m}$ which is compatible
with $\varphi^{\otimes m}$, \ie such that
$\eta\mid_Z=\varphi^{\otimes m}$;
\item[-] $H^0(V, \G_{m,(V,Z)})$ is the subgroup of those elements in
 $H^0(V,\cO_{V}^*)$ yielding 1 in $H^0(Z,\cO_{Z}^*)$;
\item[-] the map $u$ is defined by taking a unit $a\in
H^0(V,\G_{m,(V,Z)})$ to $(\cO_V,1,a^{-1})$;
\item[-] the map $p$ takes a triple $(\ccL,
\varphi, \eta)$ to the pair $(\ccL, \varphi)$ which is an $m$-torsion
element of $\Pic (V,Z)$.
\end{description}
\end{propose}
\begin{proof} The description above can be easily obtained by modifying the
original argument for absolute $\Pic$ (\cf \cite[III.4]{MI}). \end{proof}

We can regard $H^*_{\et}(V,j_!(\mu_m))$ as ``relative \'etale cohomology''
groups of the pair $(V,Z)$ for which we adopt the notation
$H^*_{\et}(V,Z;\mu_m)$. \'Etale homology groups
$H_*^{\et}(V,\mu_m^{\otimes t})$ are defined, for an $n$-dimensional
$V$ and integer $t$, to be the cohomology groups of the following (dual)
complex
$$\RHom(\RG(V,\mu_m^{\otimes -t}),\mu_m^{\otimes (-n)}[-2n])$$
in the (twisted) derived category of \'etale sheaves of $\Z/m$-modules;
these homology groups, in general, are not the same as Borel-Moore \'etale
homology groups defined by the dualizing complex. We then have the following
result (\cf Lemma~\ref{duality}).

\begin{lemma}\label{etduality} Let $V$ be an $n$-dimensional
proper smooth variety over $k=\bar k$ of characteristic 0.
Let $A+B$ be a normal crossing divisor in $V$ such that $A\cap B
=\emptyset$. Then there is a functorial duality isomorphism
$$H^r_{\et}(V-A,B;\mu_m)\cong H_{2n-r}^{\et}(V-B,A;\mu_m^{\otimes
(1-n)})$$
which is compatible with Poincar\'e-Lefschetz duality.
\end{lemma}
\begin{proof}
The same proof of Lemma~\ref{duality} applies here to the \'etale sheaf
$\mu_m$.\end{proof}

Suppose that $\bar X$ is a good n.c. compactification, with boundary $Y$,
of a resolution $\tilde{X}$ of an equidimensional $n$-dimensional
$k$-variety $X$, where $k=\bar k$, char.~$k=0$, and let $S,\tilde{S},\bar
S$ be as before. Let $(D,\ccL)\in \Div_{\bar S/S}(\bar X, Y)\times
\Pic^0(\bar X,Y)$; by definition (see Section~\ref{pre} for details)
$$T_{\sZ/m}(\Pic^-(X)) =\frac{\{(D,\ccL)\mid
\eta_D:\ccL^{\otimes m}\cong\cO_{\bar X}(-D)\}}{\{(mD,\cO_{\bar
X}(-D))\}}.$$
We have a canonical map
$$\rho_m^-: T_{\sZ/m}(\Pic^-(X))\to
H_{2n-1}^{\et}(X,\mu_m^{\otimes (1-n)})$$
defined as follows. Let $D$ be a divisor in $\Div_{\bar S/S}^0(\bar X, Y)$
and let $Z$ be its support. If $(D,\ccL)$ is a pair in
$T_{\sZ/m}(\Pic^-(X))$ then $(\ccL\mid_{\bar X-Z},1,\eta_D\mid_{\bar
X-Z})$
belongs to $H^1_{\et}(\bar X -Z,Y;\mu_m)$ by relative Kummer theory (\ie
the description in Proposition~\ref{relkum}); furthermore, the image of the
triple $(\ccL\mid_{\bar X-Z},1,\eta_D\mid_{\bar X-Z})$
under the boundary map
$$H^1_{\et}(\bar X -Z,Y,\mu_m)\to H^2_{\et,Z}(\bar X,Y,\mu_m)$$
is the class of $D$.

We have the following commutative diagram with exact rows
$$\begin{array}{ccccccc}
0&\to&H^1_{\et}(\bar X,Y,\mu_m)&\to &H^1_{\et}(\bar X -Z,Y,\mu_m)
&\to & H^2_{\et,Z}(\bar X,Y,\mu_m)\\
 & &\cong\ \downarrow\quad & &\downarrow & &\downarrow\\
0&\to&H_{2n-1}^{\et}(\tilde X,\mu_m^{\tensor 1-n})&\to&
H_{2n-1}^{\et}(\tilde X,\tilde S,\mu_m^{\tensor 1-n})&\to&
H_{2n-2}^{\et}(\tilde S,\mu_m^{\tensor 1-n})\\
 & &\downarrow & &\quad \downarrow\ \cong & &\downarrow\\
0&\to&H_{2n-1}^{\et}(X,\mu_m^{\tensor 1-n})&\to&
H_{2n-1}^{\et}(X,S,\mu_m^{\tensor 1-n})&\to&
H_{2n-2}^{\et}(S,\mu_m^{\tensor 1-n})
\end{array}$$
We then can define $\rho_m^-(D,\ccL)$ to be the image of
$(\ccL\mid_{\bar X-Z},1,\eta_D\mid_{\bar X-Z})$ in
$H_{2n-1}^{\et}(X,\mu_m^{\otimes (1-n)})$.

We let
$$H_{2n-1}^{\et}(X,\hat{\Z}(1-n))\df \liminv{m}
H_{2n-1}^{\et}(X,\mu_m^{\otimes (1-n)}).$$
We can show the following.
\begin{thm}\label{ellminus} Let $X$ be a variety of dimension $n$
over an algebraically closed field $k$ of characteristic zero. Then
$$\hat{T}(\Pic^-(X))\cong H_{2n-1}^{\et}(X,\hat{\Z}(1-n)).$$
\end{thm}
\begin{proof} As in the proof of Theorem~\ref{minus}, we reduce
immediately to the case when $X$ is equidimensional. Now we fix a choice
of resolution $f:\tilde{X}\to X$, good compactification $\bar X$, {\it
etc\/}.

By definition, $\Pic^-(X)$ is given by the 1-motive
$[\Div_{\bar S/S}^0(\bar X, Y)\to \Pic^0(\bar  X,Y)]$.
We have the following commutative diagram
$$\begin{array}{ccccccc}
0\to &\hat{T}(\Pic^0(\bar  X,Y))&\to&\hat{T}(\Pic^-(X))&\to&
\hat{T}(\Div_{\bar S/S}^0(\bar X, Y)[1])&\to 0 \\
&\bar\rho_{\et}\ \downarrow\quad & &\quad\downarrow\ \rho_{\et} &
&\quad\quad\downarrow\ \rho_{\et}^0 & \\
0\to& H^1_{\et}(\bar X, Y;\hat{\Z}(1))&\to &
H_{2n-1}^{\et}(X,\hat{\Z}(1-n))&\to & \hat{\Z}^{\oplus r} &\to 0
\end{array}$$
where the bottom row is given by the Mayer-Vietoris sequence for
\'etale homology and the duality Lemma~\ref{etduality} ($r$ is a
certain non-negative integer), and the top exact sequence is given by
(\ref{etexseq}) in Section~\ref{pre}. We get the
mapping $\rho_{\et}$ above by taking limit of $\rho_m^-$, and
$\bar{\rho}_{\et}$ is the induced map.

Note that $\bar{\rho}_{\et}$ may be
viewed as the analogue of $\rho_{\et}$ for the variety $\tilde{X}$.
 It is also easy to see from the definitions that
$\hat{T}(\Div_{\bar S/S}^0(\bar X, Y)[1]) \cong \hat{\Z}^{\oplus r}$ as
well, such
that $\rho_{\et}^0$ is an isomorphism. Granting this, we are left to show
our claim holds true for smooth schemes, \ie that $\bar \rho_{\et}$
is an isomorphism. The latter follows from the fact that the
relative Neron-Severi group of $(\bar X,Y)$ is finitely generated, whence
$\hat{T}(\Pic^0(\bar X, Y))= \liminv{m}\Pic(\bar X, Y)_{m-{\rm tors}}$,
and, by Proposition~\ref{rel90} and the Kummer sequence in
Proposition~\ref{relkum}, we have $\Pic(\bar X, Y)_{m-{\rm
tors}}=H^1_{\et}(\bar X, Y;\mu_m)$, since $H^0(X,\G_{m,(\bar
X,Y)})$ is divisible.
\end{proof}
\begin{rmk}\label{remellminus}
{\rm Theorem~\ref{ellminus} can also be used to show that $\Pic^-$ is
independent of the choices of resolutions and compactifications. In fact,
after Proposition~\ref{fhat}, the induced isomorphism on \'etale
realizations lifts to 1-motives. But, as remarked before, we consider this
proof to be ``not in the spirit of the theory of 1-motives''.}
\end{rmk}

\subsection{De Rham realization of $\Pic^-$}
Let $k$ be a field of characteristic 0. Let $\bar{X}$ be a smooth
$k$-variety, with smooth compactification $\bar X$  and normal crossing
boundary $Y =\bar X -X$. Let $\pi:\tilde{Y}\to Y$ be the normalization,
and $i:Y\into \bar X$ the inclusion. Note that $\tilde{Y}$ is a smooth
proper $k$-variety as well.

If $(\ccL,\nabla)$ is a line bundle on $\bar X$ equipped with an
integrable ($k$-linear) connection, then restriction to $Y$ yields a
connection on $\ccL\mid_Y$ with values in $\pi_*\Omega^1_{\tilde{Y}}$,
\[\tilde{\nabla}:i^*\ccL\to
i^*\ccL\tensor_{\cO_Y}\pi_*\Omega^1_{\tilde{Y}}\]
defined as the composition of the restricted connection $i^*\ccL\to
i^*\ccL\tensor_{\cO_Y}\Omega^1_Y$ with the natural $\cO_Y$-linear map
$i^*\ccL\tensor\Omega^1_Y\to i^*\ccL\tensor\pi_*\Omega^1_{\tilde{Y}}$.
There is also a connection $\tilde{d}:\cO_Y\to \pi_*\Omega^1_{\tilde{Y}}$,
similarly defined using the exterior derivative map $d:\cO_Y\to
\Omega^1_Y$; this is just $\tilde\nabla$ in the case when $\ccL=\cO_{\bar
X}$ and $\nabla=d$.

We will denote by $\Pic^{\natural}(\bar X,Y)$ the group of
isomorphism classes of triples $(\ccL,\nabla,\varphi)$, where
$(\ccL,\nabla)$ is a line bundle on $\bar X$ with an integrable
connection,  and
\[\varphi: (i^*\ccL,\tilde{\nabla})\cong (\cO_Y,\tilde{d})\]
is a trivialization on $Y$ as connections with values in
$\pi_*\Omega^1_{\tilde{Y}}$; equivalently, we have a trivialization of
$i^*\ccL$ such that the induced trivialization of $\pi^*i^*\ccL$ is given
by a flat section, for the induced connection on $\pi^*i^*\ccL$ (in the
standard sense) obtained from $\nabla$.

We can consider the relative $\natural$-Picard functor on
the category of schemes over $k$, which we
denote by $\Pic_{(\bar X,Y)/k}^{\natural}$, and is defined to be
the $fpqc$-sheaf associated to the functor
$$T\longmapsto \Pic^{\natural}
(\bar X\times_k T,Y\times_k T).$$
We clearly have the following commutative square
$$\begin{array}{ccc}
\Pic^{\natural}(\bar X, Y) &\to &\Pic^{\natural}(\bar X) \\
\downarrow & &\downarrow\\
\Pic(\bar X, Y)&\to &\Pic (\bar X)
\end{array}$$
which is functorial as well.

Let $i:Y\into \bar X$ be the inclusion of the normal crossing boundary,
and let $\pi:\tilde{Y}\to Y$ be the normalization. We have an
induced relative dlog map given by the following diagram
\begin{equation}\label{reldlog}
\begin{array}{ccccccc}  0\to & \Omega^1_{\bar X}(\log Y)(-Y)&\to
&\Omega^1_{\bar X}& \to&i_*\pi_*\Omega^1_{\tilde Y}&\to 0\\
&\  \ \uparrow \s{\rm dlog}& &\  \
\uparrow \s{\rm dlog}& & \ \ \uparrow\s{\rm dlog} & \\
 0\to &\cO^*_{(\bar X,Y)}&\to & \cO^*_{\bar X} &\to &i_*\cO_{Y}^*&\to 0
\end{array}\end{equation}

We now have the following.
\begin{propose}\label{rignabla}  Let $(\bar X , Y)$ be any pair as above
over $k=\bar k$, and let $Y_i$ ($i= 1,2,\ldots$) denote the (smooth)
irreducible components of the normal crossing boundary divisor $Y$.
\begin{description}
\item[{\it a)}] There is a functorial isomorphism
$$\Pic^{\natural}(\bar X,Y)\cong
\HH^1(\bar X,\cO^*_{(\bar X , Y)}\by{\rm dlog}
\Omega^1_{\bar X}(\log Y)(-Y)).$$
\item[{\it b)}] There is an extension
$$1\to \frac{H^0(Y,\cO_Y^*)}{H^0(\bar X,\cO_{\bar X}^*)}
\to \Pic^{\natural}(\bar X, Y)^0\to \ker^0(\Pic^{\natural}(\bar
X)^0\to\oplus\Pic^{\natural}(Y_i)^0)\to 0 $$
where $(\Pic^{\natural})^0$ denotes the pull-back of $\Pic^0$ in
$\Pic^{\natural}$.
\item[{\it c)}] The universal $\G_a$-extension of the
semi-abelian variety $\Pic^0(\bar X, Y)$ is given
by the $k$-group scheme $(\Pic^{\natural}_{(\bar X,  Y)/k})^0$,
\ie in the notation of Section~1.4, we have an isomorphism
$$(\Pic^0(\bar X, Y))^{\natural}\cong
(\Pic^{\natural}_{(\bar X, Y)/k})^0.$$
\item[{\it d)}] We have an isomorphism
$$\Lie \Pic^{\natural}(\bar X, Y)^0\cong \HH^1(\bar X,\cO_{\bar
X}(-Y)\to \Omega^1_{\bar X}(\log Y)(-Y)).$$
\end{description}
\end{propose}
\begin{proof} In order to show part {\it a)}\, we consider the
canonical mapping which associates to any line bundle $\ccL$ with an
integrable  connection $\nabla$, trivialized along $Y$ (in the appropriate
sense), the cohomology class of a \v Cech cocycle given by the transition
functions defining $\ccL$ and the induced forms. Since the following
sequence (defined by the obvious maps)
$$0\to H^0(\bar X, \Omega^1_{\bar X}(\log Y)(-Y))\to
\Pic^{\natural}(\bar X, Y)\to
\Pic(\bar X, Y)\to H^1(\bar X, \Omega^1_{\bar X}(\log Y)(-Y))$$
is exact, we get the claimed isomorphism: note that
$$H^0(\bar X, \Omega^1_{\bar X}(\log Y)(-Y))\subset H^0(\bar X,
\Omega^1_{\bar X})$$ consists of closed 1-forms, since char.~$k=0$.

The exact sequence in  {\it b)}\, is obtained by the exact sequence of
complexes given by the columns in (\ref{reldlog}): in fact,
the following equation holds
$$\ker^0(\Pic^{\natural}(\bar X)^0\to \HH^1(Y,\cO_{Y}^*\to
\pi_*\Omega^1_{\tilde Y})) =
\ker^0(\Pic^{\natural}(\bar X)^0\to\oplus\Pic^{\natural}(Y_i)^0)$$
by the Proposition~\ref{relpic}.

—From the above discussion we get the following diagram with exact rows
and columns
$$\begin{array}{ccccccc}
&&&0&&0&\\
&&&\uparrow&&\uparrow&\\
0\to & \frac{H^0(Y,\cO_Y^*)}{H^0(\bar X,\cO_{\bar X}^*)}
&\to &\Pic^0(\bar X,Y)&\to &
\ker^0(\Pic^{0}(\bar X)\to\oplus\Pic^{0}(Y_i))&\to 0\\
&\veq&&\uparrow&&\uparrow&\\
0\to& \frac{H^0(Y,\cO_Y^*)}{H^0(\bar X,\cO_{\bar X}^*)}&\to &
\Pic^{\natural}(\bar X, Y)^0&\to &
\ker^0(\Pic^{\natural}(\bar X)^0\to\oplus\Pic^{\natural}(Y_i)^0)
&\to 0\\ &&&\uparrow&&\uparrow&\\
&&&H^0(\bar X,\Omega_{\bar X}(\log Y)(-Y))&\mbox{\large $=$}&\ker
(H^0(\bar X,\Omega_{\bar X})\to \oplus_i H^0(Y_i, \Omega^1_{Y_i}))&\\
&&&\uparrow&&\uparrow&\\
&&&0&&0&
\end{array}$$
Therefore we see that $\Pic^{\natural}(\bar X, Y)^0$ is the group
of $k$-points of the pull-back of the group scheme
$$\ker^0((\Pic^{\natural}_{\bar
X/k})^0\to\oplus(\Pic^{\natural}_{Y_i/k})^0)$$
The latter is the universal extension of the abelian variety
$$\ker^0(\Pic^{0}(\bar X)\to\oplus\Pic^{0}(Y_i)),$$
therefore $\Pic^{\natural}(\bar X, Y)^0$ is the universal extension of
the semi-abelian variety $\Pic^0(\bar X, Y),$ and
{\it c)}\, is proved.

The part {\it d)}\, is standard (\eg can be obtained in a manner similar to the
corresponding result for the usual Picard functors, by computing
$k[\varepsilon]$-points as in \cite{MM}).
\end{proof}

Let $X$ be equidimensional over $k=\bar k$, and fix a resolution
$f:\tilde{X}\to X$ and good normal crossing compactification $\bar X$ with
boundary $Y$, as usual. Let $Z$ denote the union of all compact
components in $\tilde S$. By our choice of resolution and
compactification, $Z$ has normal crossings, and $Z\cap Y=\emptyset$.
Denote by $Z_j$ ($j= 1,2,\ldots$) its smooth irreducible components.
Recall (Definition~\ref{defminus}) that $Z$ yields a 1-motive defined as
follows
$$\Pic^+(\bar X -Z,Y)\df [\Div^0_Z(\bar X,Y)\longby{u_Z}\Pic^0(\bar
X,Y)].$$
By definition, $\Pic^-(X)$ is a sub-1-motive of $\Pic^+(\bar X -Z,Y)$.

Correspondingly, we define the group $\Pic^{\natural-\log}(\bar X -Z,Y)$ as the
group of isomorphism classes of triples $(\ccL,\nabla^{\log},\varphi)$ where
$\ccL$ is a line bundle on $\bar X$, $\nabla$ is an integrable connection on
$\ccL$ with log  poles along $Z$, \ie
$$\nabla^{\log}:\ccL\to \ccL\otimes \Omega_{\bar X}^{1}(\log Z)$$
and $\varphi: (i^*\ccL,\tilde{\nabla})\cong (\cO_Y, \tilde{d})$ is a
trivialization (note that we are assuming $Z\cap Y =\emptyset$).

There is a lifting
$$u^{\natural}_Z:\Div^0_Z(\bar X,Y)\to  \Pic^{\natural-log}(\bar X -Z,
Y)^0$$
of $u_Z: \Div^0_Z(\bar X,Y)\to\Pic^0(\bar X,Y)$. The lifting
$u^{\natural}_Z$ is obtained from the fact that given a divisor
$D\in\Div_Z(\bar X, Y)$, the line bundle $\cO_{\bar X}(D)$ comes
equipped also with a canonical connection with log poles along
$\supp(D)\subset Z$. The connection is characterized by the property that
the tautological meromorphic section, with divisor $D$, is flat (\ie if
$s$ is this section, $\nabla (s)=0$ defines a connection on the open
complement of $\supp(D)$, which one verifies, by local calculation, has a
unique meromorphic extension with log poles along $\supp(D)$). We then
have the following result.

\begin{lemma}\label{riglog}  Let $\bar X , Y$ and $Z$ be as above.\\[0.2cm]

{\it a)}\, There is an isomorphism
$$\Pic^{\natural-\log}(\bar X,Y)\cong
\HH^1(\bar X,\cO^*_{(\bar X , Y)}\by{\rm dlog}
\Omega^1_{\bar X}(\log (Y+Z))(-Y)).$$

{\it b)}\, There is an extension
$$\hspace{-10pt}
0\to \Pic^{\natural}(\bar X, Y)
\to \Pic^{\natural-\log}(\bar X-Z, Y)\to
\ker\left(\oplus_j H^0(Z_j,\cO_{Z_j})\to H^1(\bar X,\Omega^1_{\bar X}
(\log (Y+Z))(-Y))\right)\to 0 $$

{\it c)}\, We then have that
$$\Pic^+(\bar X-Z, Y)^{\natural}\cong
[\Div^0_Z(\bar X,Y)\by{u^{\natural}_Z}\Pic^{\natural-\log}(\bar
X-Z,Y)^0].$$

{\it d)}\, We have an isomorphism
$$\Lie \Pic^{\natural-\log}(\bar X-Z, Y)^0\cong \HH^1(\bar X,\cO_{\bar
X}(-Y)\to \Omega^1_{\bar X}(\log (Y+ Z))(-Y)).$$
\end{lemma}
\begin{proof} The proofs of parts {\it a)}\, and {\it b)}\, are very similar to
those in Proposition~\ref{rignabla}. In fact we have relative residue sequences
given by the first row of the following commutative diagram (\cf
\cite[2.3]{EV})
\begin{equation}\label{relres}
\begin{array}{ccccccc}
&0&&0&&&\\
&\downarrow&&\downarrow&&&\\
0\to & \Omega^1_{\bar X}(\log Y)(-Y)
&\to &\Omega^1_{\bar X}(\log (Y+Z))(-Y)&\to &
\oplus_j \cO_{Z_j}(-Z_j\cap Y)&\to 0\\
&\downarrow&&\downarrow&&\quad\downarrow\cong&\\
0\to& \Omega^1_{\bar X}&\to &
\Omega^1_{\bar X}(\log Z)&\to &
\oplus_j \cO_{Z_j}&\to 0\\
&\downarrow&&\downarrow&&&\\
&\oplus_i\Omega^1_{Y_i}&\longby{\cong}&\oplus_i\Omega^1_{Y_i}(\log
Z\cap Y_i)&&&\\
&\downarrow&&\downarrow&&&\\
&0&&0&&&
\end{array}
\end{equation}
where the isomorphisms are because $Y\cap Z=\emptyset$. Here recall that
$Y_i$ are the irreducible components of $Y$.

For the latter claims {\it c)}\, and {\it d)}\, we proceed
as follows. Let
$$K \df \ker (\oplus_jH^0(Z_j,\cO_{Z_j})\to
H^1(\bar X,\Omega^1_{\bar X}(\log Y)(-Y)))$$ From (\ref{relres}) above
we get the following push-out diagram
$$\begin{array}{ccccccc}
&0&&0&&&\\
&\downarrow&&\downarrow&&&\\
0\to &H^0(\bar X,\Omega^1_{\bar X}(\log Y)(-Y))& \to &
\Pic^{\natural}(\bar X, Y)^0&\to& \Pic^0(\bar X, Y)&\to 0\\
&\downarrow&&\downarrow&&\veq&\\
0\to &H^0(\bar X,\Omega^1_{\bar X}(\log (Y+Z))(-Y))& \to
& \Pic^{\natural-log}(\bar X -Z, Y)^0&\to& \Pic^0(\bar X,Y)&\to 0\\
&\mbox{\ }\downarrow\s{\rm res}&&\downarrow\delta&&&\\
&K&\mbox{\large =}&K&&&\\
&\downarrow&&\downarrow&&&\\
&0&&0&&&
\end{array}$$
where ``res'' is the ordinary residue of forms and $\delta$ is the residue of
connections. Therefore we are left to show that the canonical induced map
\begin{equation}\label{vectmin}
\Ext (\Pic^+(\bar X-Z,Y),\G_a)^{\vee} \longby{\cong} H^0(\bar X,\Omega^1_{\bar
X}(\log (Y+Z))(-Y))
\end{equation}
is an isomorphism: in fact, granting (\ref{vectmin}), {\it c)}\, follows from
the above push-out diagram, Proposition~\ref{rignabla} and the construction of
the universal extension as being given in our Section~\ref{pre}.

In order to show the isomorphism in (\ref{vectmin}) we consider the following
commutative diagram
$$\begin{array}{ccc}
0&&0\\
\downarrow&&\downarrow\\
\Ext (\Pic^0(\bar X, Y),\G_a)^{\vee}& \by{\cong} & H^0(\bar X,
\Omega^1_{\bar X}(\log (Y))(-Y))\\
\downarrow&&\downarrow\\
\Ext (\Pic^+(\bar X-Z,Y),\G_a)^{\vee}&\to& H^0(\bar X,
\Omega^1_{\bar X}(\log (Y+Z))(-Y))\\
\downarrow&&\mbox{\ }\downarrow\s{\rm res}\\
\Hom (\Div^0_Z(\bar X,Y),\G_a)^{\vee}&\by{\cong}&K\\
\downarrow&&\downarrow\\
0&&0\\
\end{array}$$
where: by the Proposition~\ref{rignabla} we know
that  $\Ext (\Pic^0(\bar X,Y),\G_a)^{\vee}\cong H^0(\bar
X,\Omega^1_{\bar X}(\log Y)(-Y))$ and
$$\Hom (\Z[Z]^0,\G_a)^{\vee}\cong K\df \ker (\oplus_j
H^0(Z_j,\cO_{Z_j})\to
H^1(\bar X,\Omega^1_{\bar X}(\log Y)(-Y)))$$
is the restriction of the canonical isomorphism
$\Hom(\Z[Z],\G_a)^{\vee}\cong \oplus_jH^0(Z_j,\cO_{Z_j})$.

Thus the relative residue sequence (\ref{relres}) yields (\ref{vectmin}) as
well.
The Lie algebra computation yielding {\it d)}\, is then straightforward.
\end{proof}

\begin{rmk}\label{derham} (De Rham cohomology and homology)

{\em For a variety $X$ over a field $k$ of characteristic zero (not
necessarily algebraically closed), let $X_{\d}\to X$ be a smooth proper
hypercovering, and let $X_{\d}\into \bar X_{\d}$ be a smooth
compactification with normal crossing boundary $Y_{\d}$, as in \cite{D}.
We define the De Rham cohomology of $X$ as follows:
$$H_{DR}^*(X)\df \HH^*(\bar X_{\d},\Omega_{\bar
X_{\d}}^{\d}(\log Y_{\d})).$$
The Hodge-De Rham filtration is that induced by truncations, as usual.
This definition is in accordance with Deligne's definition in
\cite[10.3.15]{D}, determining De Rham cohomology as a filtered vector space.
Similarly, we define the Tate twist $H_{DR}^*(X)(m)$ to be the underlying
vector space of $H_{DR}^*(X)$, with the obvious shift in indexing of its
Hodge-De Rham filtration. Relative cohomology may also be defined in a
similar way, as in \cite{D} (see also \cite{ST}).

De Rham homology, denoted by $H_*^{DR}(X)\df H_{DR}^*(X)^{\vee}$, is defined to
be the dual (filtered) vector space: it differs, in general, from
Hartshorne (Borel-Moore) De Rham homology \cite{HAD}.

It can be shown, by comparison with the case $k =\C$ and cohomological descent,
that
\begin{enumerate}
\item[(i)] the underlying $k$-vector space of $H_{DR}^*(X)$ is naturally
identified with Hartshorne's algebraic De Rham cohomology \cite{HAD}
\item[(ii)]
$(H_{DR}^*(X),F^{\cdot})$ is independent of the choice of the
hypercovering and compactification $X_{\d}\into \bar
X_{\d}$
\item[(iii)] if $f:X\to Y$ is a morphism of $k$-varieties, the induced map
$f^*:H_{DR}^*(Y)\to H_{DR}^*(X)$ is strictly compatible with the
respective Hodge-De Rham filtrations; in particular, if the underlying
linear transformation is an isomorphism of vector spaces, then it is an
isomorphism of filtered vector spaces
\item[(iv)] if $\dim X=n$, then $H_{DR}^i(X)=0$ for $i>2n$, and for
irreducible $X$ over $k=\bar k$, $H_{DR}^{2n}(X)$ is either 0 (if $X$ is
not proper over
$k$) or 1-dimensional; if $X^{(n)}$ is the union of the $n$-dimensional
irreducible components of $X$, then $H_{DR}^i(X)\to H_{DR}^i(X^{(n)})$ is
an isomorphism for $i\geq 2n-1$
\item[(v)] $H_{DR}^*$ has other standard properties, like the excision
isomorphism, and the Mayer-Vietoris exact sequence; these are valid in the
category of filtered vector spaces and strictly compatible linear maps.
\end{enumerate}

If $X$ is smooth over $k$, and $\bar X$ is a smooth compactification with
normal crossing boundary $Y$, we may regard $\bar X$ and $Y$ as
``constant'' simplicial schemes, so that we obtain
$$H_{DR}^*(X)=\HH^*(\bar X,\Omega^{\d}_{\bar X}(\log Y)).$$

More generally, if $\bar X$ is a proper smooth $k$-variety and $Y$, $Z$
are disjoint normal crossing divisors, we get
$$H^*_{DR}(\bar X -Z,Y)\df \HH^*(\bar X, \Omega^{\cdot}_{\bar X}(\log
(Y+Z))(-Y))$$
with Hodge-De Rham filtration defined by truncation of the (twisted) log
De Rham complex.
}
\end{rmk}

Now Lemma~\ref{riglog} implies the following.
\begin{cor}\label{dercor}
 Let $\bar X$ be a non-singular proper $k$-variety, $Z$ and $Y$ disjoint
normal crossing divisors in $\bar X$.
Then there is a natural isomorphism of filtered $k$-vector spaces
$$T_{DR}(\Pic^+(\bar X -Z,Y))\cong H^1_{DR}(\bar X-Z,Y)(1).$$
\end{cor}

We now have  the following duality result.
\begin{lemma}\label{deduality} Let $V$ be an $n$-dimensional
proper smooth algebraic variety over a field of characteristic zero. Let
$A+B$ be a normal crossing divisor in $V$ such that $A\cap B
=\emptyset$. Then there is a functorial duality isomorphism
$$H^r_{DR}(V-A,B)(-n)\cong H_{2n-r}^{DR}(V-B,A)$$
which is compatible with the Hodge-De Rham filtration and
Poincar\'e-Lefschetz
duality.
\end{lemma}
\begin{proof} We can consider the following pairing
$$\begin{array}{c}
\Omega_V^{\d}(\log
(A+B))(-A)\otimes\Omega_V^{\d} (\log (A+B))(-B)\\
\downarrow \\
\Omega_V^{\d}(\log (A+B))(-A-B)\\
\downarrow\\
\Omega_V^{\d}
\end{array}$$
It will suffices to show that such a pairing yields non degenerate
pairings on hypercohomology
$$\HH^r(V,\Omega_V^{\d}(\log(A+B))(-A))\tensor_k
\HH^{2n-r}(V,\Omega_V^{\cdot}(\log(A+B))(-B))\to
\HH^{2n}(V,\Omega^{\d}_V)=H^n(V,\Omega^n_V)=k.\]
Since we are in characteristic zero we are left to show it for $k =\C$ for
which it is clear from the proof of Lemma~\ref{duality} and
Remark~\ref{remduality}. Alternately, one can deduce the duality
isomorphism, as in the
proof of Poincar\'e duality for algebraic De Rham cohomology, by reducing
to Serre duality (\cf \cite[III.8]{HAA}).
\end{proof}

\begin{thm}\label{deminus} Let $X$ be any $n$-dimensional $k$-variety,
where $k$ is algebraically closed of characteristic 0. Then
$$T_{DR}(\Pic^-(X))\cong H_{2n-1}^{DR}(X)(1-n).$$
\end{thm}\begin{proof}
As usual, we can reduce to the case when $X$ is equidimensional.
Fix a resolution $f:\tilde{X}\to X$ with good normal crossing
compactification $\bar X$ and boundary $Y$.

As above, let $Z$ be the union of all compact components of $\tilde S$.
We clearly have the following (see (\ref{relres})) relative residue
sequence
$$ 0\to  \Omega^1_{\bar X}(\log Y)(-Y) \to \Omega^1_{\bar X}(\log
(Y+Z))(-Y)\to \oplus_{j} \cO_{Z_j}\to 0$$
where $Z_j$ are the smooth irreducible components of $Z$. Moreover
$$\displaylines{
H_{2n-1}^{DR}(X)(1-n)\into H_{2n-1}^{DR}(X, S)(1-n)\cong
H_{2n-1}^{DR}(\tilde X,\tilde S)(1-n) \cr
\cong H_{2n-1}^{DR}(\tilde X, Z)(1-n)
\cong H^{1}_{DR}(\bar X - Z, Y)(1)\cr}$$
by excision and duality, \ie Lemma~\ref{deduality}, and we have the following
pull-back diagram (compare with (\ref{pbhr}))
\begin{equation}\label{pbdr}
\begin{array}{ccccccc}
0\to &H^{1}_{DR}(\bar X, Y)(1)&\to&H^{1}_{DR}(\bar X - Z, Y)(1)&\by{\rm
res}& \Div_{Z}^0(\bar X,Y)\otimes k&\to 0 \\
&\veq & &\uparrow & &\uparrow & \\
0\to&H^1_{DR}(\bar X, Y)(1)&\to & H_{2n-1}^{DR}(X)(1-n) &\to &\Div_{\bar
S/S}^0(\bar X,Y)\otimes k&\to 0  \end{array}
\end{equation}
by duality and the Mayer-Vietoris sequence for De Rham homology.

Consider the following pull-back diagram of 1-motives
$$\begin{array}{ccccccc}
0\to &\Pic^0(\bar  X, Y)&\to&\Pic^+(\bar X -Z,Y)&\to&
\Div_{Z}^0(\bar X,Y)[1]&\to 0 \\
&\veq & &\uparrow& &\uparrow &\\
0\to &\Pic^0(\bar  X, Y)&\to&\Pic^-(X)&\to&
\Div_{\bar S/S}^0(\bar X,Y)[1]&\to 0.
\end{array}$$
We then get the following commutative diagram (whose middle column implies
the theorem)
\begin{equation}\label{relder}
\begin{array}{ccccccc}
0\to &T_{DR}(\Pic^0(\bar  X, Y))&\to&T_{DR}(\Pic^+(\bar X -Z,Y))&\to&
T_{DR}(\Div_{Z}^0(\bar X,Y)[1])&\to 0 \\ &\cong\uparrow\quad & &\uparrow &
&\uparrow & \\
0\to&H^1_{DR}(\bar X, Y)&\to & H_{2n-1}^{DR}(X) &\to &\Div_{\bar
S/S}^0(\bar X,Y)\otimes k&\to 0 \\
&\cong\uparrow\quad & &\cong\uparrow\quad & &\cong\uparrow\quad &\\
0\to &T_{DR}(\Pic^0(\bar  X, Y))&\to&T_{DR}(\Pic^-(X))&\to&
T_{DR}(\Div_{\bar S/S}^0(\bar X,Y)[1])&\to 0
\end{array}
\end{equation}
where: {\it i)}\, top and bottom rows are obtained by applying $T_{DR}$ to
the earlier diagram of 1-motives, and are exact by construction (\cf
Section~\ref{pre}), {\it ii)}\, the second row is exact according to
(\ref{pbdr}) {\it iii)}\, the vertical isomorphisms are then obtained by
applying Proposition~\ref{rignabla}, Lemma~\ref{riglog} and
Corollary~\ref{dercor}, yielding the top row of (\ref{pbdr}) as the top
row of De Rham realizations of (\ref{relder}).
\end{proof}

\section{Cohomological Albanese 1-motive: $\Alb^+$}\label{cohomalb}
We keep the same notations and hypotheses of the previous section.

\subsection{Definition of $\Alb^+$}
Let $X$ be a variety over an algebraically closed field $k$ of
characteristic 0. To define our cohomological Albanese $\Alb^+(X)$, we
just take the Cartier dual of $\Pic^-(X)$. We are then left with finding a
``more explicit'' description of $\Alb^+$, if possible; this is given by
Proposition~\ref{picdual}, when $X$ is smooth, and by (\ref{formula}),
when $X$ is proper.

\begin{defn}{\rm For an algebraic variety $X$ over an algebraically closed
field $k$ of characteristic zero we define the following 1-motive
$$\Alb^+(X)\df\Pic^-(X)^{\vee}= [\Div^0_{\bar S/S}(\bar X, Y)\to
\Pic^0(\bar X, Y)]^{\vee}.$$

We call $\Alb^{+}(X)$ the {\it cohomological Albanese 1-motive\,} of
$X$. Since $\Pic^-$ is independent of the choices of resolutions and
compactifications so is $\Alb^+$.}\end{defn}

We recall that Deligne's definition of ``motivic cohomology'' of a curve
$C$ (see \cite{D}, \cf \cite{DM}) is the 1-motive
$$H^1_m(C)(1) \df [\Div^0_F(\bar C')\to \Pic^0(\bar C')]$$
where $\bar C'$ is a compactification of the semi-normalisation $C'$ of
the given curve $C$, such that $F= \bar C'-C'$ is a finite set of
non-singular points. We can relate Deligne's definition to ours.

\begin{propose}\label{delmot}
If $C$ is a curve (\ie a purely 1-dimensional variety) over an
algebraically closed field $k$ of characteristic 0, we have a canonical
identification
$$H^1_m(C)(1)\cong \Alb^+(C).$$
\end{propose}
\begin{proof} The normalisation $\tilde C$ of $C$ clearly factors through
the semi-normal curve $C'$ and the morphism $C'\to C$ is  bijective on
points, and so induces an isomorphism on the groups of Weil divisors.
We therefore have that $\Alb^+(C)=\Alb^+(C')$. On the other hand,
$H^1_m(C)(1)=H^1_m(C')(1)$ as well, from Deligne's definition.

We then can assume $C=C'$ itself to be semi-normal; let $\pi :
\tilde C\to C$ be the normalisation. First consider the compact case, \ie
$\bar C'= C'=C$. We then have a canonical quasi isomorphism
$$[\cO_{C}^*\to i_*\cO_S^*]\cong [R\pi_*\cO_C^*\to R\pi_*\tilde
i_*\cO_{\tilde S}^*]$$
where $i: S\into C$ is the imbedding of the finite set $S$ of singular
points and $\tilde i:\tilde S\into \tilde C$ the imbedding of the inverse
image of $S$: therefore, we get an isomorphism
$$\Pic (C, S) \cong \Pic (\tilde C, \tilde S).$$ From the exact sequences
(\ref{longpic}) we get the following diagram
 $$\begin{array}{ccccccc}
&0& & & & &\\
&\downarrow & & &\\
0\to & \frac{H^0(S,\cO_S^*)}{H^0(C,\cO_C^*)}&\to & \Pic^0(C, S)&\to
&\Pic^0
(C)&\to 0\\
&\downarrow & &\downarrow\veq& &\downarrow& \\
0\to & \frac{H^0(\tilde S,\cO_{\tilde S}^*)}{H^0(\tilde
C,\cO_{\tilde C}^*)}&\to & \Pic^0 (\tilde C,\tilde S)&\to &\Pic^0
(\tilde C)&\to 0\\
&\downarrow & &\downarrow & &\veq & \\
0\to & \Div_{\tilde S/S}^0(\tilde C)^{\vee}& \to &
\Pic^0(C)&\to & \Pic^0(\tilde C)& \to 0\\
&\downarrow & &\downarrow& && \\
&0 & & 0& &&
\end{array}$$
showing that $[\Div_{\tilde S/S}^0(\tilde C)\to \Pic^0(\tilde C)]$
is Cartier dual of $\Pic^0 (C)=H^1_m(C)(1)$ (\cf \cite{BZ}).

If $C$ is not compact, let $\bar C$ be a smooth compactification of the
normalization $\tilde{C}$, and set $F=\bar C-C$; then
$[\Div^0_F(\bar C)\to \Pic^0(\bar C)]$ dualizes to
$$0\to  \frac{H^0(F,\cO_F^*)}{H^0(\bar C,\cO_{\bar C}^*)} \to
\Pic^0(\bar C,F)\to \Pic^0(\bar C) \to 0$$
One can then see that the symmetric avatars of $\Alb^+(C)$ and
$H^1_m(C)(1)$ are the same, \eg by making use of the ``classical''
Lemma~\ref{dual}. \end{proof}

The proof of the following fact is left as an exercise for the reader.
\begin{lemma}\label{dual} Let $C$ be a non-singular projective curve. Let
$S$ and $T$
be disjoint finite sets of closed points. Then we have the following duality
isomorphism between 1-motives
$$[\Div_S^0(C)\to \Pic^0(C,T)]^{\vee} = [\Div_T^0(C)\to \Pic^0(C,S)].$$
\end{lemma}

We have that $\Alb^{+}(X)$ is a semi-abelian variety whenever $X$ is
proper over $k$; in fact, in this case $\bar X = \tilde X$, \ie
$Y=\emptyset$, and
$\Alb^{+}(X)$ is given by the following Cartier dual
$$\Alb^{+}(X) = [\Div^0_{\tilde S/S}(\tilde X)\to\Pic^0(\tilde X)]^{\vee}$$
Thus, if $X$ is a proper $k$-variety, $\Alb^{+}(X)$ can be
represented as an extension
\begin{equation}\label{formula}
0\to T(\tilde S/S) \to \Alb^{+}(X) \to  \Alb (\tilde X)\to 0
\end{equation}
where the torus $T(\tilde S/S)$ has character group
$\Div^0_{\tilde S/S}(\tilde X)$ (\cf Section~\ref{pre}).
Therefore, we can regard $\Alb^{+}(X)$ as a $\G_m$-bundle over the abelian
variety $\Alb (\tilde X)$.

If $X$ is a smooth variety over $k=\bar k$, we then have that $\tilde X =
X$, \ie $\tilde S = S =\emptyset$, whence $\Pic^{-}(X)$ is a semi-abelian
variety, and $\Alb^{+}(X)$ is given by a homomorphism from a lattice to an
abelian variety. It is natural to ask what these are, ``concretely''.

Let $\bar X$ be a non-singular proper variety over $k=\bar k$, and
$Y\subset \bar X$ a normal crossing divisor. Denote by $\Z^Y$ and
$\Z^{\bar X}$ the free abelian groups generated by the connected
components of $Y$ and $\bar X$ respectively. Then there is a canonical
homomorphism $\gamma : \Z^Y\to \Z^{\bar X}$ induced by the mapping that
takes a component of $Y$ to the component of $\bar X$ to which it belongs.
The kernel of $\gamma$ is generated by classes $[Y_I]-[Y_J]$ where $Y_I$
and $Y_J$ are distinct connected components of $Y$ contained in the same
component of $\bar X$.

Let $Y_I$ and $Y_J$ be distinct connected components of $Y$, contained
in the same component of $\bar X$, and choose (closed) points $y_I\in
Y_I$ and $y_J\in Y_J$.  Then we consider $a_{\bar X}(y_I-y_J)\in\Alb(\bar
X)$, where $a_{\bar X}:\cZ_0(\bar X)_0\to \Alb(\bar X)$ denotes the
Albanese mapping for zero-cycles of degree zero. If $\tilde{Y}\to Y$ is
the normalization, then $\tilde{Y}\to \bar{X}$ is a morphism
between smooth and proper varieties, and so yields a morphism
$\Alb(\tilde{Y})\to \Alb (\bar X)$ of abelian varieties.
Note that $\tilde{Y}=\coprod_iY_i$, and $\Alb(\tilde{Y})=\oplus_i
\Alb(Y_i)$, where $Y=\cup_iY_i$ is the decomposition into irreducible
components.
\begin{propose} \label{picdual}
Let $\bar X$ be a smooth proper $k$-variety, and $Y$ a normal crossing
divisor in $\bar X$. The Cartier dual of $\Pic^0(\bar X, Y)$ is the
1-motive given by the lattice
$$\Z^{(\bar X, Y)}\df \ker (\Z^Y\by{\gamma} \Z^{\bar X}) = T(\bar X,
Y)^{\vee},$$
the abelian variety
$$\frac{\Alb (\bar X)}{\im (\oplus_i \Alb (Y_i))}=
(\ker^0(\Pic^0(\bar X)\to \oplus_i \Pic^0(Y_i)))^{\vee}$$
and the homomorphism of group schemes
$$u_{X} : \Z^{(\bar X, Y)} \to \frac{\Alb (\bar X)}{\im (\oplus \Alb (Y_i))}$$
defined by $$u_{X}(Y_I-Y_J) = a_{\bar X}(y_I-y_J)
\pmod{\im\oplus_i \Alb (Y_i)}$$
where $Y_I$, $Y_J$ lie in the same component of $X$, and $y_I\in Y_I$,
$y_J\in Y_J$ are any closed points. Therefore,
$$\Alb^{+}(X) =
[\Z^{(\bar X, Y)} \by{u_X} \frac{\Alb(\bar X)}{\im(\oplus_i\Alb(Y_i))} ].$$
\end{propose}
\begin{proof} We first note that the homomorphism $u_{X}$ is well-defined;
in fact if $y'_I, y'_J$ is another such pair of points, then we easily
see that $a_{\bar X}(y_I-y_J)-a_{\bar X}(y'_I-y'_J)$ lies in the
image of $\oplus_i \Alb (Y_i)\to \Alb(\bar X)$ (first we consider the case
when the pair of points $y'_I, y_I$, as well as $y'_J,y_J$, each lie in
an irreducible component of $Y$; then we can deduce the general case).

By (\ref{longpic}) it is clear that the character group of the torus
$T(\bar X, Y)$ is given by the lattice $\Z^{(\bar X, Y)}$. The following
pull-back homomorphism between abelian varieties
$$\Pic^0(\bar X)\longby{\rho} \oplus_i \Pic^0(Y_i)$$
is dual to the following push-forward homomorphism
$$\oplus_i \Alb (Y_i)\to \Alb (\bar X)$$
Thus $$\coker\left( \oplus_i \Alb (Y_i)\to \Alb (\bar X)\right) =
(\ker^0\rho)^{\vee}$$
as claimed.

In order to check that the map $u_X$ is Cartier dual to $\Pic^0(\bar
X, Y)$, it suffices to show that $u_X$ coincides, on each generator
$[Y_I]-[Y_J]$ of $\Z^{(\bar X, Y)}$, with the analogous
homomorphism for the Cartier dual 1-motive. Choosing points $y_I\in Y_I$,
$y_j\in Y_J$ which are smooth on $Y$, one can reduce (by considering the
normalization of an irreducible curve passing through the pair of points,
and standard functoriality for Picard and Albanese varieties) to checking
the duality assertion when $\bar X$ is a smooth connected projective
curve, and $Y$ consists of 2 points, for which it is ``classical'' (see
\cite{BZ} for a more general statement; see also
\cite[Exemple,~pg.11-04]{SER}, and \cite{OJ}).
\end{proof}

We can now show that the Albanese 1-motive $\Alb^+$ is a birational
invariant of normal proper varieties, and that in fact it is given by the
Albanese variety of any resolution of singularities of $X$. More
generally, we have the following.
\begin{propose}\label{normplus} If $X$ is a {\rm normal} $k$-variety the
Albanese 1-motive $\Alb^+(X)$ is the Cartier dual of $\Pic^0(\bar X, Y)$.
In particular, if $X$ is also proper, then $\Alb^+(X)=\Alb(\bar X)$.
\end{propose}
\begin{proof} First consider the case when $X$ is a proper, normal
surface. The proposition is true in this case because  the intersection
matrix of the exceptional divisor of a desingularization of a normal
surface singularity is known \cite{Mum} to be negative definite:
the group $\Div^0_{\tilde S}(\tilde X)$ is zero since any non-zero linear
combination of compact components of $\tilde S$ cannot be numerically
equivalent to zero.

For higher dimensional proper $X$, we take $\bar X$ to be smooth and
projective; now by choosing successive hyperplane sections, we can find a
complete intersection smooth surface $T$ in $\tilde X$ and a
commutative square
$$\begin{array}{ccc} \Div_{\tilde S}(\tilde X)&\to & NS
(\tilde X)\\ \downarrow & &\downarrow \\
\Div_{\tilde S\cap T}(\tilde T)&\to & NS (T)
\end{array}$$
where $\tilde S\cap T\subset T$ is a reduced normal crossing divisor.
Since $T$ is general $\Div_{\tilde S}(\tilde X)$ injects into
$\Div_{\tilde S\cap T}(T)$. If $T_0$ is the normalization of the image
of $T$ in $X$, then $T\to T_0$ is a resolution of singularities of a
normal proper surface, with exceptional divisor $\tilde S\cap T$; hence
$\Div_{\tilde S\cap T}^0(T)=0$ by the case of surfaces considered above,
and so $\Div_{\tilde{S}}^0(\tilde{X})=0$ as well.

If $X$ is open we just notice that $\Div_{\bar S}^0(\bar X, Y)$ is
contained in $\Div_{\bar  S}^0(\bar X)$; however, the latter group
can be assumed to vanish, since $\bar X$ can be chosen to be
a projective resolution of a normal compactification of $X$.
\end{proof}

\begin{rmk}{\rm After Proposition~\ref{normplus}, we have the
following alternative  description of $\Pic^-(X)$, for a proper
$k$-variety $X$.

Let $X_n$ be the normalization of $X$, $S_n$ be the pullback of
the singular locus, and $\tilde X$ a resolution of the normalization. We
then have an exact sequence
$$0\to\Div_{\tilde S/S_n}(\tilde X)\to \Pic(\tilde X)\to \Cl (X_n)\to 0$$
where $\Cl$ denotes the divisor class group, and
$\Div_{\tilde{S}/S_n}(\tilde{X})$ is the group generated by exceptional
divisors for $\tilde{X}\to X_n$. Equivalently,
$\Div_{\tilde{S}/S_n}(\tilde{X})$ is the kernel of the push-forward map
$\Div(\tilde X)\to \Div(X_n)$; it is also the kernel of the pushforward
map $\Div_{\tilde S}(\tilde X)\to \Div_{S_n}(X_n)$.

We have $\Div_{\tilde S/S_n}^0(\tilde X)=0$, by
Proposition~\ref{normplus}. Hence $\Div_{\tilde{S}/S_n}(\tilde{X})$
has no intersection with $\Pic^0(\tilde X)$, and so
$\Pic^0(\tilde X)$ injects into $\Cl (X_n)$; denote its image by
$\Cl^0(X_n)$. Let $\Div_{S_n/S}(X_n)$ denote the group of Weil
divisors on $X_n$ have vanishing push-forward in $X$; these divisors
are necessarily supported on $S_n$. Let $\Div^0_{S_n/S}(X_n)$ be the
inverse image of $\Cl^0(X_n)$ under the obvious map
$\Div_{S_n/S}(X_n)\to \Cl(X_n)$ which send a Weil divisor to its divisor
class. We can now define a {\it class group}\, 1-motive of $X$ to be the
following 1-motive:
$$[\Div_{S_n/S}^0(X_n)\to \Cl^0(X_n)].$$
We then have that the homological Picard 1-motive $\Pic^-(X)$ is
canonically isomorphic to the class group 1-motive
$$\Pic^-(X)\cong [\Div_{S_n/S}^0(X_n)\to \Cl^0(X_n)].$$}
\end{rmk}

\subsection{Albanese mappings to $\Alb^+$}

Let $X$ be an equidimensional proper $k$-variety of dimension $n$, where
$k$ is algebraically closed of characteristic 0. Let $X_{\rm reg}$ denote
the set of smooth points of $X$.  We may also consider $X_{\rm reg}$ as an
open subscheme of any given resolution of singularities $\tilde{X}$.
Let $X_{\rm reg}=\coprod_jU_j$ be the decomposition into irreducible (or
equivalently connected) components. If $\tilde{X}\to X$ is a resolution,
then the Zariski closures $\bar{U_j}\subset \tilde{X}$ are the irreducible
(equivalently, connected) components of $\tilde{X}$.

Choose base points $x_j\in U_j$ for each $j$, and let ${\bf x}=\{x_j\}_j$.
Let $a_{\bf x}:\tilde{X}\to\Alb(\tilde{X})$ be the corresponding
Albanese mapping. Since $X$ is proper over $k$, $\Alb^+(X)$ is a torus
bundle over $\Alb(\tilde{X})$. Consider the following pull-back square
$$\begin{array}{ccc} \Alb^{+}(X) &\to &\Alb (\tilde X)\\
\tilde{a}_{\bf x}\ \uparrow\quad & &\quad\uparrow
a_{\bf x} \\
\Alb^{\dag}(X)&\to &\tilde X
\end{array}$$
Then $\Alb^{\dag}(X)$ is a torus bundle on $\tilde X$, with toric fiber
$$T(\tilde S/S)\df \Hom (\Div^0_{\tilde S/S}(\tilde X),\G_m)$$
We claim that the restriction of the torus bundle $\Alb^{\dag}(X)\to
\tilde{X}$ to the open subset $X_{\rm reg}\subset \tilde{X}$ has a natural
trivialization. In fact, dually, any divisor $D$ in $\Div^0_{\tilde
S/S}(\tilde X) = T(\tilde S/S)^{\vee}$ is mapped to the class in
$\Pic^0(\tilde{X})$ of the line bundle $\cO(D)$, which is canonically
trivialized on $X_{\rm reg}$, since $\supp(D)\cap X_{\rm reg}=\emptyset$.
Therefore, by a ``classical'' argument due to Severi (\cf \cite[\S1]{SER})
there is a section $\sigma: X_{\rm reg}\to \Alb^{\dag}(X)$. By composing
$\sigma$ with $\tilde{a}_{\bf x}$ we get the Albanese mapping
\begin{equation}
\label{plusalbmap} a_{\bf x}^+: X_{\rm reg}\to \Alb^+(X). \end{equation}
It is easy to see that $a_{\bf x}^+$ is independent of the choice of
the resolution of singularities $\tilde X$ of $X$.

If $X$ is not equidimensional, let $X^{(n)}$ denote the union of its
$n$-dimensional irreducible components. We define $X_{\rm reg}$ to be the
intersection of $X^{(n)}$ with the locus of smooth points of $X$. Since
$\Alb^+(X)=\Alb^+(X^{(n)})$, while $X_{\rm reg}\subset X^{(n)}_{\rm reg}$,
we
obtain an Albanese mapping $a_{\bf x}^+:X_{\rm reg}\to \Alb^+(X)$ by
restricting that of $X^{(n)}$, if the base points $x_j$ are chosen in
$X_{\rm reg}$.

\subsection{Hodge, \'etale and De Rham realizations of $\Alb^+$}

Let $X$ be an $n$-dimensional variety over $\C$. We recall that Cartier
duality for 1-motives is compatible, under the Hodge realization, with the
canonical involution $H\leadsto \Hom(H,\Z(1))$ on the category of mixed
Hodge structures. We thus have the following consequence of
Theorem~\ref{minus}.
\begin{cor} Let $X$ be as above. Then
$$T_{Hodge}(\Alb^+(X))\cong H^{2n-1}(X,\Z(n))/{\rm (torsion)}$$
\end{cor}
\begin{proof} We have the formula
$$\Hom(H_{2n-1}(X,\Z(1-n)),\Z(1)) = H^{2n-1}(X,\Z(n))/{\rm
(torsion)}$$
in the category of mixed Hodge structures. Cartier duality for
1-motives and Theorem~\ref{minus} then yield the result.
 \end{proof}

We let
$$J^n(X) \df \frac{H^{2n-1}(X,\C(n))}{F^0H^{2n-1}(X,\C(n)) + \im
H^{2n-1}(X,\Z(n))}.$$
We then have:
\begin{cor} Let $X$ be a proper variety over $\C$ and $n=\dim X$.
The Albanese 1-motive $\Alb^+ (X)$ is canonically isomorphic to the
semi-abelian variety $J^n(X)$, given as an algebraic extension
$$0\to T \to J^n(X)\by{f^*} J^n(\tilde X)\to 0$$
where $f:\tilde X\to X$ is any resolution of singularities, and the torus
$T$ is given by
$$ \frac{H^{2n-2}(\tilde S,\Z)}{\im (H^{2n-2}(S,\Z) \oplus H^{2n-2}(\tilde
X,\Z))} \otimes \C^*$$
\end{cor}
\begin{proof}
This follows from the Mayer--Vietoris sequence of mixed Hodge structures
$$H^{2n-2}(S,\Z(n)) \oplus H^{2n-2}(\tilde X,\Z(n))\to H^{2n-2}(\tilde
S,\Z(n))\to H^{2n-1}(X,\Z(n))\by{f^*} H^{2n-1}(\tilde X,\Z(n))\to 0$$
where $H^{2n-2}(\tilde S,\Z(n))$ is pure of weight $-2$ and
$H^{2n-1}(\tilde X,\Z(n))$ is pure of weight $-1$. In fact, the Deligne
1-motive canonically associated to $H^{2n-1}(X,\Z(n))$ is exactly the
claimed semi-abelian variety but, by the Theorem~\ref{minus}, the Hodge
realization of $\Alb^+ (X)$ is $H^{2n-1}(\tilde X,\Z(n))$ and the Hodge
realization functor is fully faithful. \end{proof}

We now let $X$ be a variety over an algebraically closed field $k$ of
characteristic zero.
\begin{cor} Let $X$ be as above and $n = \dim (X)$. Then
$$\hat{T}(\Alb^+(X))\cong H^{2n-1}_{\et}(X,\hat{\Z}(n))/{\rm (torsion)}$$
\end{cor}
\begin{proof} This follows from the formula
$$\Hom(H_{2n-1}^{\et}(X,\hat{\Z}(1-n)),\hat{\Z}(1)) =
H^{2n-1}_{\et}(X,\hat{\Z}(n))/{\rm (torsion)}$$
and Theorem~\ref{ellminus}.
\end{proof}

Let $X$ be any $n$-dimensional variety over an algebraically closed field of
characteristic zero as above. Recall (\ref{derham})
that $H^*_{DR}(X)\df \HH^*(\bar X_{\d}, \Omega^{\cdot}_{\bar X_{\d}}(\log
(Y_{\d})))$, the De Rham cohomology (filtered) $k$-vector spaces of $X$,
where $\bar X_{\d}$ is any smooth compactification of a proper smooth
hypercovering $X_{\d}$ with normal crossing boundary $Y_{\d}$.
\begin{cor} Let $X$ be as above. Then
$$T_{DR}(\Alb^+(X))\cong H^{2n-1}_{DR}(X)(n)$$
\end{cor}
\begin{proof} By Theorem~\ref{deminus}, as above, we get the result.
\end{proof}

\section{Cohomological Picard 1-motive: $\Pic^+$}\label{cohompic}
We first extend some results from the folklore on the Picard
functors into the language of simplicial schemes. Presumably, these are known
to experts, though we do not have any reference for these facts.

\subsection{Simplicial Picard functor}
Let $\pi: V_{\d}\to S$ be a simplicial scheme over a base scheme $S$. We
will denote by $\bPic (V_{\d})$ the group of isomorphism classes
of simplicial line bundles on $V_{\d}$ (\ie of invertible
$\cO_{V_{\d}}$-modules). We have the following description of $\bPic
(V_{\d})$. Denote by $d^{i}_k: V_i\to V_{i-1}$ the faces map of the given
simplicial scheme $V_{\d}$, and consider the following set of data and
conditions:
\begin{description}
\item[{\it - a line bundle\,}] $\ccL$ on $V_0$;
\item[{\it - an isomorphism\,}]
$\alpha :(d^{1}_0)^*(\ccL)\by{\cong}(d^{1}_1)^*(\ccL)$
on $V_1$;
\end{description}
satisfying the condition
\begin{description}
\item[{\it - cocycle condition\,}] the following composite
$$((d_1^{2})^*(\alpha))^{-1}\ssp{\circ}((d_2^{2})^*(\alpha))
\ssp{\circ}((d_0^{2})^*(\alpha))$$
yields $1\in \Gamma (V_2,\G_m)$, \ie if we let
$$f_0\df d^{1}_0d^{2}_0 = d^{1}_0d^{2}_1$$
$$f_1\df d^{1}_0d^{2}_2 = d^{1}_1d^{2}_0$$
$$f_2\df d^{1}_1d^{2}_2 = d^{1}_1d^{2}_1$$
then we want that the following diagram
\[\begin{TriCDV}
{f_2^*(\ccL)}{\> (d_0^{2})^*(\alpha)>>}{f_1^*(\ccL)}
{\ssp{ (d_1^{2})^*(\alpha)}\SE E E}{\SW W W \ssp{(d_2^{2})^*(\alpha)}}
{f_0^*(\ccL)}
\end{TriCDV}\]
commutes.
\end{description}

We clearly then have the following.
\begin{propose}\label{simpline} Let $V_{\d}$ be a simplicial scheme.
Elements of
$\bPic (V_{\d})$ isomorphically  corresponds to isomorphism classes of pairs
$(\ccL,\alpha)$ as above, satisfying the cocycle condition. Moreover,
there is a functorial isomorphism
$$\bPic (V_{\d})\cong \HH^1(V_{\d},\cO^*_{V_{\d}}).$$
\end{propose}
\begin{proof} The identification of $\bPic(V_{\d})$ with isomorphism
classes of pairs $(\ccL,\alpha)$ is easy, and left to the reader. For
a proof of the cohomological description, see Appendix~\ref{rep}.
\end{proof}

We now consider the {\it simplicial}\, Picard functor on
the category of schemes over $S$, which we denote as follows
$$T \leadsto \bPic_{V_{\d}/S}(T)$$
obtained by sheafifying the functor $$T\leadsto  \bPic
(V_{\d}\times_S T)$$ with respect to the $fpqc$-topology.
This means that if $\pi:V_{\d}\times_ST\to T$, then
$$\bPic_{V_{\d}/S}(T)\cong
H^0_{fpqc}(T,R^1\pi_*(\cO^*_{V_{\d}\times_ST})).$$

As usual, if $\pi_*(\cO_{V_{\d}}^*)=\cO_S^*$, the Leray spectral sequence
along $\pi$ and descent yields an exact sequence $$0\to
\Pic (S)\to \bPic (V_{\d})\to \bPic_{V_{\d}/S}(S)\to
H^2(S,\G_m)\to \HH^2(V_{\d},\cO^*_{V_{\d}}).$$
Furthermore, if there is a section of $\pi$, we have that $$
\bPic_{V_{\d}/S}(S)\cong \frac{\bPic (V_{\d})}{\Pic (S)}.$$

We are mainly interested in the case when $S$ is the
spectrum of a field $k$ and $X_{\d}$ is a proper (smooth) simplicial scheme
over $k$; the previous description for $k$-points of $\bPic_{X_{\d}/k}$ (\ie
the formula $\bPic_{V_{\d}/k}(k)\cong \bPic(V_{\d})$) applies in the
geometric case (\ie when $k$ is algebraically closed), since
$H^i(k,\G_m)=0$ for $i=1,2$ in that case; here, we do not need the
assumption that $\pi_*\cO_{X_{\d}}=k$.

In order to give another description of the simplicial Picard functor, which is
more suitable for our purposes, we consider the canonical spectral
sequence
\begin{equation}\label{simpcomp}
E^{p,q}_1 = H^q(X_p,\cO^*_{X_p})\implies
\HH^{p+q}(X_{\d},\cO^*_{X_{\d}})
\end{equation}
Let $\pi_i: X_i\to k$ denote the structure morphisms.
The spectral sequence yields the following exact sequence of
$fpqc$-sheaves:
\begin{equation}\label{semisimp} 0\to \frac{\ker ((\pi_1)_*\G_{m,X_1}\to
(\pi_2)_*\G_{m,X_2})}{\im ((\pi_0)_*\G_{m,X_0}\to (\pi_1)_*\G_{m,X_1})}\to
\bPic_{X_{\d}/k} \to \ker (\Pic_{X_0/k}\to\Pic_{X_1/k})
\end{equation}
We have the following facts.
\begin{lemma}\label{simpic} If $X_{\d}$ is smooth and proper over a field
$k$, then the simplicial Picard functor $\bPic_{X_{\d}/k}$ is
representable by a group scheme locally of finite type over $k$.
\end{lemma}
\begin{proof} See Appendix~\ref{rep}.
\end{proof}

For smooth proper simplicial schemes over $k=\bar k$ we have the
following description.
\begin{propose} Let $X_{\d}$ be smooth and proper over $k=\bar k$ of
characteristic zero. The sequence (\ref{semisimp}) yields a semi-abelian group
scheme over $k$, which can be represented as an extension
\begin{equation} 1\to T( X_{\d})\to \bPic^0( X_{\d})\to A(
X_{\d})\to 0 \end{equation}
where:
\begin{description}
\item[{\it (i)}] $\bPic^0( X_{\d})$ is the
connected component of the identity of $\bPic_{X_{\d}/k}$;
\item[{\it (ii)}] $T( X_{\d})$ is the $k$-torus defined by
$$ T( X_{\d}) \df \frac{\ker ((\pi_1)_*\G_{m,X_1}\to
(\pi_2)_*\G_{m,X_2})}{\im ((\pi_0)_*\G_{m,X_0}\to (\pi_1)_*\G_{m,X_1})}$$
where $\pi_i:X_i\to k$ are the structure morphisms;
\item[{\it (iii)}] $A( X_{\d})$ is the following abelian variety
$$A( X_{\d}) \df \ker^0(\Pic^0(X_0)\to\Pic^0(X_1))$$
obtained as the connected component of the identity of the kernel.
\end{description}
\end{propose}
\begin{proof} From Lemma~\ref{simpic}, by taking connected components of the
identity of the group schemes in (\ref{semisimp}), where $T( X_{\d})$
is connected, we claim that $\bPic^0(\bar
X_{\d})$ surjects onto the abelian variety
$\ker^0(\Pic_{X_0/k}^0\to\Pic_{X_1/k}^0)$: by the spectral sequence
(\ref{simpcomp}), the image of $\bPic(\bar X_{\d})$ is the kernel of the
following edge homomorphism $$\ker(\Pic_{X_0/k}\to\Pic_{X_1/k})\to\frac{\ker
((\pi_2)_*\G_{m,X_2}\to (\pi_3)_*\G_{m,X_3})}{\im  ((\pi_1)_*\G_{m,X_1}\to
(\pi_2)_*\G_{m,X_2})}$$ which vanishes on the connected component of the
identity of the domain. \end{proof}

\subsection{Definition of $\Pic^+$}
Now let $X_{\d}$ be a smooth simplicial $k$-variety.
This $X_{\d}$ can be regarded as obtained from a
simplicial pair $(\bar X_{\d},Y_{\d})$ such that  $X_{\d}
= \bar X_{\d}-Y_{\d}$, $\bar X_{\d}$ is a proper smooth
simplicial scheme and $Y_{\d}$ has components $Y_i$ which
are normal crossing divisors in $\bar X_i$. We then have a
spectral sequence
 $$E^{p,q}_1 = H^q_{Y_p}(\bar X_p,\cO^*_{\bar X_p})\implies
\HH^{p+q}_{Y_{\d}}(\bar X_{\d},\cO^*_{\bar X_{\d}}).$$
Since each component of $\bar X_{\d}$ is smooth we have that
$H^q_{Y_p}(\bar X_p,\cO^*_{\bar X_p})\neq 0$ if and only if
$q=1$ and  we clearly have that
$$H^1_{Y_p}(\bar X_p,\cO^*_{\bar X_p})\cong
\Div_{Y_p}(\bar X_p).$$ From the above spectral sequence we
then have
\begin{equation} \label{simpdiv}
\HH^{1}_{Y_{\d}}(\bar X_{\d},\cO^*_{\bar X_{\d}})\cong  \ker
(\Div_{Y_0}(\bar X_0)\by{d_0^*-d_1^*}\Div_{Y_1}(\bar X_1))
\end{equation}
We will denote by $\Div_{Y_{\d}}(\bar X_{\d})$ the subgroup of
divisors on $\bar X_0$ given by the right side of (\ref{simpdiv}).

The canonical mapping
\begin{equation}\label{map}
\Div_{Y_{\d}}(\bar X_{\d})=\HH^{1}_{Y_{\d}}(\bar X_{\d},\cO^*_{\bar
X_{\d}})\to \HH^{1}(\bar X_{\d},\cO^*_{\bar X_{\d}})\cong \bPic_{{\bar
X_{\d}}/k}(k)
\end{equation}
is compatible with the restriction of the map
taking a divisor on $\bar X_0$ to the associated line bundle.

In order to define our $\Pic^+$, we let
$$\bPic^0(\bar X_{\d})\df \bPic^0_{{\bar X_{\d}}/k}(k)\subset
\bPic(\bar X_{\d}),$$
and let $\Div^0_{Y_{\d}}(\bar X_{\d})$
denote the inverse image of $\bPic^0(\bar X_{\d})$ under the above
mapping (\ref{map}).

Now let $X$ be an algebraic variety over a field
$k=\bar{k}$ of characteristic zero. Let $\pi: X_{\d}\to X$ be a smooth
proper hypercovering of $X$, and choose a simplicial pair $(\bar X_{\d},
Y_{\d})$ as above (\ie $\bar X_{\d} - Y_{\d} = X_{\d}$ and $Y_{\d}$ has
normal crossings.)
\begin{defn} {\rm With the hypothesis and notation as above we
define the 1-motive
$$\Pic^+(X) \df [\Div_{Y_{\d}}^0(\bar X_{\d})\to \bPic^0(\bar
X_{\d})].$$
We call $\Pic^+(X)$ the {\it cohomological Picard 1-motive\,} of $X$.}
\end{defn}
\begin{rmk}{\rm If $X$ is smooth, let $\bar{X}$ be a smooth
compactification with normal crossing boundary $Y$. We may take
$\bar{X}_{\d}$ to be the constant simplicial scheme associated to $\bar
X$. Then we see easily that
$\Pic^+(X)\cong [\Div^0_Y(\bar{X})\to \Pic^0(\bar X)]$.

On the other hand, if $X$ is proper over $k$, then $\bar X_{\d} = X_{\d}$,
and $\Pic^{+}(X)=\bPic^0(X_{\d})$ is a semi-abelian variety.
}
\end{rmk}

\subsection{Hodge realization of $\Pic^+$}
Let $(\bar X_{\d}, Y_{\d})$ be a simplicial pair as above. For $k=\C$, by the
simplicial exponential sequence on $(\bar X_{\d})_{\rm an}$ and GAGA, we have
an isomorphism
\[\bPic(\bar X_{\d})=\HH^1(\bar X_{\d},\cO^*_{\bar X_{\d}})\cong
\HH^1((\bar X_{\d})_{\rm an},\cO^*_{(\bar X_{\d})_{\rm an}})\]
and a simplicial cycle map
$$c\ell_{\d}: \bPic (\bar X_{\d}) \to \HH^2(\bar X_{\d},\Z(1)).$$
\begin{lemma}\label{simpjac}
Let $\bar X_{\d}$ be as above and $k=\C$. Then
$${\bf J}^1(\bar X_{\d})\df \frac{\HH^1(\bar X_{\d},\C (1))}{F^0 +
\HH^1(\bar X_{\d},\Z (1))}\cong \bPic^0(\bar X_{\d})$$
and  $$ \Div_{Y_{\d}}^0(\bar X_{\d})\cong \ker
(\HH^2_{Y_{\d}}(\bar X_{\d},\Z (1))\to \HH^2(\bar X_{\d},\Z (1)))$$
Under these isomorphisms the canonical mapping
$\Div_{Y_{\d}}^0(\bar X_{\d})\to \bPic^0(\bar X_{\d})$
defined above is identified with an appropriate extension class map for
mixed Hodge structures on $\HH^1(X_{\d},\Z (1))$.
\end{lemma}
\begin{proof} From the simplicial exponential sequence,
since the complex $\Z(1)_{\d}\to\cO_{\bar X_{\d}}$ is
quasi-isomorphic to $\cO_{\bar X_{\d}}^*[-1]$ on $(\bar X_{\d})_{\rm an}$, we
have that $${\bf J}^1(\bar X_{\d}) \cong \ker c\ell_{\d}$$ because
$$\HH^1(\bar X_{\d},\cO_{\bar X_{\d}})\cong \frac{\HH^1(\bar X_{\d},
\C(1))}{F^0}.$$ Since we have a spectral sequence
$$E^{p,q}_1 = H^q_{Y_p}(\bar X_p,\Z(1))\implies
\HH^{p+q}_{Y_{\d}}(\bar X_{\d},\Z(1))$$
such that $E^{p,*}_1 =0$ for $q=0,1$, we obtain $\HH^1_{Y_{\d}}(\bar X_{\d},
\Z(1))=0$, and moreover  $$\HH^{2}_{Y_{\d}}(\bar
X_{\d},\Z(1))\cong \ker (H^{2}_{Y_0}(X_0,\Z(1))\by{d_0^*-d_1^*}
H^{2}_{Y_1}(X_1,\Z(1)))$$ whence $\HH^{2}_{Y_{\d}}(\bar X_{\d},\Z(1))\cong
\Div_{Y_{\d}}(\bar X_{\d})$. The following diagram
$$\begin{array}{c}
\Div_{Y_{\d}}^0(\bar X_{\d})\into \HH^{1}_{Y_{\d}}(\bar
X_{\d},\cO_{\bar X_{\d}}^*)\by{\cong}\HH^{2}_{Y_{\d}}(\bar
X_{\d},\Z(1))\\
\downarrow\hspace*{2cm}\downarrow\hspace*{2cm}\downarrow\hspace*{0.5cm}\\
\bPic^0(\bar X_{\d})\into \bPic(\bar X_{\d}) \by{c\ell_{\d}}
\HH^{2}(\bar X_{\d},\Z(1))
\end{array}$$
commutes, showing the claimed description of $\Div_{Y_{\d}}^0(\bar
X_{\d})$ (note that $\HH^{1}_{Y_{\d}}(\bar
X_{\d},\cO_{\bar X_{\d}}^*)$ is computed using the Zariski topology).

To show that the cycle class coincides with the extension class for the
mixed Hodge structure on $\HH^1(X_{\d},\Z(1))$, we consider the following
commutative diagram  of cohomology groups having exact rows and columns
$$\begin{array}{ccccccc} &&&&&&0\\
&&&&&&\downarrow\\
&&&&&&\bPic^0(\bar X_{\d})\\
&&&&&&\downarrow\\
&&&&\HH^1_{Y_{\d}}(\bar X_{\d},\cO_{\bar X_{\d}}^*)&\to&\HH^1(\bar X_{\d},
\cO_{\bar X_{\d}}^*)\\
&&&&\downarrow&&\downarrow\\
\HH^1(\bar X_{\d}\Z (1))&\to & \HH^1(X_{\d},\Z (1)) &\to &
\HH^2_{Y_{\d}}(\bar X_{\d},\Z (1))&\to & \HH^2(\bar X_{\d},\Z (1))\\
\downarrow&&\downarrow&&\downarrow&&\\
\HH^1(\bar X_{\d},\C (1))/F^0 &\to & \HH^1(X_{\d},\C (1))/F^0
&\to &\HH^2_{Y_{\d}}(\bar X_{\d},\C (1))/F^0& &\\
\downarrow&&&&&&\\
{\bf J}^1(\bar X_{\d})&&&&&&\\
\downarrow&&&&&&\\
0&&&&&&
\end{array}$$
The result then follows from a diagram chase (\cf the proof of
Lemma~\ref{piccus} and  \cite[Lemma~2.8]{BSL}).
\end{proof}
\begin{thm}\label{plus} Let $X$ be defined over $\C$. Then
$$T_{Hodge}(\Pic^+(X))\cong H^{1}(X,\Z(1)).$$
\end{thm}
\begin{proof} We have an exact sequence of mixed Hodge structures
$$0\to \HH^1(\bar X_{\d},\Z(1))\to H^1(X,\Z(1))\to \Div_{Y_{\d}}^0(\bar
X_{\d})\to 0$$
where $H^1(X,\Z(1))\cong \HH^1(X_{\d},\Z(1)_{\d})$ by universal
cohomological descent: the claim then follows from the
Lemma~\ref{simpjac}. \end{proof}

\subsection{\'Etale realization of $\Pic^+$}
Let $V_{\d}$ be any simplicial $k$-scheme. We first need to
recall the existence of the following long exact sequence
$$ \cdots\to\HH^0_{\et}(V_{\d},\G_{m})\to \HH^1_{\et}(V_{\d},\mu_{m})\to
\HH^1_{\et}(V_{\d},\G_{m})\by{m} \HH^1_{\et}(V_{\d},\G_{m})\to\cdots $$
and a ``simplicial Hilbert's Theorem 90''.

\begin{propose} \label{simp90} There is an isomorphism
$$\HH^1_{\et}(V_{\d},\G_{m})\cong \bPic (V_{\d})$$
\end{propose}
\begin{proof} Consider the Leray spectral sequence along
$\varepsilon_{\d} :(V_{\d})_{\et}\to (V_{\d})_{\rm Zar}$. Since
$(\varepsilon_{\d})_*(\G_m)= \G_m$ we then have a canonical functorial map
$$\varepsilon_{\d}^* :\HH^1(V_{\d},\cO^*_{V_{\d}})\to
\HH^1_{\et}(V_{\d},\G_{m}).$$ Consider the canonical spectral sequence
$$E^{p,q}_1 = H^q_{\et}( V_p,\G_m)\implies\HH^{p+q}_{\et}(V_{\d},\G_m)$$
A similar spectral sequence is clearly available for Zariski cohomology
groups, and $\varepsilon_{\d}^*$ is compatible with a morphism between the
respective spectral sequences.

Since we have that $H^q_{\et}(V_p,\G_m)=H^q_{\rm Zar}(V_p,\cO^*_{V_p})$
for all $p\geq 0$ and $q=0, 1$, via $\varepsilon_p:(V_p)_{\et}\to (V_p)_{\rm
Zar}$, we then get that $\varepsilon^*_{\d}$ is an isomorphism.
\end{proof}

We then can see that the cohomology group
$\HH^{1}_{\et}(V_{\d},\mu_m)$ is isomorphic to the group given by the
isomorphism classes of pairs $(\ccL_{\d}, \eta_{\d})$ where $\ccL_{\d}$
is a simplicial line bundle and $\eta_{\d}$ is an isomorphism $\eta_{\d}:
\ccL_{\d}^{\otimes m}\cong \cO_{V_{\d}}$. Moreover, we get the following
fact.
\begin{propose}\label{simpkum}
We have the ``simplicial Kummer sequence''
$$ 0\to \HH^0(V_{\d},\cO^*_{V_{\d}})/m\by{u}
\HH^1_{\et}(V_{\d},\mu_m)\by{p}
\bPic (V_{\d})_{m-{\rm tors}}\to 0$$
where: \begin{description}
\item[-] $\HH^1_{\et}(V_{\d},\mu_m)$ can be regarded as the group of
isomorphism classes of triples $(\ccL, \alpha, \eta)$ given by a line
bundle $\ccL$ on $V_0$, an isomorphism
$\alpha :(d_0)^*(\ccL)\by{\cong}(d_1)^*(\ccL)$
on $V_1$ satisfying the cocycle condition, an isomorphism
$\eta: \cO_{V_0}\by{\cong}\ccL^{\otimes m}$ on $V_0$ which is compatible
with $\alpha^{\otimes m}$ on $V_1$, \ie such that the composite of the
following  isomorphisms
$$\cO_{V_1}=(d_0)^*(\cO_{V_0})\by{(d_0)^*(\eta)}(d_0)^*(\ccL^{\otimes m})
\by{\alpha^{\otimes m}}(d_1)^*(\ccL^{\otimes m})
\by{(d_1)^*(\eta)^{-1}}(d_1)^*(\cO_{V_0})=\cO_{V_1}$$
is the identity on $\cO_{V_1}$ (here $d_0$ and $d_1$ from $V_1$ to $V_0$ are
the face maps of the simplicial scheme);
\item[-] $\HH^0(V_{\d},\cO^*_{V_{\d}})$
is given by those units  $u_0\in H^0(V_0,\cO^*_{V_{0}})$ such that
$d_0^*(u_0)=d_1^*(u_0)$ on $V_1$;
\item[-] the map $u$ is defined by taking a unit $u_0$ to the triple
$(\cO_{V_0}, 1, u_0^{-1})$;
\item[-] the map $p$ is defined by taking a triple $(\ccL, \alpha, \eta)$ to
a the torsion pair $(\ccL, \alpha)$ in the simplicial Picard group.
\end{description}
\end{propose}
\begin{proof} Taking into account
Propositions~\ref{simpline}~and~\ref{simp90}, the proof is an easy modification
of \cite[III.4]{MI}. \end{proof}

Now let $X$ be a $k$-variety, where $k$ is algebraically closed of
characteristic 0. Fix a smooth proper hypercovering $X_{\d}$ and a
normal crossing compactification $\bar X_{\d}$ with boundary $Y_{\d}$.
For $(D,\ccL_{\d})\in \Div_{Y_{\d}}^0(\bar
X_{\d})\times \bPic^0(\bar X_{\d})$ as above, by definition
$$T_{\sZ/m}(\Pic^+(X)) =\frac{\{(D,\ccL_{\d})\mid
\eta_{\d}:\ccL_{\d}^{\otimes m}
\cong\cO_{\bar X_{\d}}(-D)\}}{\{(mD,\cO_{\bar X_{\d}}(-D))\}}.$$
We have a canonical map
$$\rho_m: T_{\sZ/m}(\Pic^+(X))\to \HH^{1}_{\et}(X_{\d},\mu_m)$$
defined as follows
$$\rho_m(D,\ccL_{\d})\df (\ccL_{\d}, \eta_{\d})_{\mid X_{\d}}.$$
Note that $\eta_{\d}$ is well-defined up to multiplication by an element of
$k^*$, which is $m$-divisible, so that the isomorphism class of
$(\ccL_{\d}, \eta_{\d})_{\mid X_{\d}}$ is well-defined.

We let $$\HH^{*}_{\et}(X_{\d},\hat{\Z}(1))\df \liminv{m}
\HH^*_{\et}(X_{\d},\mu_{m})$$
We can show the following.
\begin{thm}\label{ellplus} Let $X$ be defined over $k$ which is algebraically
closed of characteristic zero. Then
$$\hat{T}(\Pic^+(X))\cong H^{1}_{\et}(X,\hat{\Z}(1)).$$
\end{thm}
\begin{proof} If we let $\Pic^+(X)$ be given by the 1-motive
$[\Div_{Y_{\d}}^0(\bar X_{\d})\to \bPic^0(\bar  X_{\d})]$ for a choosen
hypercovering and compactification with normal crossing boundary, we get the
following commutative diagram
$$\begin{array}{ccccccc}
0\to &\hat{T}(\bPic^0(\bar  X_{\d}))&\to&\hat{T}(\Pic^+(X))&\to&
\hat{T}(\Div_{Y_{\d}}^0(\bar X_{\d})[1])&\to 0 \\
&\bar\rho_{\et}\ \downarrow & &\downarrow\ \rho_{\et} & &\downarrow\
\rho_{\et}^0 & \\
0\to&\HH^1_{\et}(\bar X_{\d},\hat{\Z}(1))&\to &
H^1_{\et}(X,\hat{\Z}(1))&\to
&\ker (\HH^2_{\et,Y_{\d}}(\bar X_{\d},\hat{\Z}(1))\to\HH^2_{\et}(\bar
X_{\d},\hat{\Z}(1)))&\to 0
\end{array}$$
where (i) the bottom row is just the exact sequence of cohomology with
supports, (ii) we have $$\HH^1_{\et}(X_{\d},\hat{\Z}(1))\cong
H^1_{\et}(X,\hat{\Z}(1)),$$
since $X_{\d}\to X$ is a universal cohomological descent morphism, and
(iii) the top exact sequence is given by (\ref{etexseq}) in
Section~\ref{pre}.
We get the mapping $\rho_{\et}$ above by taking the inverse limit of
$\rho_m$; $\bar{\rho}_{\et}$ is the induced map, which can also be
regarded as the analogue of $\rho_{\et}$ for the case when $X$ is proper
(\ie $X_{\d}$ is proper and smooth over $k$). From the above description
the mapping $\rho_{\et}^0$ is an isomorphism: in fact, is easy to see
that we have an isomorphism
$\hat{T}(\Div_{Y_{\d}}(\bar X_{\d})) \cong
\HH^2_{\et,Y_{\d}}(\bar X_{\d},\hat{\Z}(1))$
such that the following diagram
$$\begin{array}{c}
\hat{T}(\Div_{Y_{\d}}(\bar X_{\d}))\longby{\cong}\HH^{2}_{\et, Y_{\d}}(\bar
X_{\d},\Z(1))\\ \downarrow\hspace*{2cm}\downarrow\hspace*{0.5cm}\\
\hat{\bPic(\bar X_{\d})} \longby{\hat{c\ell_{\d}}}
\HH^{2}_{\et}(\bar X_{\d},\hat{\Z}(1))
\end{array}$$
commutes (here $\hat{T}(\Div_{Y_{\d}}(\bar X_{\d}))$ and $\hat{\bPic(\bar
X_{\d})}$ are the profinite completions of $\Div_{Y_{\d}}(\bar X_{\d})$
and $\bPic(\bar X_{\d})$, respectively).

Granting this, we are left to show our claim is true for proper smooth
simplicial $k$-schemes, \ie that $\bar \rho_{\et}$ is an isomorphism. The
latter follows from the fact that the Neron-Severi group of such a
scheme (\ie the group of connected components of $\bPic(\bar X_{\d})$) is
finitely generated, whence $\hat{T}(\bPic^0(\bar X_{\d}))=
\liminv{m}\bPic(\bar X_{\d})_{m-{\rm tors}}$ and, by the simplicial
variants of Hilbert's Theorem 90 and Kummer theory (see (\ref{simp90})
and (\ref{simpkum})), we have
$\bPic(\bar X_{\d})_{m-{\rm tors}}=\HH^1_{\et}(\bar X_{\d},\mu_m).$
\end{proof}

\begin{rmk}\label{rempicplus}
{\rm From Theorem~\ref{ellplus} and
Proposition~\ref{fhat}, we
can see that the definition of $\Pic^+(X)$ is independent of choices of
hypercoverings and compactifications. In fact, given two smooth proper
hypercoverings we can always find a third one mapping to both, see
\cite[Expos\'e  V bis, 5.1.7 and 5.2.4]{SGA4}. Now let $X_{\d}'$ be a
smooth proper  hypercovering of $X$ with smooth compactification $\bar
X_{\d}'$, and let $X_{\d}''$ be another one, with smooth compactification
$\bar X_{\d}''$
mapping to $\bar X_{\d}'$ compatibly with the normal crossings boundaries.
Then, we get a map of 1-motives $$[\Div_{Y_{\d}'}^0(\bar X_{\d}')\to
\bPic^0(\bar  X_{\d}')]\to  [\Div_{Y_{\d}''}^0(\bar X_{\d}'')\to \bPic^0(\bar
X_{\d}'')]$$ by pulling-back cycles and simplicial line bundles. By the
Theorem~\ref{ellplus} $$\hat{T}([\Div_{Y_{\d}'}^0(\bar X_{\d}')\to \bPic^0(\bar
X_{\d}')])\cong \hat{T}([\Div_{Y_{\d}''}^0(\bar X_{\d}'')\to \bPic^0(\bar
X_{\d}'')]).$$
By Proposition~\ref{fhat} this isomorphism lifts to an isomorphism of
1-motives.

However, as for the case of $\Pic^-$, one would like to see directly, by a
geometric argument, that the above map of 1-motives is an isomorphism.}
\end{rmk}

\subsection{De Rham realization of $\Pic^+$}
Let $k$ be a field of characteristic 0. For any simplicial $k$-scheme
$X_{\d}$ we will denote by $\bPic^{\natural}(X_{\d})$ the group of
isomorphism classes of pairs $(\ccL_{\d},\nabla_{\d})$, where $\ccL_{\d}$
is a simplicial line bundle and $\nabla_{\d}$ is a simplicial
integrable connection
$$\nabla_{\d}:\ccL_{\d}\to \ccL_{\d}\otimes_{\cO_{X_{\d}}}
\Omega_{X_{\d}}^{1}.$$
We can consider the {\it simplicial}\, $\natural$-Picard functor on
the category of $k$-schemes, which we denote by
$\bPic_{X_{\d}/k}^{\natural}$, obtained by sheafifying the
functor
$$T\leadsto  \bPic^{\natural} (X_{\d}\times_k T)$$
with respect to the $fpqc$-topology.

For a given pair $(\ccL_{\d},\nabla_{\d})$  we clearly get a pair
$(\ccL,\nabla)$ on $X_0$ and an isomorphism
$\alpha:d_0^*(\ccL,\nabla)\by{\cong} d_1^*(\ccL,\nabla)$,
\ie $\alpha$ is an isomorphism $(d_0)^*(\ccL)\by{\cong}(d_1)^*(\ccL)$
which is compatible with the connections,
and, moreover, $\alpha$ satisfies the cocycle condition (\cf
Proposition~\ref{simpline}). In fact, we have the following description.

\begin{propose}\label{simpnabla}
Let $X_{\d}$ be any smooth proper simplicial $k$-scheme. Elements
of $\bPic^{\natural}(X_{\d})$ are in natural bijection with
isomorphism classes of triples $(\ccL,\nabla, \alpha)$ consisting of
an invertible sheaf $\ccL$ on $X_0$, with an integrable connection
$\nabla$, and an isomorphism $\alpha :d_0^*(\ccL,\nabla)\by{\cong}
d_1^*(\ccL,\nabla)$ satisfying the cocycle condition. There is a
functorial isomorphism
$$\bPic^{\natural}(X_{\d})\cong
\HH^1(X_{\d},\cO^*_{X_{\d}}\by{\rm dlog}\Omega^1_{X_{\d}})$$
\end{propose}
\begin{proof} It follows from the Proposition~\ref{simpline} and a
simplicial version of \cite[Sections~3-4]{MM} according to the
general hint given by Deligne in \cite[10.3.10]{D}.
\end{proof}

We have the following exact sequence
$$0\to\HH^0(X_{\d},\Omega^1_{X_{\d}}) \to  \bPic^{\natural}(X_{\d})\to
\bPic (X_{\d})\to \HH^1(X_{\d},\Omega^1_{X_{\d}})$$
which is obtained from the exact sequence of complexes of
simplicial sheaves
$$ 0 \to \Omega^1_{X_{\d}}[-1]\to [\cO^*_{X_{\d}} \by{\rm
dlog}\Omega^1_{X_{\d}}]\to \cO^*_{X_{\d}}\to 0$$
using the Proposition~\ref{simpnabla}.

Since $X_{\d}$ is smooth and proper over $k$, the
semi-abelian variety $\bPic^0(X_{\d})$ is mapped to zero in
$\HH^1(X_{\d},\Omega^1_{X_{\d}})$;
we thus get an extension
\begin{equation}\label{exnat} 0\to
\HH^0(X_{\d},\Omega^1_{X_{\d}}) \to
\bPic^{\natural}(X_{\d})^0\to \bPic^0(X_{\d})\to 0
\end{equation}
by pulling back along the inclusion $\bPic^0\into \bPic$. The
group extension (\ref{exnat}) is the group of $k$-points
of the universal $\G_a$-extension of the semi-abelian scheme
$\bPic^0_{X_{\d}/k}$; in fact we have the following.
\begin{lemma}\label{decomp}
Let $X_{\d}$ be a smooth proper simplicial $k$-scheme, where $k$ is
algebraically closed of characteristic 0. We have that
$$(\bPic^0_{X_{\d}/k})^{\natural}\cong(\bPic^{\natural}_{X_{\d}/k})^0,$$
and we have a canonical isomorphism
$$\Lie(\bPic^{\natural}_{X_{\d}/k})^0 \cong
\HH^1(X_{\d},\cO_{X_{\d}}\to\Omega^1_{X_{\d}}).$$
\end{lemma}
\begin{proof} The universal $\G_a$-extension of any semi-abelian
scheme is obtained as a pullback from the universal extension of
its abelian quotient. The abelian quotient of $\bPic^0_{X_{\d}/k}$
is $$A_{X_{\d}} \df \ker^0(\Pic^0_{X_0/k}\to\Pic^0_{X_1/k}).$$ By
\cite[Sections~3-4]{MM} it is easy to see that the universal
$\G_a$-extension of $A_{X_{\d}}$ is given by the group scheme (\cf
Section~\ref{pre})
$$\ker^0((\Pic^{\natural}_{X_0/k})^0\to(\Pic^{\natural}_{X_1/k})^0)$$
and we then have that
$$\Ext(A_{X_{\d}},\G_a)^{\vee} \cong
\ker(H^0(X_0,\Omega_{X_0}^1)\to H^0(X_1,\Omega_{X_1}^1)).$$

Everything then follows from the following diagram with exact rows
and columns,
$$\begin{array}{ccccccc}
&&&0&&0&\\
&&&\uparrow&&\uparrow&\\
0\to &T(X_{\d})&\to &\bPic^0(X_{\d})&\to &
\ker^0(\Pic^0_{X_0/k}\to\Pic^0_{X_1/k})&\to 0\\
&\veq&&\uparrow&&\uparrow&\\
0\to& T(X_{\d})&\to & \bPic^{\natural}(X_{\d})^0&\to &
\ker^0((\Pic^{\natural}_{X_0/k})^0\to(\Pic^{\natural}_{X_1/k})^0)
&\to 0\\ &&&\uparrow&&\uparrow&\\
&&&\HH^0(X_{\d},\Omega_{X_{\d}}^1)&\mbox{\large $=$}&\ker
(H^0(X_0,\Omega_{X_0}^1)\to H^0(X_1,\Omega_{X_1}^1))&\\
&&&\uparrow&&\uparrow&\\
&&&0&&0&
\end{array}$$
where $T(X_{\d})$ is the toric part of $\bPic^0_{X_{\d}/k}$ and
the middle row and column are exact by Proposition~\ref{simpnabla} and
(\ref{exnat}). Therefore, by taking associated sheaves, we see that
$(\bPic^{\natural}_{X_{\d}/k})^0$ is representable by the
pull-back of the universal extension of $A_{X_{\d}}$. Finally,
since the Lie algebra of  $\bPic^{0}_{X_{\d}/k}$ is
$H^1(X_{\d},\cO_{X_{\d}})$, from (\ref{exnat}), we get the last claim by
taking Lie algebras (\cf \cite{MM}). \end{proof}

More generally, let $X_{\d}$ be a smooth simplicial $k$-variety, where
$k$ is a field of characteristic 0. Let $\bar X_{\d}$ be a smooth
compactification with normal crossing boundary $Y_{\d}$. We then can
define $\bPic^{\natural-\log}(X_{\d})$ to be the group of isomorphism
classes of pairs $(\ccL_{\d},\nabla_{\d}^{\log})$, where $\ccL_{\d}$ is a
simplicial line bundle on $\bar X$ and $\nabla_{\d}^{\log}$ is a
simplicial integrable connection with log poles along $Y_{\d}$, \ie
$\nabla_{\d}^{\log}$ is a $k$-linear simplicial sheaf homomorphism
$$\nabla_{\d}^{\log}:\ccL_{\d}\to \ccL_{\d}\otimes_{\cO_{X_{\d}}}
\Omega_{\bar X_{\d}}^{1}(\log Y_{\d})$$
satisfying the Leibniz product rule (\cf \cite[Section~2]{EV},
\cite[Section~3]{MM}, \cite{DD}).
We clearly have a natural injective homomorphism
$$\bPic^{\natural}(\bar X_{\d})\longto
\bPic^{\natural-\log}(X_{\d})$$ and we have the following
cohomological description.
\begin{propose}\label{lognabla}
Let $X_{\d}$ be any smooth simplicial $k$-variety. Elements
of $\bPic^{\natural-\log}(X_{\d})$ are in bijection with
isomorphism classes of triples $(\ccL,\nabla^{\log}, \alpha)$
consisting of an invertible sheaf $\ccL$ on $\bar X_0$, with an
integrable connection with log poles
$$\nabla^{\log}: \ccL \to \ccL \otimes_{\cO_{X_0}}
\Omega_{\bar X_0}^{1}(\log Y_0)$$ and an
isomorphism $\alpha :d_0^*(\ccL,\nabla^{\log})\by{\cong}
d_1^*(\ccL,\nabla^{\log})$, satisfying the cocycle condition. There
is a functorial isomorphism
$$\bPic^{\natural-\log}(X_{\d})\cong
\HH^1(\bar X_{\d},\cO^*_{\bar X_{\d}}\by{\rm dlog}
\Omega^1_{\bar X_{\d}}(\log Y_{\d}))$$
\end{propose}
\begin{proof} A variant of the proof
of Proposition~\ref{simpnabla} (\cf Lemma~\ref{riglog}).
\end{proof}

Now let $X$ be a $k$-variety, where $k$ is a field of characteristic 0.
Let $X_{\d}\to X$ be a smooth proper hypercovering, and $\bar X_{\d}$ a
smooth proper compactification with normal crossing boundary $Y_{\d}$. We
recall (\ref{derham}) that by the De Rham cohomology of $X$ we mean the
graded, filtered $k$-vector space
$$H_{DR}^*(X)\df \HH^*(\bar X_{\d},\Omega_{\bar X_{\d}}^{\cdot}(\log
Y_{\d})).$$

\begin{thm}\label{deplus} Let $X$ be a $k$-variety, where $k$ is
algebraically closed of characteristic 0. Then
$$T_{DR}(\Pic^+(X))\cong H^{1}_{DR}(X)(1).$$
\end{thm}
\begin{proof}
Let $\Pic^+(X)$ be given by the 1-motive
$[\Div_{Y_{\d}}^0(\bar X_{\d})\to \bPic^0(\bar  X_{\d})]$
for a choosen hypercovering and compactification with normal crossing
boundary $Y_{\d}$. We have the following exact sequence of complexes
$$\begin{array}{ccccccc} 0\to &
\Omega^1_{\bar X_{\d}}&\to &\Omega^1_{\bar X_{\d}}(\log Y_{\d})&
\to &\cQ_{\d} &\to 0\\
&\  \ \uparrow \s{\dlog}& &\  \ \uparrow
\s{\dlog}& & \uparrow &\\
0\to &\cO^*_{\bar X_{\d}}&\mbox{\large $=$}&
\cO^*_{\bar X_{\d}} &\to & 0&
\end{array}$$
where $\cQ_{\d}$ is just the quotient $\Omega^1_{\bar X_{\d}/k}(\log
Y_{\d})/\Omega^1_{\bar X_{\d}/k}$.
We therefore have the following push-out diagram
 $$\begin{array}{ccccccc}
0\to &\HH^0(\bar X_{\d},\Omega^1_{\bar X_{\d}})& \to &
\bPic^{\natural}(\bar X_{\d})^0&\to& \bPic^0(\bar X_{\d})&\to 0\\
&\downarrow&&\downarrow&&\veq&\\
0\to &\HH^0(\bar X_{\d},\Omega^1_{\bar X_{\d}}(\log Y_{\d}))& \to
& \bPic^{\natural-log}(X_{\d})^0&\to& \bPic^0(\bar X_{\d})&\to 0
\end{array}$$
where the top row is (\ref{exnat}) and the bottom row is
obtained from the dlog map as well. From the latter we are then
left to show that
\begin{equation}\label{vectlog}
\Ext (\Pic^+,\G_a)^{\vee} =
\HH^0(\bar X_{\d},\Omega^1_{\bar X_{\d}}(\log Y_{\d})).
\end{equation}
In fact, granting (\ref{vectlog}), we have that, by the push-out diagram and
the Lemma~\ref{decomp}, the universal $\G_a$-extension of
$\Pic^+(X)$ is given by
$$\Pic^+(X)^{\natural}\cong [\Div_{Y_{\d}}^0(\bar X_{\d})\by{u^{\natural}}
\bPic^{\natural-log}(X_{\d})^0]$$
where the lifting $u^{\natural}$ above of $u: \Div_{Y_{\d}}^0(\bar X_{\d})\to
\bPic^{0}(\bar X_{\d})$ can be described as in Lemma~\ref{riglog} {\it
via}\, the Proposition~\ref{lognabla}. Therefore
$$\begin{array}{rl}
T_{DR}(\Pic^+(X)) &\df \Lie\bPic^{\natural-log}(X_{\d})^0\\
& \cong \Lie(\HH^1(\bar X_{\d},\cO^*_{\bar X_{\d}}\by{\rm
dlog}\Omega^1_{\bar X_{\d}}(\log Y_{\d}))\\
& \cong  \HH^1(\bar X_{\d},\cO_{\bar X_{\d}}\to
\Omega^1_{\bar X_{\d}}(\log Y_{\d}))\\
&\cong \HH^1(\bar X_{\d},\Omega^{\cdot}_{\bar X_{\d}}(\log Y_{\d}))\\
&\df H^1_{DR}(X)
\end{array}$$
Moreover, this isomorphism is clearly compatible with the Hodge
filtrations, provided we shift the index of the filtration on the right
by 1.

In order to show (\ref{vectlog}) we consider the
following commutative diagram with exact columns
$$\begin{array}{ccc}
0&&0\\
\downarrow&&\downarrow\\
\Ext (\bPic^0(\bar X_{\d}),\G_a)^{\vee}& \by{\cong} & \HH^0(\bar
X_{\d},\Omega^1_{\bar X_{\d}/k})\\
\downarrow&&\downarrow\\
\Ext (\Pic^+(X),\G_a)^{\vee}&\to& \HH^0(\bar X_{\d},\Omega^1_{\bar
X_{\d}/k}(\log Y_{\d}))\\
\downarrow&&\mbox{\ }\downarrow\s{\rm res}\\
\Hom (\Div_{Y_{\d}}^0(\bar X_{\d}),\G_a)^{\vee}&\by{\cong}&
\ker (\HH^0(\bar X_{\d},\cQ_{\d})\to
\HH^1(\bar X_{\d},\Omega^1_{\bar X_{\d}/k}))\\
\downarrow&&\downarrow\\
0&&0\\
\end{array}$$
where the horizontal maps are the canonical maps induced by
universality; from the previous Lemma~\ref{decomp} we know that  the
horizontal map on top is an isomorphism, so that we are left to show
that the horizontal map at the bottom is an isomorphism.

If $Y_{\d}$ is smooth this last claim is clear since we have
a simplicial surjective Poincar\'e residue map
$${\rm res}_{\d}:\Omega^1_{\bar X_{\d}}(\log Y_{\d})\to \cO_{Y_{\d}},$$
and therefore
$\HH^0(\bar X_{\d},\cQ_{\d})\cong\HH^0(Y_{\d},\cO_{Y_{\d}})$.
In general, since the subschemes $Y_i\subset X_i$ are
normal crossing divisors, for each $i\geq 0$,
we have exact sequences (\cf \cite[2.3]{EV})
$$0\to  \Omega^1_{\bar X_{i}}\to \Omega^1_{\bar X_{i}}(\log Y_{i})
\to\oplus_{ji}\cO_{Y_{ji}} \to 0$$
where the index $ji$ ($i$ fixed) runs over the smooth components
of $Y_i$, \ie $Y_i=\cup_{ji} Y_{ji}$. These sequences are compatible {\it
via}\, the face and degeneracy maps of the simplicial scheme $\bar
X_{\d}$. Because of this construction, and the definition of global
sections of a simplicial sheaf, we clearly get a canonical identification
$$\HH^0(\bar X_{\d},\cQ_{\d})=\ker (\oplus_{j0}H^0(Y_{j0},\cO_{Y_{j0}})\to
\oplus_{j1}H^0(Y_{j1},\cO_{Y_{j1}}))\cong \Div_{Y_{\d}}(\bar
X_{\d})\otimes k.$$
We finally then get
$\Div_{Y_{\d}}^0(\bar X_{\d})\otimes k\cong
\ker (\Div_{Y_{\d}}(\bar X_{\d})\otimes k\to
\HH^1(\bar X_{\d},\Omega^1_{\bar X_{\d}}))$
as claimed.
\end{proof}

\begin{rmk}{\rm From Theorem~\ref{deplus}, we obtain an ``algebraic
proof'' (\ie independent of base change to $\C$ and comparison with the
analytic topology) that $(H_{DR}^1(X),F^{\cdot})$ is independent of the
choices of $X_{\d}$ and its compactification $\bar X_{\d}$, since the
1-motive $\Pic^+$ is independent of these choices, as we saw earlier using
\'etale realizations, as a consequence of Theorem~\ref{ellplus}.}
\end{rmk}

\begin{rmk}{\rm For a given singular variety $X$ we can consider a singular
compactification $\bar X$ in such a way that $\bar X_{\d}$, in our above
considerations, is a hypercovering of $\bar X$. By the previous argument,
in the proof of the Theorem~\ref{deplus}, we constructed the following
extension
$$\begin{array}{ccccccc} 0\to &T_{DR}(\bPic^0(\bar
X_{\d}))&\to&T_{DR}(\Pic^+(X))&\to& T_{DR}(\Div_{Y_{\d}}^0(\bar
X_{\d})[1])&\to 0
\\ &\cong\uparrow\quad & &\cong\uparrow\quad & &\cong\uparrow\quad & \\
0\to&H^1_{DR}(\bar X)(1)&\to &
H^1_{DR}(X)(1)&\by{\rm res} &\Div_{Y_{\d}}^0(\bar X_{\d})\otimes k&\to 0
\end{array}$$
The resulting bottom row can be regarded as obtained from an exact sequence of
cohomology with supports as well as a Poincar\'e ``residue'' map
compatible
with the Hodge filtration.}
\end{rmk}

\section{Homological Albanese 1-motive: $\Alb^-$}\label{homalb}
We keep the notations and hypotheses from the previous section.
\subsection{Definition of $\Alb^-$}
In order to define our homological Albanese $\Alb^-(X)$ we just take
the Cartier dual of $\Pic^+(X)$.

\begin{defn}{\rm If $X$ is a variety over an algebraically
closed field $k$ of characteristic zero, define the following
1-motive
$$\Alb^-(X) \df \Pic^{+}(X)^{\vee} = [\Div_{Y_{\d}}^0(\bar X_{\d})\to
\bPic^0(\bar  X_{\d})]^{\vee},$$
where $X_{\d}\to X$ is a smooth proper hypercovering, and $\bar X_{\d}$
a smooth compactification of $X_{\d}$ with normal crossing boundary
$Y_{\d}$. We call $\Alb^-(X)$ the {\it homological Albanese 1-motive\,} of
$X$.} \end{defn}

\begin{lemma} \label{surj}
If $X$ is proper over $k=\bar k$ of characteristic 0, and $\pi:X_{\d}\to
X$ is any proper hypercovering, then the natural homomorphism between
connected algebraic $k$-groups
$$\pi^*:\Pic^0(X)\onto\bPic^0(X_{\d})$$
is a surjection with torsion free kernel.
\end{lemma}
\begin{proof} In fact, by Kummer theory and cohomological descent we get the
following commutative square of isomorphisms
$$\begin{array}{ccc}
\HH^1_{\et}(X_{\d},\mu_m)&\by{\cong}&\bPic (X_{\d})_{m-{\rm tors}}\\
\cong\uparrow\quad & &\cong\uparrow\quad \\
H^1_{\et}(X,\mu_m)&\by{\cong}&\Pic (X)_{m-{\rm tors}}
 \end{array}$$
Therefore, since the Neron-Severi groups are finitely generated, the Tate
module of $\Pic^0(X)$ is isomorphic to $\hat{T}(\bPic^0(X_{\d}))$. To
conclude we remark that $\bPic^0(X_{\d})$ is the group of $k$-points of a
semi-abelian variety, in which torsion points are Zariski dense (\cf
Proposition~\ref{fhat}).
\end{proof}
\begin{rmk}{\rm
As a consequence, we see that for any smooth proper hypercovering
$X_{\d}\to X$ of a proper $k$-variety $X$, the simplicial Picard variety
$\bPic^0(X_{\d})$ is the semi-abelian quotient of the
connected commutative algebraic group $\Pic^0(X)$.
}
\end{rmk}

For a proper smooth hypercovering $\pi:X_{\d}\to X$ of a proper $k$-scheme $X$,
where $k=\bar k$, we let $\Z^{X_i}$ denote the free abelian group on
connected components of $X_i$. Let
$$L_{X_{\d}}\df \frac{\ker(\Z^{X_1}\to\Z^{X_0})}
{\im (\Z^{X_2}\to\Z^{X_1})} \pmod{\rm torsion}$$
and consider the following abelian variety
$$(\ker^0(\Pic^0(X_0)\longby{d_0^*-d_1^*} \Pic^0(X_1)))^{\vee} =
\frac{\Alb (X_0)}{(d_0-d_1)_*\Alb (X_1)}.$$
Here $(d_0-d_1)_*\df (d_0)_*-(d_1)_*: \Alb (X_1)\to \Alb (X_0)$. Let
$(d_0)_*-(d_1)_*:\cZ_0(X_1)\to \cZ_0(X_0)$ denote as well the induced map
between  zero-cycles, and let $a_{X_0}: \cZ_0(X_0)_0\to \Alb (X_0)$ be the
Albanese map, defined on the subgroup $\cZ_0(X_0)_0\subset \cZ_0(X_0)$ of
zero-cycles on $X_0$ which have degree zero on each component of $X_0$.
The Albanese variety of the smooth, proper $k$-variety $X_0$ is here
defined to be the product of the Albanese varieties of its connected
components (note that these connected components are irreducible, smooth,
proper $k$-varieties, which need not have the same dimension).

\begin{propose} Let $X$ be proper over
$k =\bar k$. Then $\Alb^-(X)$ coincides with the 1-motive
$$u_{X_{\d}}^{\;}:L_{X_{\d}}\to \frac{\Alb (X_0)}{(d_0-d_1)_*\Alb (X_1)}$$
where the map $u_{X_{\d}}^{\;}$ is defined as follows. For each connected
component $X^c$ of $X_1$ choose a closed point $x_c\in X^c$. Then, for
$\sum n_c X_c\in \ker(\Z^{X_1}\to\Z^{X_0})$, we have
\[\sum n_c ((d_0)_*(x_c)-(d_1)_*(x_c))\in \cZ_0(X_0)_0,\]
 and we define
$$u_{X_{\d}}^{\;}(\sum n_c X_c \pmod{ \im \Z^{X_2}})\df a_{X_0}(\sum n_c
(d_0)_*(x_c)-(d_1)_*(x_c))
\pmod{(d_0-d_1)_*\Alb (X_1)}.$$
If $X$ is also {\em normal} then $L_{X_{\d}}=0$ and
$$\Alb^-(X)\cong \Pic^0(X)^{\vee}$$
is an abelian variety.
\end{propose}
\begin{proof} To check that $u_{X_{\d}}^{\;}$ is well defined is left as
an exercise. We recall that $\bPic^0(X_{\d})$ is an extension of the
abelian variety
$\ker^0(\Pic^0(X_0)\to \Pic^0(X_1))$ by the torus $T(X_{\d})$ (see
(\ref{semisimp}) in Section~\ref{cohompic} for a description of the
torus). Now, $L_{X_{\d}}$ is the character group of the torus $T(X_{\d})$,
$\Alb(X_0)/\im \Alb (X_1)$ is the dual abelian variety; the claimed
map between them is obtained from Cartier duality -- as in the proof of
the corresponding assertion of Proposition~\ref{picdual}, using standard
functoriality properties of Albanese and Picard varieties, one can reduce
to the case of the standard smooth proper hypercovering of an irreducible
projective curve with 1 node; now we further reduce to determining the
Cartier dual of $\Pic^0$ of this singular curve, which is treated in
\cite{BZ}.

If $X$ is normal, then $\pi_*(\cO_{X_{\d}}^*)=\cO_{X}^*$, and so
$$\pi^*:\Pic(X)\into\bPic(X_{\d})$$
is injective, by the Leray spectral sequence for the sheaf
$\cO_{X_{\d}}^*$ along $\pi$; therefore, from Lemma~\ref{surj} we get
$$\pi^*:\Pic^0(X)\cong\bPic^0(X_{\d}).$$
Since $X$ is normal, $\Pic^0(X)$ is an abelian variety \cite{Ch},
therefore $T(X_{\d})=0$.
\end{proof}

\begin{rmk}{\rm
If $X$ is smooth, possibly open, then $X_{\d}$, $\bar X_{\d}$, $Y_{\d}$
can be taken to be the constant simplicial schemes associated to $X$,
$\bar X$ and $Y$, respectively, where $\bar X$ is a smooth
compactification of $X$ with normal crossing boundary $Y$. In this case
$\Alb^-(X)$ is a semi-abelian variety, which can be represented by an
extension
$$1\to T(Y)\to \Alb^-(X)\to \Alb (\bar X)\to 0$$
where, by definition, $T(Y)$ is the $k$-torus with character group
$\Div_{Y}^0(\bar X)$; see Proposition~\ref{serre} below.
}\end{rmk}

\begin{rmk}
{\rm As a consequence of Lemma~\ref{surj} we have that
$$H^1(X,\cO_X)\onto H^1(X,\pi_*(\cO_{X_{\d}}))$$
is always a surjection, and the following edge homomorphism
$$\HH^1(X_{\d},\cO_{X_{\d}})\longby{\rm zero} H^0(X,R^1\pi_*(\cO_{X_{\d}}))$$
is the zero map. In fact, since $H^1(X,\cO_X) =\Lie \Pic^0(X)$ and
$\HH^1(X_{\d},\cO_{X_{\d}}) =\Lie \bPic^0(X_{\d})$, we see that $H^1(X,\cO_X)$
always surjects onto $\HH^1(X_{\d},\cO_{X_{\d}})$; moreover,
$$H^1(X,\cO_X)\cong \HH^1(X_{\d},\cO_{X_{\d}})$$
if $X$ is normal, and we then have that
$$\Lie \Alb^-(X) = H^1(X,\cO_X)^{\vee}.$$}
\end{rmk}

\subsection{Albanese mappings to $\Alb^-$}
Let $X_{\rm reg}$ be the smooth locus of an equidimensional $k$-variety
$X$, where $k$ is algebraically closed of chaacteristic 0. We then have
that $X_{\rm reg}=\tilde X - \tilde S =\bar X - (Y\cup \bar S)$ for a
resolution of singularities $\tilde X$ and a good normal crossing
compactification $\bar X$, with boundary divisor $Y$; also $Y\cup
\bar S$ is a normal crossing divisor in $\bar X$. We then  have a
commutative square of 1-motives
$$\begin{array}{ccc}
[0\to\Pic^0(\bar X, Y)] &\to &[\Div_Y^0(\bar X)\to\Pic^0(\bar X)]\\
\downarrow & &\downarrow \\{}
[\Div_{\bar S/S}^0(\bar X, Y)\to\Pic^0(\bar X, Y)]&\to& [\Div_{\bar S\cup
Y}^0(\bar X)\to\Pic^0(\bar X)]
 \end{array}$$
which we may rewrite as
$$\begin{array}{ccc} \Pic^-(\tilde
X)&\longto &\Pic^+(\tilde X) \\ \downarrow & &\downarrow \\
\Pic^-(X) & \longto & \Pic^+(X_{\rm reg}).\end{array}$$
By taking Cartier duals, we obtain the following commutative square.
$$\begin{array}{ccc} \Alb^-(\tilde X)&\longto &\Alb^+(\tilde X) \\
\uparrow & &\uparrow \\
\Alb^-(X_{\rm reg}) & \longto & \Alb^+(X).
\end{array}$$
In particular we get a canonical mapping
\begin{equation}\label{tau}
\tau^{-}_{+}:\Alb^-(X_{\rm reg})\to \Alb^+(X).
\end{equation}
Let $a_{\bf x}:\bar X\to \Alb(\bar X)$ be the Albanese mapping, obtained
by choosing a base point $x_c$ in each component $X^c$ of $X_{\rm reg}$
(note that the components of $X_{\rm reg}$ and $\bar X$ are in bijection).
Since $\Alb^-(X_{\rm reg})$ is a torus bundle on $\Alb (\bar X)$ we can
consider the following pull-back (\cf Section~\ref{cohomalb})
$$\begin{array}{ccc}
\Alb^{-}(X_{\rm reg}) &\to &\Alb (\bar X)\\
\bar a_{}\ \uparrow & &\mbox{\ \ } \uparrow a_{} \\
\Alb^T(X)&\to &\bar X
\end{array}$$
One can see that the torus bundle $\Alb^T(X)\to \bar X$ is
trivial when restricted to the open subset $X_{\rm reg}\subset \bar X$ (in
fact, the same argument in Section~\ref{cohomalb} applies here, \cf
\cite[\S1]{SER}). Hence we get a section $\sigma_{}^{\rm reg}: X_{\rm
reg}\to
\Alb^T(X)$. By composing $\sigma_{}^{\rm reg}$ and $\bar a_{}$
we
get the Albanese  mapping
\begin{equation} \label{regalbmap}
a^-_{\bf x}: X_{\rm reg}\to \Alb^-(X_{\rm reg}).
\end{equation}

\begin{propose}\label{serre}
For any equidimensional variety $X$ over $k=\bar k$, the morphism
$$a^-_{\bf x}: X_{\rm reg}\to \Alb^-(X_{\rm reg})$$
is universal among (base point preserving) morphisms to semi-abelian
varieties, in the sense of Serre \cite{SEM}.
If $X$ is a {\rm normal} proper $k$-variety, then $\Alb^-(X_{\rm reg})$ is
the Albanese variety of any resolution of singularities of $X$.
 \end{propose}
\begin{proof} It follows from the explicit construction by Serre in
\cite{SER} that $\Alb^-(X_{\rm reg})$ is equal to its ``Albanese variety''
in the sense of \cite{SEM},
and the morphism $a_{\bf x}^-$ is then universal by \cite[Th\'eor\`eme
1]{SER}, \ie  any torus bundle on $\Alb (\bar X)$ which is trivial on
$X_{\rm reg}$ is a push-out of $\Alb^-(X_{\rm reg})$.
If moreover $X$ is normal and proper, and $\tilde{X}\to X$ is a projective
resolution of singularities, then $\Div^0_{\tilde{S}}(\tilde{X})=0$, as
seen in the proof of Proposition~\ref{normplus}; therefore the character
group of the torus vanishes.
\end{proof}

If $X$ is proper then $\Alb^+(X)$ is semi-abelian and the Albanese map
$a^+_{\bf x}:X_{\rm reg} \to \Alb^+(X)$ defined in (\ref{plusalbmap}) can
be obtained by composing  $a^-_{\bf x}$ and $\tau^{-}_{+}$ defined in
(\ref{tau}). Since $\Alb^-(X_{\rm reg})$ is universal, $\tau^{-}_{+}$ can
be also be regarded as being induced by the universal property (note that
$\tau^{-}_{+}$ is affine and surjective).
\begin{propose}\label{reg}
Let $X$ be proper over $k$. Then there is an extension
$$0\to T(S)\to \Alb^-(X_{\rm reg})\to \Alb^+(X)\to 0$$
with kernel the torus $T(S)$ whose character group is the
quotient lattice
$$\frac{\Div_{\tilde S}^0(\tilde X)}{\Div_{\tilde S/S}^0(\tilde X)}.$$
This is a sub-lattice of the lattice of divisors on $X$ which are supported on
the singular locus $S$ of $X$: in particular, $T(S)=0$ if $X$ is
non-singular in codimension  one.
\end{propose}
\begin{proof} Since $\bar X =\tilde X$, the claimed torus bundle is obtained as
the Cartier dual of the following injective map of 1-motives
$$[\Div_{\tilde S/S}^0(\tilde X)\to\Pic^0(\tilde X)]\to
[\Div_{\tilde S}^0(\bar X)\to\Pic^0(\tilde X)] $$
Since $\Div_{\tilde S/S}(\tilde X) \df \ker (\Div_{\tilde S}(\tilde X)
\to \Div_{S}(X))$ the description of $T(S)$ is clear.
\end{proof}

\subsection{Hodge, \'etale and De Rham realizations of $\Alb^-$}
An immediate consequence of the Theorem~\ref{plus} is the following.
\begin{cor} Let $X$ be defined over $\C$. Then
$$T_{Hodge}(\Alb^-(X))\cong H_{1}(X,\Z)/{\rm (torsion)}$$
\end{cor}
\begin{proof} It follows from Cartier duality and the isomorphism of mixed
Hodge
structures $$\Hom(H^{1}(X,\Z(1)),\Z(1))\cong H_{1}(X,\Z)/{\rm
(torsion)}$$ because of Theorem~\ref{plus}.
 \end{proof}

As a consequence of Theorems~\ref{ellplus} and~\ref{deplus} we have:
\begin{cor} Let $X$ be defined over an algebraically closed field of
characteristic zero. Then
$$\hat{T}(\Alb^-(X))\cong H_{1}^{\et}(X,\hat{\Z})/{\rm (torsion)}$$
and
$$T_{DR}(\Alb^-(X))\cong H_{1}^{DR}(X).$$
\end{cor}

We then have the following corollary deduced from the properties
of $\Alb^-$ and $\Alb^+$ obtained so far (see Proposition~\ref{reg} and
Lemma~\ref{surj}). Of course this may also be proved directly by
topological arguments (and is in fact well known to experts).
\begin{cor}
Let $X$ be a {\em normal} proper $k$-variety. Then
$$H_1^{\et}(X_{\rm reg},\hat{\Z})/{\rm (torsion)}\cong H_1^{\et}(\tilde
X,\hat{\Z})/{\rm (torsion)}$$
and $H_{1}^{\et}(X,\hat{\Z})/{\rm (torsion)}$ is a quotient of
$H_1^{\et}(\tilde X,\hat{\Z})$.

If $k=\C$ then
$$H_1(X_{\rm reg},\Z)/{\rm (torsion)}\cong
H_1(\tilde X,\Z)/{\rm (torsion)}$$
are isomorphic Hodge structures, pure of weight $-1$, and
$H_{1}(X,\Z)/{\rm (torsion)}$ is a quotient Hodge structure of $H_1(\tilde
X,\Z)$.
\end{cor}

\section{Motivic Abel-Jacobi and Gysin maps}\label{geometry}

We now obtain some further properties of our Albanese and Picard
1-motives. We will give algebro-geometric (= ``motivic'') constructions of
some cohomological operations.

\subsection{Functoriality}
Contravariant (resp. covariant) functoriality of $\Pic^+$ (resp. $\Alb^-$)
is true, essentially by construction, and is valid for every morphism.
For $\Pic^-$ (resp. $\Alb^+$) we can expect covariance (resp.
contravariance) only for morphisms between varieties of the same
dimension, yielding the zero map if the morphism does not have dense image
in some irreducible component. We will work throughout over an
algebraically closed base field $k$ of characteristic 0.

\begin{propose} \label{pushminus}
Let $f:X\to X'$ be any morphism between $k$-varieties such that $\dim
X'=\dim X$. We then have a push-forward $f_*:\Pic^-(X)\to \Pic^-(X')$ and,
dually, a pull-back $f^*:\Alb^+(X')\to \Alb^+(X)$.
\end{propose}
\begin{proof} We will assume, for simplicity of exposition, that $X$
and $X'$ are irreducible; we leave the necessary modifications (mainly
notational) for the general case to the reader.

If the morphism $f$ is not dominant we define $f_*$ to be the zero
homomorphism. If $f$ is dominant, we choose resolutions $\tilde X\to X$,
$\tilde X' \to X'$ and good compactifications $\tilde X\into \bar X$ and
$\tilde X' \into \bar X'$ with normal crossing boundaries $Y\subset \bar X$,
$Y'\subset \bar X'$, such that there is a morphism $\bar f :\bar X\to \bar
X'$ compatible with $f$, and hence satisfying $\bar{f}^{-1}(Y')\subset Y$.

Let $D\in \Div_{\bar S/S}(\bar X, Y)$. The push-forward $f_*(D)$, as a Weil
divisor, clearly belongs to $\Div_{\bar S'/S'}(\bar X', Y')$.  We therefore
just need to show that there is an induced push-forward of relative line
bundles which is compatible with the push-forward of Weil divisors. This is the
content of the following lemma.\end{proof}

\begin{lemma} Let $f:\bar{X}\to \bar{X'}$ be a proper surjective morphism
between $n$-dimensional integral smooth proper varieties over an  algebraically
closed field $k$ of characteristic 0. Let  $\partial\bar{X}\subset\bar{X}$ and
$\partial\bar{X'}\subset\bar{X'}$  be reduced, normal crossing divisors such
that  $f^{-1}(\partial\bar{X'})_{\rm red}$ is a  normal crossing divisor in
$\bar{X}$ which is contained in  $\partial\bar{X}$.

Then there is a homomorphism of algebraic  groups
$$f_*:\Pic^0(\bar{X},\partial\bar{X})\to\Pic^0(\bar{X'},\partial\bar{X'})$$
such that
\begin{enumerate}
\item[(i)] $f_*$ is compatible with the natural
homomorphism  $f_*:\Pic^0(\bar{X})\to\Pic^0(\bar{X'})$ induced by the cycle
theoretic  direct image (push forward) on divisors
\item[(ii)] the assignment
$f\mapsto f_*$ is compatible with composition  of appropriate proper
maps
\item[(iii)] if $D$ is any divisor on $\bar{X}$ with support disjoint from
$\partial{\bar{X}}$, and $[D]\in \Pic^0(\bar{X},\partial\bar{X})$ is the  class
of the pair $(\cO_{\bar{X}}(D),s_D)$ (where $s_D$ is the tautological
meromorphic section of $\cO_{\bar{X}}(D)$ with divisor  $D$), then
$f_*[D]=[f_*D]\in \Pic(\bar{X'},\partial\bar{X'})$, where  $f_*D$ is
the cycle theoretic direct image (push forward) of $D$ under  the proper
map $f$, which is a divisor on $\bar{X'}$ whose support is  disjoint from
$\partial\bar{X'}$.
\end{enumerate}
\end{lemma}
\begin{proof}
By considering the obvious map $\Pic^0(\bar{X},\partial\bar{X})\to
\Pic^0(\bar{X},f^{-1}(\partial\bar{X'}))$,  we reduce immediately to the
case when  $\partial\bar{X}=f^{-1}(\partial\bar{X'})$. Now we can construct
a Stein factorization diagram
\[ \begin{TriCDV}
{\bar{X}}{\> f >>}{\bar{X'}}
{\SE h E E }{\NE E g E}{Y}
\end{TriCDV}\]
where $Y$ is a normal, proper variety of dimension $n$, $g:Y\to \bar{X'}$
is a finite, surjective morphism, and $h$ is birational and proper with
connected fibres. Further, $h_*\cO_{\bar{X}}=\cO_Y$, and
$f_*\cO_{\bar{X}}=g_*\cO_Y$.

Define
$\partial Y=g^{-1}(\partial\bar{X'})_{\rm red}$,
so that
$\partial\bar{X}=h^{-1}(\partial Y)_{\rm red}$.
Let $Z_1$ denote the union of the components of $Y_{sing}$ which are not
contained in $\partial Y$. Let $Z'=g(Z_1)$, $Z=g^{-1}(Z')$,
$\bar{Z}=h^{-1}(Z')$. Then $\bar{Z}$, $Z$ and $Z'$ are each closed subsets
of $\bar{X}$, $Y$ and $\bar{X'}$, respectively, which have codimension
$\geq 2$. Let $U=Y-Z$, $V=\bar{X}-\bar{Z}$, $W=\bar{X'}-Z'$, so that we
have an induced commutative triangle of proper morphisms
\[ \begin{TriCDV}
{V}{\> f >>}{W}
{\SE h E E }{\NE E g E}{U}
\end{TriCDV}\]
which is the Stein factorization of $f:V\to W$. Also define $\partial
V=V\cap \partial\bar{X}$, $\partial U=Y\cap\partial Y$, $\partial W=
W\cap\partial \bar{X'}$.

We now make the following claims.
\begin{enumerate}
\item[(i)] There is a homomorphism
$\alpha:\Pic^0(\bar{X},\partial\bar{X})\to\Pic(U,\partial U)$,
which fits into a commutative triangle
\[ \begin{TriCDV}
{\Pic^0(\bar{X},\partial\bar{X})}{\>{\rm restriction}>>}
{\Pic(V,\partial V)}{\SE \alpha E E }{\NE E h^* E}{\Pic(U,\partial U)}
\end{TriCDV}\]
\item[(ii)] There is a norm map $g_*:\Pic(U,\partial U)\to
\Pic(W,\partial W)$,
such that
(a) the composition $g_*\circ g^*$ is multiplication by $\deg g$, and
(b) $g_*[D]=[g_*D]$ for the class of any Weil ($=$ Cartier) divisor
$D$ on $U$ with support disjoint from $\partial U$.
\item[(iii)] The natural restriction map
$\rho :\Pic(\bar{X'},\partial\bar{X'})\to  \Pic(W,\partial W)$ is an
isomorphism.
\end{enumerate}

Granting these claims, the desired map $f_*$ is the composition
\[\Pic^0(\bar{X},\partial X)\by{\alpha}\Pic(U,\partial
U)\by{g_*}\Pic(W,\partial
W)\longby{(\rho)^{-1}}\Pic(\bar{X'},\partial\bar{X'}).\]
This obviously factors through the subgroup
$\Pic^0(\bar{X'},\partial\bar{X'})$,  which is the maximal divisible
subgroup of
$\Pic(\bar{X'},\partial\bar{X'})$.

We now proceed to prove the claims, in the order stated. First, we
consider the map $h^*:\Pic(U,\partial U)\to \Pic(V,\partial V)$. We have
that
\[\Pic(U,\partial U)=H^1(U,\cO^*_{(U,\partial U)}),\]
\[\Pic(V,\partial V)=H^1(V,\cO^*_{(V,\partial V)}),\]
where for a scheme $A$ and a closed subscheme $B$, we let
$\cO_{(A,B)}^*\df\ker(\cO_A^*\to\cO_B^*)$. By the Leray spectral
sequence for $h$, we obtain an exact sequence
\[0\to H^1(U,h_*\cO^*_{(V,\partial V)})\longby{h'}\Pic(V,\partial V)\to
H^0(U,R^1h_*\cO^*_{(V,\partial V)}),\]
and $h^*:\Pic(U,\partial U)\to \Pic(V,\partial V)$ is the composition of
$h'$ with the natural map
\[\Pic(U,\partial U)=H^1(U,\cO^*_{(U,\partial U)})\to
H^1(U,h_*\cO^*_{(V,\partial V)}).\]
In fact
\[\cO^*_{(U,\partial U)}=h_*\cO^*_{(V,\partial V)},\]
since $\cO^*_U=h_*\cO^*_V$, and the natural map $\cO^*_{\partial U}\to
h_*\cO^*_{\partial V}$ is injective. This means we have an exact sequence
\[0\to \Pic(U,\partial U)\by{h^*}\Pic(V,\partial V)\to
H^0(U,R^1h_*\cO^*_{(V,\partial V)}).\]
So to construct the map $\alpha$ and the commutative triangle in
Claim~(i), it suffices to prove that the natural map
\[\Pic^0(\bar{X},\partial \bar{X})\to H^0(U,R^1h_*\cO^*_{(V,\partial V)})\]
vanishes.

Thus it suffices to prove that for each closed point $x\in U$,
the map to the stalk at $x$
\[\Pic^0(\bar{X},\partial \bar{X})\to (R^1h_*\cO^*_{(V,\partial V)})_x\]
vanishes. We may identify $(R^1h_*\cO^*_{(V,\partial V)})_x$ with
$H^1(V_x,\cO^*_{(V_x,\partial V_x)})$, where
\[V_x=V\times_U\Spec \cO_{x,U},\;\;\partial V_x=\partial
V\times_U\Spec\cO_{x,U}.\]
So we want to show that the maps
\[\Pic^0(\bar{X},\partial \bar{X})\to \Pic(V_x,\partial V_x)\]
vanish, for all $x\in U$.

If $x\nin\partial U$, then $\partial V_x=\emptyset$. It suffices to see
that the natural map $\Pic^0(\bar{X})\to\Pic(V_x)$ vanishes. Now $x$ is a
non-singular point of $U$. Thus we can find a non-singular proper variety
$\bar{U}$, containing $U_{\rm reg}$ as a dense open subset; we can find a
non-singular proper variety $\bar{V}$ containing $h^{-1}(U_{\rm reg})$ as
a dense open set, and dominating $\bar{X}$. Then
$\Pic^0(\bar{X})\cong\Pic^0(\bar{V})\cong\Pic^0(\bar{U})$,
and evidently the map
\[\Pic^0(\bar{U})\to\Pic(V_x)\]
vanishes, as it factors through $\Pic(\Spec\cO_{x,U})=0$.

So we may take $x\in\partial U$. Now the fiber $h^{-1}(x)$ is contained
in $\partial V$. Let $\hat{\cO}_{x,U}$ be the completion of $\cO_{x,U}$,
and let
\[\hat{V_x}=V\times_U\Spec \hat{\cO}_{x,U},\;\;\partial\hat{V_x}=\partial
V\times_U\Spec\hat{\cO}_{x,U}.\]
Then we have a natural homomorphism
\[\Pic(V_x,\partial V_x)\to\Pic(\hat{V_x},\partial \hat{V_x}).\]
Since $\cO_{x,U}\to \hat{\cO}_{x,U}$ is faithfully flat, we see easily
that this homomorphism is {\em injective}. So we are reduced to proving
that $\Pic^0(\bar{X},\partial\bar{X})\to
\Pic(\hat{V_x},\partial\hat{V_x})$ vanishes.

For each $n\geq 1$, let $V_x^n\subset V_x$ be the closed subscheme
defined by the $n$-th power of the ideal sheaf of the reduced fiber
$h^{-1}(x)_{\rm red}$. Let $\partial V_x^n$ denote the scheme
theoretic intersection $\partial V\cap V_x^n$. Then $V_x^1=\partial
V_x^1=h^{-1}(x)_{\rm red}$, since $h^{-1}(x)\subset\partial V$. There
is a natural homomorphism
\[\Pic(\hat{V_x},\partial\hat{V_x})\to\liminv{n}\Pic(V_x^n,\partial
V_x^n).\]
We claim that it is an isomorphism. This follows, using the five
lemma, from the Grothendieck Existence Theorem~\cite{SGA2}, which gives
isomorphisms
\[\Pic(\hat{V_x})\cong\liminv{n}\Pic(V_x^n),\;\;
\Pic(\partial\hat{V_x})\cong\liminv{n}\Pic(\partial V^n_x),\]
and analogous isomorphisms on unit groups.

Hence we are reduced to proving that for each $n$, the natural
restriction maps
\[\Pic^0(\bar{X},\partial\bar{X})\to\Pic(V_x^n,\partial V_x^n)\]
are zero. This is clear for $n=1$ since $V_x^1=\partial V_x^1$, so that
$\Pic(V_x^1,\partial V_x^1)=0$. For $n>1$, one has that
\[\Pic(V_x^n,\partial V_x^n)=\ker\left(\Pic(V_x^n,\partial V_x^n)\to
\Pic(V_x^1,\partial V_x^1)\right)\]
is an affine algebraic group which is purely of additive type (\ie is a
vector group) \cite[Section~4]{BOU}. Hence any homomorphism from a semi-abelian
variety to $\Pic(V_x^n,\partial V_x^n)$ must vanish. This completes the
proof of Claim~(i).

Now we construct the norm map $g_*:\Pic(U,\partial U)\to\Pic(W, \partial
W)$ of Claim~(ii). First note that $R^1g_*\cO^*_{(U,\partial U)}=0$, since
the relative Picard group of a semi-local pair vanishes. Hence we have an
identification
\[\Pic(U,\partial U)=H^1(W,g_*\cO^*_{(U,\partial U)}).\]
Since $U,W$ are integral and normal, and $g$ is finite surjective, the norm
map on functions induces a homomorphism $N_{U/W}:g_*\cO_U^*\to \cO_W^*$. We
claim this induces a map on subsheaves $N_{U/W}:g_*\cO^*_{(U,\partial
U)}\to \cO^*_{(W,\partial W)}$, or equivalently, that the composition
\[g_*\cO^*_{(U,\partial U)}\into
g_*\cO^*_U\longby{N_{U/W}}\cO^*_W\to\cO^*_{\partial W}\]
vanishes. Since $\cO^*_{\partial W}$ injects into the direct sum of
constant sheaves associated to its stalks at the generic points of
$\partial W$, it suffices to show that for any such generic point
$\eta\in\partial W$, the map on stalks
\[(g_*\cO^*_{(U,\partial U)})_{\eta}\to \cO^*_{\eta,\partial W}\]
vanishes. Now $\cO_{\eta,\partial W}$ is the function field of an
irreducible component of $\partial W$, and is the residue field of the
discrete valuation ring $\cO_{\eta,W}$. The stalk
$(g_*\cO_U^*)_{\eta}$ is the unit group of the (semi-local) integral closure of
$\cO_{\eta,W}$ in the function field of $U$; denote this semi-local ring by
$\cO_{\eta,U}$.  The stalk $(g_*\cO^*_{(U,\partial U)})_{\eta}$ is the
subgroup of $\cO^*_{\eta ,U}$ of units congruent to $1$ modulo the
Jacobson radical (which is the ideal defining $\partial
U=g^{-1}(\partial W)_{\rm red}$ in $\cO_{\eta,U}$).

Now $\cO_{\eta,U}$ is a free module over $\cO_{\eta,W}$ of rank equal to
the degree of $g$, and for any $a\in\cO_{\eta,U}$, the norm of $a$ equals
the determinant of the endomorphism of the free $\cO_{\eta,W}$-module
$\cO_{\eta,U}$ given by multiplication by $a$. So it suffices to observe
that if $a\in \cO^*_{\eta ,U}$ is congruent to $1$ modulo the Jacobson
radical, then this endomorphism is of the form $1+A$, where the matrix
entries of $A$ lie in the maximal ideal of $\cO_{\eta,W}$; hence the
determinant of this matrix maps to $1$ in the residue field of
$\cO_{\eta,W}$. This proves that
\[(g_*\cO^*_{(U,\partial U)})_{\eta}\to \cO^*_{\eta,\partial W}\]
vanishes.

Now we define the map
\[g_*:\Pic(U,\partial U)\to \Pic(W,\partial W)\]
to be the map
\[H^1(W,g_*\cO^*_{(U,\partial U)})\to H^1(W,\cO^*_{(W,\partial W)})\]
induced by the sheaf map
\[N_{U/W}:g_*\cO^*_{(U,\partial U)}\to \cO^*_{(W,\partial W)}.\]
This evidently has the property that $g_*\circ  g^*$ is multiplication by
$\deg g$, since this is true at the sheaf level. To see the
compatibility with the push-forward for divisors $D$ with support
$\mod{D}$ disjoint from $\partial U$, we
compare the above map $g_*$ with the analogous map
\[H^1_{g(\mid{D}\mid)}(W,g_*\cO^*_{(U,\partial U)})\to
H^1_{g(\mid{D}\mid)}(W,\cO^*_{(W,\partial W)}).\]
 This completes the proof of Claim~(ii).

To prove Claim~(iii), it suffices to note that
\[\Pic(\bar{X'})\cong\Pic(W),\;\;\Pic(\partial \bar{X'})\into\Pic(\partial
W),\] and that
\[H^0(\bar{X'},\cO^*_{\bar{X'}})\cong H^0(W,\cO^*_W),\;\;H^0(\partial
\bar{X'},\cO^*_{\partial \bar{X'}})\cong H^0(\partial W,\cO^*_{\partial
W}).\]
All of these follow from the choice of the open set $W\subset \bar{X'}$,
such that $\bar{X'}-W$ has codimension $\geq 2$ in $\bar{X'}$, and $\partial
\bar{X'}-\partial W$ has codimension $\geq 2$ in $\partial \bar{X'}$ (recall
that $\bar{X'}$ is integral, non-singular and complete, and $\partial
\bar{X'}$ is a reduced, normal crossing divisor in $\bar{X'}$, and is
hence a complete, equidimensional and Cohen-Macaulay scheme; thus
$\bar{X'}$ is locally connected in codimension 2).
\end{proof}

\begin{rmk}{\rm  By making use of Propositions~\ref{pushminus}
and~\ref{fhat} we
can see the following faithfulness property of the Albanese and Picard
1-motives.
If we let $f:X\to X'$ be a generically finite morphism such that
the push-forward $f_*:\Pic^-(X)\to\Pic^-(X')$ induces an isomorphism
on \'etale realizations then $f_*$ itself is an isomorphism of 1-motives.
A similar statement holds for $f^*:\Alb^+(X')\to \Alb^+(X)$.}
\end{rmk}

\subsection{Projective bundles and vector bundles}
Let $P=\P (\cE)={\bf Proj}\,S(\cE)$ be the projective bundle associated to
a locally free sheaf $\cE$ on $X$ (here $S(\cE)$ is the
symmetric algebra of $\cE$ over $\cO_X)$.
\begin{propose}\label{proj} There are canonical
isomorphisms $\Pic^-(X)\cong \Pic^-(P)$ and $\Pic^+(X)\cong \Pic^+(P)$,
therefore, dually, $\Alb^+(X)\cong \Alb^+(P)$ and $\Alb^-(X)\cong
\Alb^-(P)$.  \end{propose}
\begin{proof} Let $\tilde P=\P (\tilde \cE)\to \tilde X$ be the pullback
along a choosen resolution of singularities $\tilde X \to X$. We can
choose a ``Nash compactification'' $\bar X$ of the resolution $\tilde X$,
\ie we can also get a locally free sheaf $\bar \cE$ on $\bar X$ which
extends $\tilde \cE$ (to construct a Nash compactification, first choose
an arbitrary one, and a coherent extension $\cF$ of $\tilde \cE$; then
resolve singularities of the Nash blow-up associated to $\cF$, on which
the pull-back of $\cF$, modulo torsion, is a locally free sheaf).

We can then assume that $\tilde P$ extends to $\bar P =\P(\bar \cE)$ on
$\bar X$, and the boundary $\bar P-\tilde P=Z$ is a normal crossing
divisor in $\bar P$ which is a projective bundle over the normal crossing
boundary $Y$ of $\bar X$. Since the Picard varieties of $\bar X$ and $\bar
P$ are also isomorphic, the exact sequence (\ref{longpic}) (\cf
Proposition~\ref{relpic}) yields an isomorphism of semi-abelian
varieties $\Pic^0(\bar X,Y)\cong \Pic^0(\bar P,Z)$. Pull-back of divisors
from $\bar X$ to $\bar P$ yields a compatible isomorphism between
lattices, giving rise to the claimed isomorphism for $\Pic^-$; that for
$\Alb^+$ follows from Cartier duality.

For $\Pic^+$ and $\Alb^-$, we argue as follows. Consider a Nash
compactification $X^*$ of $X$, \ie such that $\cE$ extends to a locally
free sheaf $\cE^*$ on $X^*$. We can find a smooth proper hypercovering
$\bar X_{\d}$ of $X^*$ such that the induced reduced hypercovering of
$X^*-X$ is a normal crossing divisor $Y_{\d}$ in $\bar X_{\d}$. Then
$X_{\d}=\bar X_{\d}-Y_{\d}$ yields a smooth proper hypercovering of $X$,
and $\bar X_{\d}$ is a smooth compactification with normal crossing
boundary. Now $P_{\d}=P\times_X X_{\d}$ is a smooth proper hypercovering
of $P$, and we can get an induced compactification of $P_{\d}$
\[\bar P_{\d}=\bar X_{\d}\times_{X^*} \P (\cE^*)\]
which has normal crossing boundary. We then see easily that $\bPic^0(\bar
X_{\d})\cong \bPic^0(\bar P_{\d})$ because of the exact sequence
(\ref{semisimp}) (\cf Proposition~\ref{simpic}); similarly the lattices
are isomorphic.
\end{proof}

Let $V=\V(\cE)={\bf Spec}\,S(\cE)$ be the geometric vector bundle
associated to a locally free sheaf $\cE$ on $X$. We have the following
homotopy invariance property.
\begin{propose} There is a canonical isomorphism
$\Pic^+(X)\cong \Pic^+(V)$ and, dually, there is an isomorphism
$\Alb^-(X)\cong \Alb^-(V)$.
\end{propose} \begin{proof} Consider a Nash compactification $X^*$
of $X$, so that $\cE$ extends to a locally free sheaf $\cE^*$ on
$X^*$, and let $V^*=\V(\cE^*)$. We let $\bar X_{\d}$ be a smooth proper
hypercovering of $X^*$ such that the reduced inverse image of $X^*-X$
is a normal crossing divisor, and let $V_{\d}^*$ be the simplicial vector
bundle on $\bar X_{\d}$ obtained by the pull-back of $V^*$ along the
hypercovering. We take
$$\bar V_{\d} = \bar X_{\d}\times_{X^*}\P(\cE^*\oplus\cO_{X^*})$$
to be the compactification of $V_{\d}^*$ with normal crossing boundary.
We then have to show that
$$[\Div_{Y_{\d}}^0(\bar X_{\d})\to \bPic^0(\bar X_{\d})]\cong
[\Div_{N_{\d}}^0(\bar V_{\d})\to \bPic^0(\bar V_{\d})]$$
where $N_{\d}$ is the normal crossing boundary of $\bar V_{\d}$,
considered as a compactification of $V_{\d}= X_{\d}\times_XV$. We have
$$N_{\d}=\bar V_{\d} - V_{\d}= Y_{\d}\times_X\P(\cE^*\oplus\cO_{X^*})\cup
\bar X_{\d}\times_{X^*}\P(\cE^*).$$
Thus it is clear that the groups of divisors supported on $N_{\d}$ and
on $Y_{\d}$, which are algebraically equivalent to zero (\ie have
classes in $\bPic^0$) on the respective proper simplicial schemes, are
naturally isomorphic; hence the lattices of our two 1-motives are
naturally isomorphic.  From the short exact sequence
$$0\to \bPic (\bar X_{\d})\to \bPic (\bar V_{\d})\to \Z\to 0$$
we conclude that $\bPic^0(\bar X_{\d})\cong \bPic^0(\bar V_{\d})$,
and we are done.
\end{proof}

\subsection{Universality and zero-cycles} We let $X$ be a {\it projective}
$n$-dimensional $k$-variety. Let $X^{(n)}$ be the union of the
$n$-dimensional irreducible components of $X$, and let $X_{\rm reg}$
denote
the locus of smooth points of $X$ which lie in $X^{(n)}$. We fix base
points $x_c\in X^c$ in each component of $X_{\rm reg}$, and let
$a_{\bf x}^+:X_{\rm reg}\to \Alb^+(X)$ be the corresponding Albanese map
(see
(\ref{plusalbmap})). We denote by
$$a^+_X:\cZ^n(X_{\rm reg})_{\deg 0}\to \Alb^+(X)$$
the induced map on the group $\cZ^n(X_{\rm reg})_{\deg 0}$ of zero
cycles on $X_{\rm reg}$ which have degree 0 on each component of
$X_{\rm reg}$; in
fact $a^+_X$ is independent of the choices of base points
$\{x_c\}={\bf x}$.

We recall that the ``cohomological'' Levine-Weibel Chow group of
zero-cycles $CH^n(X)$ is defined to be the quotient of the free abelian
group on (closed) points of $X_{\rm reg}$, \ie $\cZ^n(X_{\rm reg})$,
modulo the subgroup of zero-cycles which are divisors of appropriate
rational functions on Cartier curves on $X$ (\cf \cite{LW} and \cite{BS}).

Using $a^+_X$, we get a ``motivic'' construction  of an Abel-Jacobi
map, generalizing the Abel-Jacobi map for the Chow group of zero-cycles
of degree 0 on projective non-singular varieties, to the case of
projective varieties with arbitrary singularities (this is
done in \cite{BPW} and \cite{BS} over $\C$). For a different
algebraic construction, see \cite{ESV}.
\begin{thm} \label{universal} Let $X$ be a projective $k$-variety.  The
Albanese map $a^+_X$ yields a universal regular homomorphism
\begin{equation} a^+: CH^{n}(X)_{\deg 0}\to \Alb^+(X)
\end{equation}
from the ``cohomological'' Chow group of zero-cycles of degree zero to
semi-abelian $k$-varieties.
\end{thm}
We will prove the above theorem in several steps. We first
construct Gysin maps for ``good'' curves, defined as follows. A curve
$C\subset X$ is ``good'' if (i) $C$ is reduced, purely of dimension 1, and
$C\cap S$ is reduced of dimension 0 (ii) $C$ is a local complete
intersection in $X$ (iii) $C\subset X^{(n)}$
(iv)  if $X_n\to X$ is the normalization, and we set $C_0=X_n\times_XC$,
$S_0=X_n\times_X(C\cap S)$, then $C_0$ is also purely 1-dimensional, and
$S_0\subset C_0$ consists of smooth points of $C_0$. Note that if $C$ is
good, and $X^{<n}$ is the union of the irreducible
components of $X$ of dimension $<n$, then $C\cap X^{<n}=\emptyset$.

\begin{lemma} \label{good}
Let $i:C\into X$ be a ``good'' curve in $X$. We then have Gysin maps
$$i^*_{-}:\Pic^-(X)\to \Pic^-(C)$$
and dually
$$i_*^+:\Alb^+(C)\to \Alb^+(X).$$
 \end{lemma}
\begin{proof} We may assume without loss of generality that $X$ is
equidimensional. Let $X_n\to X$ be the normalization, $\tilde{X}\to X_n$ a
resolution of singularities, and $f:\tilde X\to X$ the induced resolution
of singularities. Since $C$ is ``good'', the scheme $C\times_X\tilde X$ is
a curve which is smooth at $f^{-1}(C\cap S)$.

Denote by $C'$ the pull-back curve $\tilde X\times_X C$. Let
$f':C'\to C$ be the restriction of $f$. Then $C'\cong C_0=X_n\times_XC$,
and the normalisation $\tilde C\to C$ of the curve $C$ clearly factors
through $f':C'\to C$. Let $\tilde i:\tilde C\to \tilde X$ be the induced
map. Then there is a natural pull-back map on Picard varieties
$\tilde i^*:\Pic^0(\tilde X)\to\Pic^0(\tilde C)$.
Thus, in order to get the claimed map $i^*_{-}$ on 1-motives, it is enough
to show that any divisor $D\in\Div_{\tilde S/S}^0(\tilde X)$ pulls back to
a divisor $\tilde i^*(D)\in\Div_{\tilde C/C}^0(\tilde C)$. Since $C'$
is smooth at the finite set of points $f^{-1}(C\cap S)$ it will suffices
to show that  $(i')^*(D)\in\Div_{ C'/C}^0(C')$ where $i':C'\into \tilde X$
is the canonical induced imbedding.

Now let $\tilde D$ denote the support of $D$. Then $\tilde D$ is mapped to
$S$, and therefore $\tilde D\times_X C$, which is the support of
$(i')^*(D)$, is mapped to $C\cap S$. We thus have the following diagram of
Fulton's homological Chow groups
 $$\begin{array}{ccccc}
CH_{n-1}(\tilde D)&\to & CH_{n-1}(S)&\to & CH_{n-1}(X)\\
\downarrow & & & &\downarrow\\
CH_{0}(\tilde D\times_X C)&\to & CH_{0}(C\cap S)& \into &CH_{0}(C)
\end{array}$$
by Fulton's compatibility result \cite[Theorem~6.1]{FU} between pull-back and
Gysin maps for locally complete intersection morphisms. Since the
push-forward of $D$ vanishes as a cycle on $S$, the pull-back of $D$
to $\tilde D\times_X C$ pushes forward to zero in $CH_{0}(C\cap S)$.
Since $C\cap S$ is a reduced 0-dimensional scheme, the latter push-forward
to $CH_{0}(C\cap S)$ is in fact zero as a cycle on $C\cap S$.
\end{proof}

We need the following compatibilities (\cf Lemma~3.3-3.4 in \cite{BS}).
\begin{lemma}\label{ajcomp}
(a) Let $C$ be a ``good'' curve as in Lemma~\ref{good}. There is a
commutative diagram
 $$\begin{array}{ccc} \cZ^1(C_{\rm reg})_{\deg 0} &\longby{a^+_C}
&\Alb^+(C)\\
i_*\downarrow & & \downarrow i^+_*\\
 \cZ^n(X_{\rm reg})_{\deg 0} &\longby{a^+_X} &\Alb^+(X).
\end{array}$$
(b) Let $f:Y\to X$ be a morphism of $n$-dimensional projective
varieties, such that $f_{\mid Y_{\rm reg}}: Y_{\rm reg}\to X_{\rm reg}$ is a
finite, flat morphism. Let $\gamma\in\cZ^n(X_{\rm reg})_{\deg 0}$ be a
zero-cycle of degree zero with inverse image
$f^*(\gamma)\in\cZ^n(Y_{\rm reg})_{\deg 0}$. We then have
\begin{equation}\label{ajpb} a^+_Y(f^*(\gamma)) =f^*(a^+_X(\gamma)).
\end{equation}
(c) If $f:Y\to X$ is a blow up at a smooth point of $X$ there is a
commutative diagram
$$\begin{array}{ccc} \cZ^n(Y_{\rm reg})_{\deg 0}&\longby{a^+_Y}
&\Alb^+(Y)\\
f_*\downarrow & &\cong \uparrow f^*\\
 \cZ^n(X_{\rm reg}) &\longby{a^+_X} &\Alb^+(X)
\end{array}$$
\end{lemma}
The proof is left as an exercise for the reader.

\begin{lemma}\label{ajcurve} Let $C$ be a reduced projective curve. The
canonical section $C_{\rm reg}\to \Alb^+(C)$ yields a universal
regular homomorphism to semi-abelian $k$-varieties
$$a^+:CH^1(C)_{\deg 0}\to \Alb^+(C),$$
which is an isomorphism when $C$ is seminormal.
\end{lemma}\begin{proof}
We recall that $CH^1(C)\cong \Pic(C)$ and $CH^1(C)_{\deg
0}\cong\Pic^0(C)$. Let $C'$ be the semi-normalization
of $C$; the canonical identification $\Pic^0(C')\cong\Alb^+(C)$ (see
Proposition~\ref{delmot}) together with the pull-back map $\Pic^0(C)\to
\Pic^0(C')$, which is just the semi-abelian quotient of $\Pic^0(C)$,
yields the result.
\end{proof}

Now, in order to show that the map $a^+:\cZ^n (X_{\rm reg})_{\deg 0}\to
\Alb^+(X)$ ($n>1$) factors through rational equivalence, by \cite{BS}, it
suffices to show that $\ker a^+$ contains all divisors $(f)_C$ where: {\it
i)}\, $C$ is a ``good'' curve in $X$, and {\it ii)}\, $f$ is a rational
function on $C$ which is a unit at points of $C\cap S$. Using our
Lemma~\ref{ajcomp} we adapt the proof of Lemma~3.5 in \cite{BS} to our
situation.

In order to show universality of $a^+:CH^n(X)_{\deg 0}\to \Alb^+(X)$, we
first note that, from the definitions, it is easy to see that $a^+$
factors through the natural surjection $CH^n(X)_{\deg 0}\to
CH^n(X^{(n)})_{\deg 0}$, since by definition $\Alb^+(X)=\Alb^+(X^{(n)})$.
So we may assume $X$ is equidimensional.

Now consider the canonical
extension
$$0\to T(S)\to \Alb^-(X_{\rm reg})\to \Alb^+(X)\to 0$$
(see Proposition~\ref{reg}). If $\psi : CH^{n}(X)_{\deg 0}\to G$ is a
regular homomorphism to a semi-abelian variety $G$, we need to find a
unique factorisation
$$\begin{array}{ccc}
CH^{n}(X)_{\deg 0}& \by{a^+}&\Alb^+(X)\\
&&\mbox{\ \ \ }\downarrow\psi^+\\
&&G
\end{array}$$
through $a^+$, for some homomorphism of algebraic groups $\psi^+$. Since
$X_{\rm reg}=\coprod_cU^c$ maps to $CH^{n}(X)_{\deg 0}$ by taking a point
$x\in U^c$ to the difference $x-x_c$ in the Chow group, we get a
map $\psi_0:X_{\rm reg}\to G$. By definition, since $\psi$ is a regular
homomorphism, $\psi_0$ is a morphism, which sends each of the base points
$x_c$ to  0. By the universal property (Proposition~\ref{serre})
of $\Alb^-(X_{\rm reg})$, $\psi_0$ factors through $\Alb^-(X_{\rm reg})$
yielding a map $\psi^-:\Alb^-(X_{\rm reg})\to G$.
Using the above-mentioned canonical extension, we need to show that
$\psi^-(T(S))=0$ in order to obtain a well defined map $\psi^+$ on the
quotient semi-abelian variety $\Alb^+(X)$; the uniqueness of $\psi^-$
will then imply that of $\psi^+$.

We have the following fact.
\begin{lemma}\label{ajreg}
 Let $i:C\into X$ be a complete intersection curve in
$X$ which is ``good'' (\ie satisfies the hypoteses of Lemma~\ref{good}), such
that $C$ meets every irreducible component of $S$, and
moreover its singular locus $F$ is exactly $C\cap S$.
We have a commutative diagram
$$\begin{array}{ccccccc}
0\to & T(S)& \to & \Alb^-(X_{\rm reg})& \to & \Alb^+(X)&\to 0\\
&\uparrow & &\uparrow i_*& &\uparrow i_*^+& \\
0\to & T(F)& \to & \Alb^-(C_{\rm reg})& \to & \Alb^+(C)&\to 0
\end{array}$$
where $T(F)\to T(S)$ is a surjection of tori.
\end{lemma}
\begin{proof} This follows easily from the dual statement, \ie that the
following diagram
$$\begin{array}{ccccccc}
0\to &\Pic^-(X)& \to & \Pic^+(X_{\rm reg})& \to & \Div_S^0(X)&\to 0\\
&\downarrow & &\downarrow& &\downarrow & \\
0\to & \Pic^-(C)& \to & \Pic^+(C_{\rm reg})& \to & \Div_F^0(C)&\to 0
\end{array}$$
commutes. Moreover, $\Div_S(X)$ injects into $\Div_F(C)$.
\end{proof}

By successive hyperplane sections we can always find a general complete
intersection curve $C$ as above; therefore by Lemma~\ref{ajcurve} and
Lemma~\ref{ajreg} we conclude as follows. Since $i:C_{\rm reg} \into
X_{\rm reg}$, we have that the composite of the following
$$CH^1(C)_{\deg 0}\by{i_*}CH^n(X)_{\deg 0}\by{\psi} G$$
yields a unique
$$\psi^+_C:\Alb^+(C)\to G$$ by the universal property for curves, \ie
Lemma~\ref{ajcurve}; whence $\psi^+_C(T(F))=0$, because the universal
morphism $C_{\rm reg}\to\Alb^-(C_{\rm reg})$ is compatible with
$\psi^+_C$. Since $T(F)$ surjects onto $T(S)$, the commuativity of
the diagram in Lemma~\ref{ajreg} implies that $\psi^-(T(S))=0$ as claimed.
Thus Theorem~\ref{universal} is proved.

\subsection{Gysin maps}
First consider the case of normal varieties.
\begin{propose} Let $f:X' \to X$ be any proper morphism of $k$-varieties,
where $X$ is normal. We then have a functorial Gysin map
$$f^*_-:\Pic^-(X)\to\Pic^-(X')$$
and, dually,
$$f^+_*:\Alb^+(X')\to\Alb^+(X)$$
\end{propose}
\begin{proof} Let $\bar f:\bar X' \to \bar X$ be the induced map on smooth
compactifications $\bar X'$ and $\bar X$, compatibly with the normal
crossing boundaries $Y'$ and $Y$. We
then have the following diagram of 1-motives
$$\begin{array}{ccc} \Pic^-(X)& &\Pic^-(X')\\
\veq\downarrow & &\quad \uparrow \bar{f}^*\\
\Pic^-(\tilde X)&\to &\Pic^-(\tilde X')
\end{array}$$
yielding the claimed map, where since $X$ is normal we have that
$\Pic^-(X)\cong \Pic^-(\tilde X) =\Pic^0(\bar X, Y)$, and we have a
pull-back map on relative line bundles $\bar f^*:\Pic^0(\bar X, Y)\to
\Pic^0(\bar X', Y')$.
\end{proof}

We recall that a morphism $f:Z \to X$ is a {\em projective local complete
intersection morphism} if can be factorized as $f=\pi\circ i$ for a
regular imbedding $i: Z\into P$ and a projection $\pi:P\to X$ from
the projective bundle $P=\P(\cE)$ associated to a locally free
$\cO_X$-module $\cE$.
\begin{thm} Let $f:Z \to X$ be a projective local complete intersection
morphism. We then get a functorial Gysin map
$$f^+_*:\Alb^+(Z)\to\Alb^+(X)$$  and, dually,
$$f^*_-:\Pic^-(X)\to\Pic^-(Z)$$
\end{thm}
\begin{proof} Since $\Alb^+(P)\cong\Alb^+(X)$ by Proposition~\ref{proj} we
are left to prove our claim for regular imbeddings. We then have the
following diagram
$$\begin{array}{ccc}
CH^{n}(Z)_{\deg 0}& \by{a^+_Z}&\Alb^+(Z)\\
i_*\downarrow&&\\
CH^{n}(P)_{\deg 0}& \by{a^+_P}&\Alb^+(P)
\end{array}$$
where $i_*$ for cycles exists trivially, and therefore, by
Theorem~\ref{universal}, the composite of $i_*$ and $a^+_P$ factors through
$\Alb^+(Z)$. In order to show that the construction is independent
of the factorisation, we observe that it is so on the \'etale
realizations (where it coincides with the Gysin map obtained via
Grothendieck-Verdier duality), and therefore, by Proposition~\ref{fhat},
we are done.
\end{proof}
\begin{rmk}\label{remgysin}
{\rm It would of course be of interest to have the same result for proper
local complete intersection morphisms as well, for which the above
strategy of comparison with the Levine-Weibel Chow group of 0-cycles is
not applicable. It would also be desirable to have a ``geometric'' proof
of independence of the Gysin map from the choice of factorization,
instead of the above one using the \'etale realization.}
\end{rmk}

\section{Rationality Questions}\label{rational}
In this section, we consider the above theory in the case when the ground
field $k$ is an arbitrary field of characteristic 0. Let $\bar k$ denote
a fixed algebraic closure of $k$; if $A$ is any ``object''
(1-motive, scheme, morphism, sheaf ...) over $k$, then $A_{\bar{k}}$ will
denote its base change to $\bar{k}$.

First, consider a 1-motive $M=[L\by{u}G]$ over $k$. By definition, this is
a homomorphism between $k$-group schemes, where $L$ is an \'etale group
scheme, and $G$ a semi-abelian scheme, such that $L_{\bar k}$ is a lattice
(free abelian group of finite rank). The lattice $L_{\bar k}$ is naturally
a module over the Galois group $\Gal(\bar k/k)$, and the \'etale group
scheme $L$ is determined by this Galois module. The Galois group operates
semi-linearly on $G_{\bar k}=G\times_k\bar k$ as well, and the morphism
$u$ is determined uniquely by the morphism $u_{\bar k}:L_{\bar k}\to
G_{\bar k}$, which is $\Gal(\bar k/k)$-equivariant. Conversely, any Galois
equivariant morphism $L_{\bar k}\to G_{\bar k}$ is necessarily of the form
$u_{\bar k}$.

Thus, to give a 1-motive over $k$ is to give (i) a semi-abelian $k$-scheme
$G$ (ii) a lattice $\bar L$ which underlies a $\Gal(\bar k/k)$-module
(iii) a 1-motive $[\bar L\by{\bar u} G_{\bar k}]$ over $\bar k$, such that
$\bar u$ is $\Gal(\bar k/k)$-equivariant, for the given module structure
on $\bar L$, and the natural semi-linear action on $G_{\bar k}$.

If $k\into \C$, then for any 1-motive $M$ over $k$, we obtain a
corresponding 1-motive $M_{\sC}$ over $\C$, which has a Hodge realization.
For the \'etale realization, note that $\hat{T}(M_{\bar k})$ is a free
$\hat{\Z}$-module of finite rank, which supports a natural action of
$\Gal(\bar k/k)$. We call this Galois module the {\em \'etale
realization} of $M$. Finally, if $M^{\natural}$ denotes the universal
$\G_a$-extension of $M$ in the category of complexes of $k$-group schemes,
then  $M^{\natural}_{\bar k}$ is the universal $\G_a$-extension of
$M_{\bar k}$ in the category of complexes of $\bar k$-group schemes, and
\[\Lie(M^{\natural})_{\bar k}\cong \Lie(M^{\natural}_{\bar k})\]
as filtered $\bar{k}$-vector spaces. We define the De Rham realization
$T_{DR}(M)$ to be the filtered $k$-vector space $\Lie(M^{\natural})$.

The aim of this section is to show that if $X$ is a $k$-variety, then
there are naturally defined 1-motives $\Pic^+(X)$, $\Pic^-(X)$,
$\Alb^+(X)$, $\Alb^-(X)$ {\em defined over $k$}, pairwise Cartier
dual, with the following properties.
\begin{enumerate}
\item[(i)] If $k'$ is an extension field of $k$, the corresponding
1-motives for $X_{k'}$ are obtained by base change from $k$ to $k'$
from the 1-motives for $X$.
\item[(ii)] The \'etale realizations coincide with appropriate the
\'etale (co)homology groups (modulo torsion) of $X_{\bar k}$ as
$\Gal(\bar k/k)$-modules, where the Galois action on \'etale (co)homology
is the standard one.
\item[(iii)] The De Rham realizations coincide, as filtered $k$-vector
spaces, with the appropriate De Rham (co)homology groups of $X$ (defined
as in (\ref{derham}) via suitable hypercoverings and compactifications
over $k$).
\end{enumerate}

The proofs of the above assertions are fairly straightforward, and
basically amount to the observation that, when we carry out the
constructions of 1-motives for $X_{\bar k}$ as in the earlier sections,
and consider the computations of realizations, these are sufficiently
natural as to be automatically compatible with the action of $\Gal(\bar
k/k)$. As such, our arguments will be a little sketchy.

First consider the construction of $\Pic^-(X)$. Let $n=\dim X=\dim
X_{\bar k}$. If $X^{(n)}_{\bar k}$ is the union of the $n$-dimensional
irreducible components of $X_{\bar k}$, then it corresponds to a unique
closed $k$-subscheme $X^{(n)}$ of $X$, which is also purely of dimension
$n$. So we reduce to the case when $X$ and $X_{\bar k}$ are
equidimensional.

Now we may choose a resolution of singularities $f:\tilde{X}\to X$, and a
compactification $\bar X$ of $\tilde{X}$, both defined over $k$, such that
$\bar X_{\bar k}$ is a good normal crossing compactification of the
resolution $\tilde{X}_{\bar k}\to X_{\bar k}$. Let $Y\subset \bar X$ be
the normal crossing boundary divisor, $S\subset X$ the singular locus,
$\bar{S}\subset \bar X$ the Zariski closure of $f^{-1}(S)$.

Lemma~\ref{grpic} gives the representability of the relative Picard
functor of the pair $(\bar X,Y)$ by a $k$-group scheme (say, $\Pic(\bar
X,Y)$), locally of finite type, whose $\bar k$-points coincide with the
relative Picard group $\Pic(\bar X_{\bar k},Y_{\bar k})$. The identity
component $\Pic^0(\bar X_{\bar k},Y_{\bar k})$ is stable under the
semi-linear $\Gal(\bar k/k)$-action on $\Pic(\bar X,Y)_{\bar k}$, and so
naturally determines a $k$-subgroup scheme $\Pic^0(\bar X,Y)$ of
$\Pic(\bar X,Y)$. The lattice
$\Div^0_{\bar{S}_{\bar k}/S_{\bar k}}(\bar X_{\bar k},Y_{\bar k})$
is evidently stable under $\Gal(\bar k/k)$, with respect to the natural
Galois action on Weil divisors on $\bar X_{\bar k}$. Finally, the
canonical map
\[\Div^0_{\bar{S}_{\bar k}/S_{\bar k}}(\bar X_{\bar k},Y_{\bar k})\to
\Pic^0(\bar X_{\bar k},Y_{\bar k})\]
is clearly Galois equivariant. Hence we obtain a well-defined 1-motive
over $k$, which we define to be $\Pic^-(X)$; by construction we then have
$\Pic^-(X_{\bar k})=\Pic^-(X)_{\bar k}$.

The isomorphism
\[\hat{T}(\Pic^-(X_{\bar k}))\to H^{\et}_{2n-1}(X_{\bar
k},\hat{\Z}(1-n))/({\rm torsion})\]
is $\Gal(\bar k/k)$-equivariant, since it ultimately rests on the
identification, via Kummer theory, of \'etale $\mu_m$-coverings of
certain open subschemes of $X_{\bar k}$ with isomorphism classes of
triples $(\ccL,\varphi,\alpha)$ (see Proposition~\ref{relkum}) consisting
of $m$-torsion line bundles $\ccL$ with additonal trivializing data; but
this identification is easily seen to be Galois equivariant, where the
Galois group operates on such triples in the obvious way (corresponding
to the natural Galois action on $T_{\sZ/m}(\Pic^-(X_{\bar k}))$), while
it acts on the collection of \'etale coverings by twisting (changing the
structure morphism to $\Spec \bar k$), which corresponds to the natual
action on \'etale (co)homology.

As for the De Rham realization, the same proof that
$T_{DR}(\Pic^-(X_{\bar k}))=\Lie(\Pic^-(X_{\bar k})^{\natural})$
coincides, as a filtered $\bar k$-vector space with $H^{DR}_{2n-1}(X_{\bar
k})(1-n)$, yields a proof that $H^{DR}_{2n-1}(X)(1-n)$ coincides with
$\Lie(\Pic^-(X)^{\natural})$ as a filtered $k$-vector space, provided we
have 1 fact: that for $Z\subset \bar X$ as in Lemma~\ref{riglog}, the
Lie algebra computation in Lemma~\ref{riglog}(d) is valid.

This is of course clear over $\bar k$, from the formula in
Lemma~\ref{riglog}(a)
\[\Pic^{\natural-log}(\bar X_{\bar k},Y_{\bar k})\cong
\HH^1\left(\bar X_{\bar k}, \cO^*_{\bar X_{\bar k},Y_{\bar k}}
\longby{\dlog}
\Omega^1_{\bar X_{\bar k}}(\log(Y_{\bar k}+Z_{\bar k}))
(-Y_{\bar k})\right).\]
The analogous formula may not be valid over $k$, since the expression on
the right side arises as the value (on $k$, or $\bar k$) of an appropriate
Picard functor, while the left side refers to the sections of the
associated $fpqc$ sheaf (these do coincide over $\bar k$, while this is
unclear over $k$). But the tangent space at the identity to the Picard
functor admits a $k$-linear transformation to the corresponding tangent
space of the representable functor given by the associated $fpqc$-sheaf.
This linear transformation, upon base change to $\bar k$, is an
isomorphism of (filtered) vector spaces. Hence it is an isomorphism over
$k$ as well. So the ``presheaf tangent space'' is the same as the true
tangent space (this applies also to the simplicial Picard functor).

The results for $\Alb^+(X)$ now follow by Cartier duality from those for
$\Pic^-(X)$.

Next, consider $\Pic^+(X)$. We can choose a smooth proper hypercovering
$X_{\d}\to X$ and a smooth compactification $\bar X_{\d}$ of $X_{\d}$ with
normal crossing boundary $Y_{\d}$, all in the category of simplicial
$k$-schemes. The $fpqc$-sheaf associated to the simplicial Picard functor
of $\bar X_{\d}$ is representable by a $k$-group scheme, locally of finite
type, whose identity component is a semi-abelian $k$-scheme
$\bPic^0(\bar X_{\d})$, such that
\[\bPic^0(\bar X_{\d})_{\bar k}=
\bPic^0((\bar X_{\d})_{\bar k}).\]
The lattice $\Div^0_{(Y_{\d})_{\bar k}}((\bar X_{\d})_{\bar k})$ is a
Galois module in an obvious way, such that the map defining the 1-motive
$\Pic^+(X_{\bar k})$ is Galois equivariant. Thus there is a well-defined
1-motive $\Pic^+(X)$, defined over $k$, such that there is
an identification $\Pic^+(X)_{\bar k}=\Pic^+(X_{\bar k})$. The dicsussion
of the properties of the \'etale and De Rham realizations is very
similar to that of $\Pic^-(X)$, and is left to the reader. The case of
$\Alb^-(X)$ follows by Cartier duality.

—From the above constructions of 1-motives over $k$, it is also clear that
various functorial properties, and Albanese mappings, are all defined
over $k$ as well, possibly after slight reformulation. For example, if
$f:X\to X'$ is a morphism between $n$-dimensional $k$-varieties, then
there is a push-forward $f_*:\Pic^-(X)\to \Pic^-(X')$, and a pull-back
$f^*:\Alb^+(X')\to \Alb^+(X)$. For an arbitrary morphism of $k$-varieties,
there is a pull-back $f^*:\Pic^+(X')\to\Pic^+(X)$ and a push-forward
$\Alb^-(X)\to\Alb^-(X')$.

For the Albanese mappings $a^+_{\bf x}:X_{\rm reg}\to \Alb^+(X)$ and
$a^-_{\bf x}:X_{\rm reg}\to \Alb^-(X_{\rm reg})$, these
exist over $k$ provided we can choose the base points $x_c\in X^c\cap
(X_{\rm reg})_{\bar k}$ to be $k$-rational points, or more generally,
if the 0-cycle $\sum_c x_c$ is defined over $k$ (\ie is $\Gal(\bar
k/k)$-invariant). However, we cannot in general choose such such
base-points. Instead, we can consider the map
$$(A^-)_{\bar k}:\coprod_cX_c\times X_c\subset (X_{\rm reg}\times
X_{\rm reg})_{\bar k}\to  \Alb^-(X_{\rm reg})_{\bar k}$$
given by $(A^-)_{\bar k}(x,y)=a^-_{\bf x}(x)-a^-_{\bf x}(y)$, which is in
fact independent of ${\bf x}=\{x_c\}_c$, and $\Gal(\bar k/k)$-equivariant,
thus yielding a map of $k$-varieties
\[A^-:U\to \Alb^-(X_{\rm reg}),\]
where $U\subset X_{\rm reg}\times X_{\rm reg}$ is the open $k$-subscheme
consisting of the union of the connected components intersecting the
diagonal, and so $U_{\bar k}=\coprod_cX_c\times X_c$. The morphism $A^-$
is universal among those $k$-morphisms from $U$ to semi-abelian
$k$-varieties such that the natural involution on $U$ intertwines with
multiplication by $-1$ on the semi-abelian variety.

In a similar way, we can define a Galois equivariant $\bar k$-morphism
$$(A^+)_{\bar k}:\coprod_cX_c\times X_c\subset (X_{\rm reg}\times
X_{\rm reg})_{\bar k}\to \Alb^+(X)_{\bar k},$$
and hence a $k$-morphism
\[A^+:U\to \Alb^+(X).\]

\appendix
\section{Appendix: Picard functors}\label{rep}
Let $\cS ch_k$ be the category of schemes over a field $k$.
We will consider contravariant functors from $\cS ch_k$
to $\cA b$, the category of abelian groups; we will
refer to such a functor as a presheaf on $\cS ch_k$. We
are interested in representing such functors, when possible, by $k$-group
schemes of finite type whose identity component is a {\it semi-abelian}
$k$-scheme, \ie an extension of an abelian scheme by a $k$-torus.

 We recall that according to Grothendieck \cite{GRO} and Murre \cite{MUR}
a presheaf  $F$ is representable by a $k$-group scheme, locally of
finite type, if and only if a certain list of 7 axioms is satisfied. This
implies the following necessary and sufficient conditions for
representability by a group scheme whose identity component is
semi-abelian, where P7${}'$ is a modification of \cite[P7]{MUR}.
\begin{description}
\item[P1] $F$ is strictly pro-representable and the local components
at rational points are noetherian;
\item[P2] if $A = \liminv{n} A/\wp^{n+1}$ is a local
$k$-algebra which is complete and separated w.r.t. the
$\wp$-adique topology, then $F(A)\cong\liminv{n}F(A/\wp^{n+1})$;
\item[P3] if $A = \limdir{\alpha} A_{\alpha}$ as $k$-algebras,
then $F(A)\cong\limdir{\alpha}F(A_{\alpha})$;
\item[P4-P5] $F$ is a $fpqc$-sheaf;
\item[P6] if $T\in\cS ch_k$ and $\xi$ is a $T$-point of $F$, then
$N(\xi)\df \{ f:T'\to T/F(f)(\xi)=0\}$ yields a closed subscheme of $T$;
\item[P7${}'$] if $\xi$ is a $V$-point of $F$, for $V=C-S$ a
Zariski open of a projective non-singular $\bar k$-curve $C$, then the
induced map of $\bar k$-points has the module ${\frak m}=S$, in the
sense of Serre \cite{SEG}.
\end{description}
Here, in P7${}'$, the condition on $\frak m$ is that if
$S=\{P_1,\ldots,P_r\}$, then the non-negative integers $n_1,\ldots,n_r$
involved in the definition of a modulus (see \cite[page 10]{SEG}) are all
taken to be 1.

Let $F$, $F'$ and $F''$ be presheaves on $\cS ch_k$ forming an extension
$$ 0\to F'\to F\to F''\to 0,$$
\ie such that $0\to F'(T)\to F(T)\to F''(T)\to 0$ is
an exact sequence of abelian groups for any $T\in\cS ch_k$, which is
natural in $T$. Denote by
$$ 0\to \tilde F'\to \tilde F\to \tilde F''\to 0$$
the exact sequence of associated sheaves for the
$fpqc$-topology (\ie the faithfully flat and quasi compact
topology). If we let
\begin{equation}\label{quot}
Q(T)\df \coker (\tilde F(T)\to \tilde F''(T))
\end{equation}
for $T\in\cS ch_k$, then $Q$ is a functor on $\cS ch_k$.
\begin{lemma} \label{extrep} Assume that the $fpqc$-sheaves $\tilde
F'$ and $\tilde F''$ satisfy the axioms {\rm P1}--{\rm P3} and {\rm P6}.
Further assume that {\it i)}\, $Q$ vanishes on the subcategory of
of artinian algebras, {\it ii)}\, $Q$ vanishes on the subcategory of
local complete algebras, and {\it iii)}\, $\limdir{\alpha}Q(A_{\alpha})$
injects into $Q(\limdir{\alpha}A_{\alpha})$.

Then the $fpqc$-sheaf $\tilde F$ satisfies the axioms {\rm P1}--{\rm P3}
and {\rm P6}.
\end{lemma}
\begin{proof} The following sequence \begin{equation}\label{qseq}
 0\to \tilde F'(T)\to \tilde F(T)\to \tilde F''(T)\to Q(T)\to 0
\end{equation}
is exact and natural in $T\in\cS ch_k$. To show pro-representability
we use Grothendieck's criterion in \cite[195-5/9]{GRO}, saying that
$\tilde F$ needs to be left exact on the subcategory of artinian
algebras. From {\it i)}\, and (\ref{qseq}) one can see that $\tilde F$ is
pro-representable. Then there is a topological algebra $\cO$ such that
$\tilde F (A) \cong \Hom_c(\cO, A)$.

We show that the local components are noetherian by using Grothendieck's
criterion \cite[Prop.5.1, 195-8]{GRO}. In fact, the local component at a point
$\xi$ is pro-represented by the localization $\cO_{\wp_{\xi}}$, where
$\wp_{\xi}=\ker(\cO\by{\xi} k)$, and in order to show that
$\cO_{\wp_{\xi}}$ is noetherian it will suffices to show that
$(\wp_{\xi}/\wp_{\xi}^2)^{\vee}$ is finite dimensional.
For any $k$-scheme $T$ the $k$-point $\xi$ of $\tilde F$ induces an element
$\xi_T\in \tilde F(T)$ by pulling back along the structural morphism; we then
get an automorphism $(+\xi)_*:\tilde F\by{\cong}\tilde F$ by adding $\xi_T$
in $F(T)$. Thus we can assume $\xi=0$, therefore we have
$$(\wp_{\xi}/\wp_{\xi}^2)^{\vee}=\ker (\tilde F(k[t]/(t^2))\to \tilde F(k)).$$
By (\ref{qseq}) we conclude that P1 is satisfied.

Axioms P2--P3 follow from a diagram chase, because of (\ref{qseq}) and the
assumptions {\it ii)--iii)}. To show P6 we proceed as follows.
Let $\xi : T\to \tilde F$ be a point, \ie $\xi\in \tilde F (T)$. We have
to show that $N(\xi)$ is a closed subscheme of $T$. Let $\xi''\in
\tilde F''(T)$ be the induced point of $\tilde F''$ and let $i:
N(\xi'')\into T$ be the closed embedding. Then $\tilde F(i)(\xi)$
actually belongs to $\tilde F' (N(\xi''))$ since it yields zero
in $\tilde F'' (N(\xi''))$. Then $N (\tilde F(i)(\xi))$ is a
closed subscheme of $N(\xi'')$ hence of $T$. We can see that
$N(\xi)=N (\tilde F(i)(\xi))$. In fact, if $\alpha: T'\to N (\xi'')$
is a point such that $\tilde F (\alpha)\tilde F(i)(\xi)=0$ then
$i\alpha : T'\to T$ belongs to $N(\xi)$. Conversely, if $\alpha :T'\to T$
is such that $\tilde F (\alpha)(\xi)=0$ then $\alpha$ belongs to
$N (\xi'')$ as well which means that $\alpha = i\beta$ where
$\beta : T'\to N (\xi'')$ whence $\beta\in N (\tilde F(i)(\xi))$.
 \end{proof}

If $F$ is any functor from $\cS ch_k^{op}$ to abelian groups, and if $x\in
X$ is a point, let $x^*:F(X)\to F(k(x))$ denote the homomorphism induced
by the inclusion morphism $\Spec k(x)\to X$ determined by $x$.

\begin{lemma}\label{modulus} Let $k$ be an algebraically closed field,
and let $F:\cS ch_k^{op}\to \cA b$ be a functor, satisfying:
\begin{enumerate}
\item[(a)] the natural map $F(k)\to
F(\P^1-\{1\})$ (induced by the structure morphism) is
surjective, and
\item[(b)] if $V$ is a non-singular
quasi-projective $k$-curve, the image of the natural map
$F(S^nV)\to F(V^n)$ is the subgroup of invariants under
the natural action of the permutation group.
\end{enumerate}
Then for any non-singular projective $k$-curve, any finite subset
$S\subset C(k)$, and any rational function $f\in k(C)$ which is regular on
$S$ and has $f(x)=1$ for all $x\in S$, we have that
\[\sum_{x\in C(k)}{\rm ord}_x(f)x^*:F(C-S)\to F(k)\]
is the zero map. Thus $F$ satisfies condition {\rm P7${}'$} above.
\end{lemma}
\begin{proof} Let $C$ be a non-singular projective $k$-curve, $S\subset
C(k)$ a finite set of closed points. For any divisor $\delta=\sum_in_ix_i$
on $C-S$, let $\delta^*:F(C-S)\to F(k)$ be the map given by
$\delta^*=\sum_in_ix_i^*$. Clearly  $\delta\mapsto\delta^*$
is a homomorphism from divisors on $C_S$ to $\Hom(F(C-S),F(k))$. If
$\delta$ is an effective divisor of degree $n$, then $\delta$ determines a
point $[\delta]\in S^n(C-S)$ in an obvious way, and hence a homomorphism
\[[\delta]^*:F(S^n(C-S))\to F(k).\]
There is a homomorphism
\[\left(\sum_{i=1}^n\pi_i^*\right):F(C-S)\to F((C-S )^n),\]
where $\pi_i:(C-S)^n\to C-S$ is the $i^{\rm th}$
projection. Clearly the image is contained in the subgroup
of invariants for the action of the permutation group
$S_n$. Hence there exists a map of sets (not necessarily
unique, or even a homomorphism)
\[\psi:F(C-S)\to F(S^n(C-S)),\;\;\xi\mapsto S^n(\xi)\]
such that $S^n(\xi)$ is a pre-image in $F(S^n(C-S)$ of
$\left(\sum_{i=1}^n\pi_i^*\right)(\xi)$, for any $\xi\in
F(C-S)$, \ie the diagram
\[\begin{TriCDV}
{F(C-S)}{\>\sum_{i=1}^n\pi_i^*>>}{F((C-S)^n)}
{\SE \psi E E}{\NE E E}
{F(S^n(C-S))}
\end{TriCDV}\]
commutes.

We claim that for any effective divisor $\delta=\sum_jn_jx_j$ of degree
$n$ on $C-S$, we have
\begin{equation}\label{eq1}
[\delta]^*\circ\psi=\delta^*\in \Hom(F(C-S),F(k)).
\end{equation}
Indeed, let $\tilde{\delta}=(x_1,\ldots,x_1,x_2,\ldots,x_2,\cdots)\in
(C-S)^n$, where $x_j$ is repeated $n_j$ times as a coordinate. Then
$\tilde{\delta}$ is a preimage in $(C-S)^n$ of $[\delta]\in S^n((C-S)^n)$.
Hence
\[\tilde{\delta}^*\circ\left(\sum_{i=1}^n\pi_i^*\right)=
[\delta]^*\circ\psi.\]
On the other hand, from the definitions, it is clear that
\[\tilde{\delta}^*\circ\left(\sum_{i=1}^n\pi_i^*\right)=\sum_jn_jx_j^*
=\delta^*.\]

Now suppose $f\in k(C)$ such that $f\mid_S=1$. Let $T=f^{-1}(1)$, and
consider $f$ as a morphism $f:C-T\to\P^1-\{1\}$. There is an induced
morphism
\[\tilde{f}:\P^1-\{1\}\to S^n(C-T)\into S^n(C-S),\]
where $n=\deg f$. The map $\tilde{f}$ has the property that if
$\delta_t=f^{-1}(t)$ as a divisor, then $\tilde{f}(t)=[\delta_t]$.

Let $\delta_0=(f)_0$, $\delta_{\infty}=(f)_{\infty}$ be
the divisors of poles and zeroes of $f$. Then the lemma asserts that
\[\delta_0^*=\delta_{\infty}^*:F(C-S)\to F(k).\]
To prove this, by (\ref{eq1}), it suffices to show that
\[[\delta_0]^*=[\delta_{\infty}]^*:F(S^n(C-S))\to F(k).\]

Since $\tilde{f}(t)=[\delta_t]$, it follows that
\[ [\delta_t]^*=t^*\circ \tilde{f}^*\mbox{\ \ }\forall\;t\in\P^1-\{1\},\]
and so we are reduced to proving that $0^*={\infty}^*:F(\P^1-\{1\})\to
F(k)$.

If $\pi:\P^1-\{1\}\to\Spec k$ is the structure morphism, then we are given
that $$\pi^*:F(k)\onto F(\P^1-\{1\}),$$ while clearly $0^*\circ\pi^*$ and
${\infty}^*\circ\pi^*$ both equal the identity on $F(k)$. Hence
$0^*={\infty}^*$ as desired.
\end{proof}

Now we can easily show that our ``relative'' and ``simplicial'' Picard
functors are representable by checking the list of axioms
P1--P7${}'$. Because of Lemma~\ref{extrep} and Lemma~\ref{modulus}
representability will follows from the representability
theorems for the classical Picard functor: we will sketch the arguments
below.

We remark that, in the particular case when we have an extension of sheaves
as above and we moreover assume that $\tilde F'$ is affine, one can then
also deduce representability of $\tilde F$ by descent, as in
Proposition~17.4 of \cite{OG}.

\subsection{Representability of the relative Picard functor}

In order to show representability of $fpqc$-sheaves one can assume that
the base field $k$ is algebraically closed (see \cite[Lemma~I.8.9]{MUR}).
The $fpqc$-sheaf associated to the relative Picard functor
$T\leadsto\Pic (\bar X\times_k T, Y\times_k T)$ in Lemma~\ref{grpic}
 will be denoted by $\Pic_{(\bar X, Y)/k}$.  The exact
sequence (\ref{longpic}) yields the following short exact sequence of
$fpqc$-sheaves $$0\to T(\bar X, Y)\to \Pic_{(\bar X, Y)/k}\to \ker
(\Pic_{\bar X/k}\to \Pic_{Y/k})\to 0$$

Since $\bar X$ is non-singular and complete, we see that to prove
representability of $\Pic_{(\bar X,Y)/k}$, we reduce immediately
to the case when $\bar X$ is connected, hence irreducible. If
$Y=\emptyset$, then \cite{MUR} yields the desired representability. If
$\bar{X}$ is irreducible and $Y\neq \emptyset$ (as we may now assume),
then pairs $(\ccL,\alpha)$ consisting of line bundles on $\bar{X}$,
trivialized along $Y$, do not admit non-trivial automorphisms. Therefore,
the functor which takes a $k$-scheme $T$ to
$\Pic (\bar X\times_k T, Y\times_k T)$ is already a
sheaf with respect to the Zariski topology and, by descent theory (see
\cite[\S 2.1]{RAS}, \cite[\S 8.1]{BLR}), even with respect to the
$fpqc$-topology.

We now apply our Lemma~\ref{extrep} and Lemma~\ref{modulus} to the
functors $F= \tilde F =\Pic_{(\bar X, Y)/k}$, $F' = \tilde F' = T(\bar X,
Y)$ and
$$F'' = \ker (\Pic_{\bar X/k}\to \Pic_{Y/k}) = \ker (\Pic_{\bar X/k}\to
\oplus\Pic_{Y_i/k})$$
where $Y_j$ for $j= 1, ... ,r$ are the connected (possibly reducible)
components of $Y$. We have that $\bar X$ integral and
$(\pi_{\bar X})_*(\cO_{\bar X}) =(\pi_{Y_j})_*(\cO_{Y_j})=k$. We then have
the following commutative diagram  with exact rows and columns
$$\begin{array}{ccccc}
 & &0& &0\\
 & &\downarrow & &\downarrow\\
0&\to& \Pic (T)&\longby{\rm diag} &\oplus_{i=1}^{r} \Pic (T)\\
\downarrow & &\downarrow & &\downarrow\\
F''(T)&\into &\Pic(\bar X\times_k T) &\to &
\oplus_{j=1}^{r} \Pic(Y_j\times_k T)\\
\downarrow & &\downarrow & &\downarrow\\
\tilde F''(T)&\into &\frac{\Pic(\bar X\times_k T)}{\Pic (T)} &\to &
\oplus_{j=1}^{r} \frac{\Pic(Y_j\times_k T)}{\Pic (T)}\\
\downarrow & &\downarrow & &\\
(\Pic (T))^{\oplus r-1}& & 0 & &\\
\downarrow & & & &\\
0 & & & &
\end{array}$$
Thus we have that the functor $Q (T)$ in (\ref{quot}) is canonically
isomorphic to $(\Pic (T))^{\oplus r-1}$. It is then easy to see that $Q$
satisfies the  hypotheses {\it i) --- iii)}\, stated in the
Lemma~\ref{extrep}.

Let $\coprod_i Y_i \to Y$ be the normalization of the normal crossing
divisor $Y$, where $Y_i$ are now the irreducible components of $Y$.
Consider the following exact sequence
$$0\to \ker (\Pic_{\bar X/k}\to \Pic_{Y/k})\to \ker (\Pic_{\bar X/k}\to
\oplus_i\Pic_{Y_i/k}) \by{\rho} \ker (\Pic_{Y/k}\to \oplus
\Pic_{Y_i/k}).$$
In the Lemma~\ref{grpic} we have shown that the map $\rho$ above vanishes
on the connected components of the identity yielding a description of the
semi-abelian scheme $\Pic^0_{(\bar X, Y)/k}$.

\subsection{Representability of the simplicial Picard functor}
Let $X_{\d}$ be a simplicial scheme. We first construct an explicit functorial
isomorphism
$$\bPic (X_{\d})\longby{\cong} \HH^1(X_{\d},\cO^*_{X_{\d}})$$
as claimed in Proposition~\ref{simpline}. We clearly can bijectively
associate to
(the isomorphism class of) a simplicial line bundle $\ccL_{\d}$ on
$X_{\d}$ (\ie to an invertible $\cO_{X_{\d}}$-module)
a pair $(\ccL,\alpha)$ consisting of a line bundle $\ccL$ on $X_0$ and an
isomorphism $\alpha :(d_0)^*(\ccL) \by{\cong}(d_1)^*(\ccL)$ on $X_1$,
satisfying the cocycle condition (as in the Section~4.1).

Assume given:
\begin{description}
\item[{\it 1)}\,] an element $\xi\in\mbox{\v {$H$}}^1(\cU,\cO_{X_0}^*)$,
for an open covering $\cU =\{U_i\}_{i\in I}$ of $X_0$, corresponding to a
line bundle $\ccL\in \Pic (X_0)$, together with trivializations
$s_i:\cO_{U_i}\by{\cong} \ccL\mid_{U_i}$; then
$\xi =\{f_{ij}\in \cO^*_{X_0}(U_i\cap U_j)\}$ with
$s_i\mid_{U_i\cap U_j} = f_{ij} s_j\mid_{U_i\cap U_j}$
\item[{\it 2)}\,] an isomorphism $\alpha :(d_0)^*(\ccL)
\by{\cong}(d_1)^*(\ccL)$ satisfying the cocycle condition.
\end{description}
Let $$V_{ij}\df d^{-1}_1(U_i)\cap d^{-1}_0(U_j).$$
Then $\{V_{ij}\}_{(i,j)\in I\times I}$ is an open covering of $X_1$.
Moreover, on
$V_{ij}$ we have trivializations $d^{*}_1(s_i)$ of $(d_1)^*(\ccL)$, and
$d^{*}_0(s_j)$ of $(d_0)^*(\ccL)$, respectively. Therefore, $\alpha$ is
uniquely
determined by $\alpha_{ij}\in \cO^*_{X_1}(V_{ij})$, satisfying
$$d^{*}_1(s_i) = \alpha_{ij}d^{*}_0(s_j)$$
on $V_{ij}$. The $\alpha_{ij}$ have to satisfy a compatibility condition: on
$V_{ij}\cap V_{kl}$, have $d^{*}_1(s_i) = d^{*}_1(f_{ik})d^{*}_1(s_k)$,
and
$d^{*}_0(s_j) = d^{*}_0(f_{jl})d^{*}_0(s_l)$, thus
$d^{*}_1(s_i) = d^{*}_1(f_{ik})d^{*}_1(s_k) =
d^{*}_1(f_{ik})\alpha_{kl}d^{*}_0(s_l)$,
but $d^{*}_1(s_i) = \alpha_{ij}d^{*}_0(s_j) = \alpha_{ij}
d^{*}_0f_{jl}d^{*}_0(s_l)$
as well, therefore $d^{*}_0(s_l)$ cancels and we obtain
\begin{equation}\label{alcomp}
d^{*}_1f_{ik}\alpha_{kl} =\alpha_{ij} d^{*}_0f_{jl}
\end{equation}
on $V_{ij}\cap V_{kl}$.

Let $\cK_i$ be the canonical (Godement) flasque sheaf of discontinuous
sections of
$\cO^*_{X_i}$ and let $\cQ_i \df \frac{\cK_i}{\cO^*_{X_i}}$ denote the quotient
sheaf. We have an exact sequence
$$0\to \Gamma (X_i,\cO^*_{X_i})\to  \Gamma (X_i,\cK_i)\to
\Gamma (X_i,\cQ_i)\to H^1(X_i,\cO^*_{X_i})\to 0$$
Choose a function $\varphi : X_0\to I$ such that
$x\in U_{\varphi (x)}$ for any $x\in X$. This determines well defined sections
$t_i\in \cK_0(U_i)$, $t_i(x) = f_{i\varphi (x)}\in \cO^*_{X_0,x}$. For any
$x\in U_i\cap U_j$, we have $f_{ij} = f_{i\varphi (x)} f_{j\varphi (x)}^{-1}\in
\cO^*_{X_0,x}$, therefore: $t_i = f_{ij} t_j$ on $U_i\cap U_j$ and the images
of $t_i$ in $\cQ_0(U_i)$ patch together to give a global section $t\in
\Gamma (X_0,\cQ_0)$. By construction,
$$[t]\in \frac{\Gamma (X_0,\cQ_0)}{\im \Gamma (X_0,\cK_0)} =
H^1(X_0,\cO^*_{X_0})$$
is the class of the given line bundle $\ccL$ on $X_0$.

Next, we have a natural element $\beta\in \Gamma (X_1,\cK_1)$, given by
$$\beta (x) = \alpha_{\varphi (d_1(x))\varphi (d_0(x))}\in \cO^*_{X_0,x};$$
note that $x\in V_{\varphi (d_1(x))\varphi (d_0(x))}$ by the definition of
$\varphi$.
We claim that (writing the group operation in $\Gamma (X_1,\cQ_1)$
multiplicatively):
\begin{equation}\label{claim}
[\beta] = \frac{d^{*}_0(t)}{d^{*}_1(t)} \in \frac{\Gamma
(X_1,\cK_1)}{\Gamma (X_1,
\cO^*_{X_1})} \subset \Gamma (X_1,\cQ_1).
\end{equation}
In fact, on $V_{ij}$ we have that $d^{*}_1(t)$ is the image of $d^{*}_1(t_{i})$
where $d^{*}_1(t_{i})(x) = d^{*}_1(f_{i\varphi (d_1(x))})\in \cO^*_{X_1,x}$;
similarly, $d^{*}_0(t)$ is the image of $d^{*}_0(t_{j})$
where $d^{*}_0(t_{j})(x) = d^{*}_0(f_{j\varphi (d_0(x))})\in
\cO^*_{X_1,x}$. From the definition of $\cQ_1$ as a quotient sheaf, the
claimed formula  (\ref{claim}) will be proved if: for any $i, j\in I$ and
$x\in V_{ij}$, $d^{*}_1(t_i)(d^{*}_0(t_j))^{-1}\beta (x)$ defines a
section  in $\cO^*_{X_1}(V_{ij})$. From the identity (\ref{alcomp}) we
have  $$d^{*}_1(t_i)(d^{*}_0(t_j))^{-1}\beta (x)= d^{*}_1(f_{i\varphi
(d_1(x))})( d^{*}_0(f_{j\varphi (d_0(x))}))^{-1}\alpha_{\varphi
(d_1(x))\varphi (d_0(x))} =\alpha_{ij} \in  \cO^*_{X_1,x}$$

Thus, given a simplicial line bundle, therefore data as in {\it 1)}\, and
{\it 2)}\, above, together with a choice of $\varphi : X\to I$ we get an
element of $\HH^1 (X_{\d}, \cO_{X_{\d}}^*)$ computed by means of the
canonical Godement resolution of the simplicial sheaf $\cO_{X_{\d}}^*$. It
is now easy to verify that this construction is independent of the
additional choices made (the local trivializations $s_i$ and the map
$\varphi$), and defines a homomorphism
$$\bPic (X_{\d})\to \HH^1(X_{\d},\cO^*_{X_{\d}}).$$

Conversely, we see that $\HH^1(X_{\d},\cO^*_{X_{\d}})$ is identified with
the $H^1$ of the following complex
$$\Gamma (X_0,\cK_0)\to \Gamma (X_1,\cK_1)\oplus \Gamma (X_0,\cQ_0)
\to \Gamma (X_2,\cK_2)\oplus \Gamma (X_1,\cQ_1).$$
Given a cycle $(\beta,t)\in \Gamma (X_1,\cK_1)\oplus \Gamma (X_0,\cQ_0)$, we
can choose an open cover $\{U_i\}$ of $X_0$ and pre-images $t_i\in \cK_0(U_i)$
of $t$, and we will then obtain $f_{ij}\in \cO^*_{X_0}(U_i\cap U_j)$ satisfying
$t_i=f_{ij}t_j$ in $\cK_0(U_i\cap U_j)$. Now one immediately verifies that
the $f_{ij}$ define an invertible sheaf $\ccL$ on $X_0$, and (reversing the
earlier arguments) $\beta$ determines an isomorphism
$\alpha:d_0^*\ccL\by{\cong} d_1^*\ccL$. Since $\beta$ maps to 0 in
$\Gamma(X_2,\cK_2)$, we  deduce that $\alpha$ satisfies the cocycle
condition.\\

We now come to the proof of the representability of the simplicial Picard
functor. Let $X_{\d}$ be a smooth and proper simplicial $k$-scheme.
To prove representability of the simplicial Picard functor, we again
reduce to the case when $k$ is algebraically closed, using
\cite[Lemma~I.8.9]{MUR}. Then we may further reduce to the case when
$(\pi_{X_{\d}})_*(\cO_{X_{\d}})=k$. Then we have
\begin{equation}\label{eq}
\bPic_{X_{\d}/k}(T)=\frac{\Pic(X_{\d}\times_kT)}{\Pic T}
\end{equation}
since we can choose a base point in $X_{\d}$.
We let $\Z^{X_a}$ denote the free abelian group on the connected
components of $X_a$, $a= 0, 1,\ldots$, and let $\pi_a:X_a\to k$ be
the structural morphism.

We then set (\cf Section~\ref{cohompic})
$$\begin{array}{ll} K\df \ker (\Z^{X_0}\to\Z^{X_1})\\ C\df \coker
(\Z^{X_0}\to\Z^{X_1})&\\ F' \df T(\bar X_{\d}) \df \frac{\ker
((\pi_1)_*\G_{m,X_1}\to (\pi_2)_*\G_{m,X_2})}{\im ((\pi_0)_*\G_{m,X_0}\to
(\pi_1)_*\G_{m,X_1})}&\\ F(T) \df \bPic (X_{\d}\times_k T) &\tilde F =
\bPic_{X_{\d}/k}\\ F''(T) \df \ker (\Pic(X_0\times_k T)\to \Pic(X_1\times_k
T)) &\tilde F''=\ker \Pic_{X_0/k}\to\Pic_{X_1/k}\\
G\df \frac{\ker ((\pi_2)_*\G_{m,X_2}\to (\pi_3)_*\G_{m,X_3})}{\im
((\pi_1)_*\G_{m,X_1}\to (\pi_2)_*\G_{m,X_2})}&\tilde{G}=\mbox{associated
$fpqc$ sheaf.}
\end{array}$$
We then have an exact sequence of pre-sheaves
$$0\to F'\to F\to F''\to G$$
and, for each $T\in\cS ch_k$, a commutative diagram of complexes
$$\begin{array}{cccccccc}
0\to & F'(T)& \to  & F(T) & \to & F''(T)& \to & G(T)\\
 &\downarrow & &\downarrow & & \downarrow & & \downarrow\\
0\to & \tilde F'(T) &\to &\tilde F(T) & \to & \tilde F''(T) & \to & \tilde
G(T)\end{array}$$
with exact top row. By (\ref{eq}), we have that $F(T) \to \tilde F(T)$ is
surjective with kernel $\Pic (T)$. Moreover we have the following
commutative diagram with exact rows and columns:
$$\begin{array}{ccccc}
 0&&0& &0\\
 \downarrow& &\downarrow & &\downarrow\\
 \Pic (T)\otimes K&\to&\oplus_{i=1}^{r} \Pic (T)&\to &\oplus_{j=1}^{s}\Pic
(T)\\
\downarrow & &\downarrow & &\downarrow\\
F''(T)&\into &\oplus_{i=1}^{r}\Pic(X_0^i\times_k T) &\to &
\oplus_{j=1}^{s}\Pic(X_1^j\times_k T)\\
\downarrow & &\downarrow & &\downarrow\\
\tilde F''(T)&\into &\oplus_{i=1}^{r}\frac{\Pic(X_0^i\times_k T)}{\Pic
(T)}&\to &
\oplus_{j=1}^{s}\frac{\Pic(X_1^j\times_k T)}{\Pic (T)}\\
\downarrow & &\downarrow & &\downarrow \\
\Pic (T)\otimes C& & 0 & &0\\
\downarrow & & & &\\
0 & & & &
\end{array}$$
Now we can see that $F'(T)\cong \tilde F'(T)$ in fact: $F'$ is of
$\G_m$-type
whence the map $F'(T)\to\tilde F'(T)$ is surjective with finite kernel but
the
finite kernel is actually zero being isomorphic to the kernel of the
injective map
$\Pic (T)\into  \Pic (T)\otimes K$.
Now we let $G'$ denote the image of $F$ in $F''$. The associated
$fpqc$-sheaf
$\tilde G'$ is representable, in fact: $\tilde G'$ is the kernel of the
homomorphism
of group schemes $\tilde F'' \to \tilde G$.
We then have the following commutative diagram with exact rows:
$$\begin{array}{cccccccc}
0\to & F'(T)& \to  & F(T) & \to & G'(T)& \to 0&\\
 &\downarrow & &\downarrow & & \downarrow & &\\
0\to & \tilde F'(T) &\to &\tilde F(T) & \to & \tilde G'(T) & \to & Q(T)\to
0
\end{array}$$
where, by definition,  $Q(T)$ is the cokernel of $\tilde F(T) \to\tilde G'(T)$
and we can apply our Lemma~\ref{extrep} and Lemma~\ref{modulus}.
As $G$ is of $\G_m$-type then the map $G(T)\onto\tilde G(T)$ is surjective
with finite kernel $G_0(T)$. Moreover we have that $\Pic(T)$ and $Q(T)$ are
respectively  the kernel and the cokernel of $G'(T)\to \tilde G'(T)$.
Considering $\tilde G'$ as the kernel of the homomorphism of group schemes
$\tilde F'' \to \tilde G$ we can see that there is a functorial exact
sequence  \begin{equation} \label{qsimp}
0\to  {\rm finite\ group}\to Q(T)\to \Pic (T)\otimes C \end{equation}
where the finite group is a subgroup of $G_0(T)$ whence it is zero whenever
$H^0(T,\cO^*_T)$ is divisible \eg if $T$ is an artinian algebra or a
strictly  Hensel local ring.

Therefore we can easily check the vanishing conditions
of Lemma~\ref{extrep}: {\it i) --- ii)}\, follow from (\ref{qsimp}), and
{\it iii)}\, follows by a diagram chase using (\ref{qsimp}) since  $\Pic$
and $H^0(T,\cO^*_T)$ commute with the relevant direct limits.

\vspace*{2cm}
\noindent{Dipartimento di Matematica, Universit\`a di Genova, Via
Dodecaneso, 35,\\16146 -- Genova, Italia\\[2mm]
School of Mathematics, Tata Institute of Fundamental Research,
Homi Bhabha Road,\\
Mumbai-400005, India.}\\[2mm]
e-mail: {\it barbieri}@{\it dima.unige.it} \ and\  {\it
srinivas}@{\it math.tifr.res.in}

\end{document}